\documentclass[fleqn]{amsart}

\usepackage[all]{xy}
\usepackage{latexsym}
\usepackage{amsfonts}
\usepackage{amssymb}
\usepackage{stmaryrd}
\usepackage{a4wide}
\usepackage[utf8]{inputenc}
\usepackage{dsfont}

\newtheorem{thmalpha}{Theorem}

\swapnumbers
\newtheorem{thm}[subsubsection]{Theorem}

\newtheorem{prop}[subsubsection]{Proposition}
\newtheorem{lem}[subsubsection]{Lemma}
\newtheorem{cor}[subsubsection]{Corollary}

\theoremstyle{definition}
\newtheorem*{rem}{\sc Remark}

\newtheorem*{pf}{\sc Proof}

\newtheorem*{ex}{\sc Examples}
\newtheorem*{ex1}{\sc Example}

\newcommand{\B}{\mathcal{B}}
\newcommand{\C}{\mathcal{C}}
\newcommand{\D}{\mathcal{D}}
\newcommand{\E}{\mathcal{E}}
\newcommand{\F}{\mathcal{F}}
\newcommand{\LbP}{\mathbb{L}_{\mathcal{P}}}
\newcommand{\R}{\mathcal{R}}

\newcommand{\Po}{\mathcal{P}}
\newcommand{\So}{\mathbb{S}}
\newcommand{\K}{\mathbb{K}}
\newcommand{\Kc}{\mathbb{O}_{\mathcal{P}}}

\newcommand{\Lb}{\mathbb{L}}
\newcommand{\M}{\mathcal{M}}

\newcommand{\RT}{\mathcal{RT}}
\newcommand{\cqfd}{\ \hfill \square}

\newcommand{\Ass}{\mathcal{A}ss}
\newcommand{\Lie}{\mathcal{L}ie}
\newcommand{\Com}{\mathcal{C}om}
\newcommand{\Dias}{\mathcal{D}ias}
\newcommand{\Leib}{\mathcal{L}eib}
\newcommand{\Zinb}{\mathcal{Z}inb}
\newcommand{\Poiss}{\mathcal{P}oiss}
\newcommand{\Prelie}{\mathcal{P}relie}
\newcommand{\Perm}{\mathcal{P}erm}

\newcommand{\ash}{\textrm{!`}}
\newcommand{\Poa}{\Po^{\ash}}
\newcommand{\epi}{\twoheadrightarrow}
\newcommand{\mono}{\rightarrowtail}
\newcommand{\iso}{\xrightarrow{\simeq}}
\newcommand{\qiso}{\xrightarrow{\sim}}
\newcommand{\Hom}{\mathrm{Hom}}

\title{André-Quillen cohomology of algebras over an operad}

\author{Joan Mill\`es}

\thanks{Joan Millès, Laboratoire J. A. Dieudonné, Université de Nice Sophia-Antipolis, Parc Valrose, 06108 Nice Cedex 02, France\\
E-mail : \texttt{joan.milles@math.unice.fr}\\
\textsc{URL :} \texttt{http://math.unice.fr/$\sim$jmilles}}

\begin{document}

\maketitle

\begin{abstract}
We study the André-Quillen cohomology with coefficients of an algebra over an operad. Using resolutions of algebras coming from the Koszul duality theory, we make this cohomology theory explicit and we give a Lie theoretic interpretation. For which operads is the associated André-Quillen cohomology equal to an Ext-functor ? We give several criteria, based on the cotangent complex, to characterize this property. We apply it to homotopy algebras, which gives a new homotopy stable property for algebras over cofibrant operads.\\

\end{abstract}

\section*{Introduction}

Hochschild \cite{Hochschild} introduced a chain complex which defines a cohomology theory for associative algebras. In $1948$, Chevalley and Eilenberg gave a definition of a cohomology theory for Lie algebras. Both cohomology theories can be written as classical derived functors (Ext-functors). Later, Quillen \cite{Quillen} defined a cohomology theory associated to commutative algebras with the use of model category structures. André gave similar definitions only with simplicial methods \cite{Andre}. This cohomology theory is not equal to an Ext-functor over the enveloping algebra in general.\\

Using conceptual model category arguments, we recall the definition of the \emph{André-Quillen cohomology (for algebras over an operad)}, in the differential graded setting, from Hinich \cite{Hinich} and Goerss and Hopkins \cite{GoerssHopkins}. Because we work in the differential graded setting, we use known functorial resolutions of algebras to make chain complexes which compute André-Quillen cohomology explicit. The first idea of this paper is to use Koszul duality theory of operads to provide such functorial resolutions. We can also use the simplicial bar construction, which proves that cotriple cohomology is equal to André-Quillen cohomology. The André-Quillen cohomology is represented by an object, called the \emph{cotangent complex} which therefore plays a crucial role in this theory. The notion of \emph{twisting morphism}, also called twisting cochain, coming from algebraic topology, has been extended to (co)operads and to (co)algebras over a (co)operad by Getzler and Jones \cite{GetzlerJones}. We make the differential on the cotangent complex explicit using these two notions of twisting morphisms all together. When the category of algebras is modeled by a binary Koszul operad, we give a Lie theoretic interpretation of the previous construction. In the review of \cite{Frabetti}, Pirashvili asked the question of characterizing operads such that the associated André-Quillen cohomology of algebras is an Ext-functor. This paper provides a criterion to answer that question.\\

When the operad is Koszul, we describe the cotangent complex and the André-Quillen cohomology for the algebras over this operad using its Koszul complex. We recover the classical cohomology theories, with their underlying chain complexes, like André-Quillen cohomology for commutative algebras, Hochschild cohomology for associative algebras and Chevalley-Eilenberg cohomology for Lie algebras. We also recover cohomology theories which were defined recently like cohomology for Poisson algebras \cite{Fresse2}, cohomology for Leibniz algebras \cite{LodayPirashvili}, cohomology for pre-Lie algebras \cite{Dzhumadil}, cohomology for diassociative algebras \cite{Frabetti} and cohomology for Zinbiel algebras \cite{Balavoine}. More generally, Balavoine introduced a chain complex when the operad is binary and quadratic \cite{Balavoine}. We show that this chain complex defines André-Quillen cohomology when the operad is Koszul. We make the new example of Perm algebras explicit. For any operad $\Po$, we can define a relax version up to homotopy of the notion of $\Po$-algebra as follows: we call homotopy $\Po$-algebra any algebra over a cofibrant replacement of $\Po$ (cf. \cite{BoardmanVogt}). Using the operadic cobar construction, we make the cotangent complex and the cohomology theories for homotopy algebras explicit. For instance, we recover the case of homotopy associative algebras \cite{Markl} and the case of homotopy Lie algebras \cite{HinichSchechtman}.\\

For any algebra $A$, we prove that its André-Quillen cohomology is an additive derived functor, an Ext-functor, over its enveloping algebra if and only if its cotangent complex is quasi-isomorphic to its module of Kähler differential forms $\Omega_{\Po}(A)$.  We define a \emph{functorial cotangent complex} and a \emph{functorial module of K\"ahler differential forms} which depend only on the operad and we reduce the study of the quasi-isomorphisms between the cotangent complex and the module of K\"ahler differential forms for any algebra to the study of the quasi-isomorphisms between the cotangent complex and the module of K\"ahler differential forms for any chain complex, with trivial algebra structure (when $\Po$ is an \emph{PBW-operad}, that is the $\Po$-algebras satisfy an analogue of the Poincaré-Birkhoff-Witt theorem and the $\Po$-Kähler differentials too). This allows us to give a uniform treatment for any algebra over an operad. Assuming that $\Po$ is a PBW-operad, we prove that the functorial cotangent complex is quasi-isomorphic to the functorial module of K\"ahler differential forms (we say sometimes concentrated in degree $0$ or acyclic), if and only if the André-Quillen cohomology theory for any algebra over this operad is an Ext-functor over its enveloping algebra, so this functorial cotangent complex carries the obstructions for the André-Quillen cohomology to be an Ext-functor. For instance, we prove that the functorial cotangent complex is acyclic for the operads of associative algebras and Lie algebras. In order to control the map between the functorial cotangent complex and the functorial module of K\"ahler differential forms, we look at its kernel. This defines a new chain complex whose homology groups can also be interpreted as obstructions for the André-Quillen cohomology theory to be an Ext-functor. In this way, we give a new, but more conceptual proof that the cotangent complex for commutative algebras is not always acyclic. Equivalently, it means that there exist commutative algebras such that their André-Quillen cohomology is not an Ext-functor over their enveloping algebra. With the same method, we show the same result for Perm algebras. We can summarize all these properties in the following theorem (Section 4 and 5).

\begin{thmalpha}\label{A}
Let $\Po$ be an operad. The following two properties are equivalent:
\begin{itemize}
\item[$(P_{0})$] the André-Quillen cohomology is an Ext-functor over the enveloping algebra $A\otimes^{\Po}\K$ for any $\Po$-algebra $A$;
\item[$(P_{1})$] the cotangent complex is quasi-isomorphic to the module of Kähler differential forms for any $\Po$-algebra $A$.
\end{itemize}
They moreover imply the following equivalent properties:
\begin{itemize}
\item[$(P_{2})$] the functorial cotangent complex $\LbP$ is quasi-isomorphic to the functorial module of Kähler differential forms ${\bf \Omega_{\Po}}$;
\item[$(P_{3})$] the module of obstructions $\Kc$ is acyclic.
\end{itemize}
When $\Po$ be a PBW-operad, the previous implication is an equivalence, that is
 $$
 (P_{0}) \Leftrightarrow (P_{1}) \Leftrightarrow (P_{2}) \Leftrightarrow (P_{3}).
 $$
\end{thmalpha}

In the case of homotopy algebras, we prove that the obstructions for the cohomology to be an Ext-functor vanish. Moreover, any $\Po$-algebra is also a homotopy $\Po$-algebra. Thus we can compute its André-Quillen cohomology in two different ways. We show that the two coincide. Hence we get the following theorem.

\begin{thmalpha}\label{B}
Let $A$ be a $\Po$-algebra and let $M$ be an $A$-module over the Koszul operad $\Po$. We have
$$\textrm{\emph{H}}^{\bullet}_{\Po}(A,\, M) \cong \emph{Ext}_{A\otimes^{\Po_{\infty}}\K}^{\bullet}(\Omega_{\Po_{\infty}}(A),\, M).$$
\end{thmalpha}

Therefore, even if the André-Quillen cohomology of commutative and Perm algebras cannot always be written as an Ext-functor over the enveloping algebra $A\otimes^{\Po}\K$, it is always an Ext-funtor over the enveloping algebra $A\otimes^{\Po_{\infty}}\K$.\\

The paper begins with first definitions and properties of differential graded (co)operads, (co)al\-ge\-bras, modules and free modules over an algebra (over an operad). In Section $1$, we recall the definition of the André-Quillen cohomology theory for dg algebras over a dg operad, from Hinich and Goerss-Hopkins. We introduce functorial resolutions for algebras over an operad, which allow us to make the cotangent complex and the cohomology theories explicit. Then, in Section $2$, we give a Lie interpretation of the chain complex defining the André-Quillen cohomology. Using the notion of twisting morphism on the level of (co)algebra over a (co)operad, we make the differential on the cotangent complex explicit (Theorem \ref{iso2}). Section $3$ is devoted to applications and examples. In Section $4$, we prove that the cotangent complex is quasi-isomorphic to the module of Kähler differential forms for any algebra if and only if the André-Quillen cohomology theory is an Ext-functor over the enveloping algebra for any algebra. Moreover, we study the André-Quillen cohomology theory for operads. In Section $5$, we introduce the functorial cotangent complex and the functorial module of K\"ahler differential forms and we finish to prove Theorem \ref{A}. In Section $6$, we study the André-Quillen cohomology for homotopy algebras and we prove Theorem \ref{B}.

\tableofcontents

\section*{Notation and preliminary}

We recall the classical notation for $\So$-module, composition product, (co)operad, (co)algebra over a (co)operad and module over an algebra over an operad. We refer to \cite{GinzburgKapranov} and \cite{GetzlerJones} for a complete exposition and \cite{Fresse} for a more modern treatment. We also refer to the books \cite{LodayVallette} and \cite{MarklShniderStasheff}.\\

In the whole paper, we work over a field $\K$ of characteristic $0$. In the sequel, the ground category is the category of graded modules, or \emph{g-modules}. For a morphism $f : O_{1} \rightarrow O_{2}$ between differential graded modules, the notation $\partial(f)$ stands for the derivative $d_{O_{2}} \circ f - (-1)^{|f|} f \circ d_{O_{1}}$. Here $f$ is a map of graded modules and $\partial(f) = 0$ if and only if $f$ is a map of dg-modules. Moreover, for an other morphism $g : O_{1}' \rightarrow O_{2}'$, we define a morphism $f\otimes g : O_{1}\otimes O_{1}' \rightarrow O_{2}'\otimes O_{2}'$ using the Koszul-Quillen convention: $(f\otimes g)(o_{1}\otimes o_{2}) := (-1)^{|g| |o_{1}|}f(o_{1})\otimes g(o_{2})$, where $|e|$ denotes the degree of the element $e$. We denote by $g\M od_{\K}$ the category whose objects are differential graded $\K$-modules (and not only graded $\K$-modules) and morphisms are maps of graded modules. We have to be careful with this definition because it is not usual. However, we denote as usual by $dg\mathcal{M}od_{\K}$ the category of differential graded $\K$-modules. In this paper, the modules are all differential graded, except explicitly stated.

\subsection{Differential graded $\So$-modules}\label{dgsm}
A \emph{dg $\So$-module} (or \emph{$\So$-module} for short) $M$ is a collection $\{M(n)\}_{n\geq 0}$ of dg modules over the symmetric group $\So_{n}$. A \emph{morphism of dg $\So$-modules} is a collection of equivariant morphisms of chain complexes $\{f_{n}: M(n) \rightarrow N(n)\}_{n\geq 0}$, with respect to the action of $\So_{n}$.\\
We define a monoidal product on the category of dg $\So$-modules by
$$(M \circ N)(n) := \bigoplus_{k \geq 0} M(k) \otimes_{\So_k} \left(\bigoplus_{i_1 + \cdots + i_k = n}
Ind_{\So_{i_{1}}\times \cdots \times \So_{i_{k}}}^{\So_{n}} (N(i_1) \otimes \cdots \otimes N(i_k))\right).$$
The unit for the monoidal product is $I := (0,\, \K,\, 0,\, \ldots)$.
Let $M$, $N$ and $N'$ be dg $\So$-modules. We define the right linear analog $M \circ (N,\, N')$ of the composition product by the following formula\\
\begin{eqnarray*}
&& [M \circ (N,\, N')](n) :=\\
&& \bigoplus_{k \geq 0} M(k) \otimes_{\So_k} \left(\bigoplus_{i_1 + \cdots + i_k = n} \bigoplus_{j = 1}^{k} Ind_{\So_{i_{1}}\times \cdots \times \So_{i_{k}}}^{\So_{n}} (N(i_1) \otimes \cdots \otimes \underbrace{N'(i_{j})}_{j^{th} \textrm{ position}} \otimes \cdots \otimes N(i_k))\right).
\end{eqnarray*}
Let $f : M \rightarrow M'$ and $g : N \rightarrow N'$ be morphisms of dg $\So$-modules. We denote by $\circ'$ the \emph{infinitesimal composite of morphisms}
$$f \circ' g : M \circ N \rightarrow M' \circ (N,\, N')$$
defined by
$$\sum_{j=1}^k f \otimes (id_{N} \otimes \cdots \otimes \underbrace{g}_{j^{th} \textrm{ position}} \otimes \cdots \otimes id_{N}).$$
Let $(M,\, d_{M})$ and $(N,\, d_{N})$ be two dg $\So$-modules. We define a grading on $M\circ N$ by
$$(M\circ N)_{g}(n) := \bigoplus_{\genfrac{}{}{0cm}{1}{k\geq 0}{e+g_{1}+\cdots +g_{k} = g}}
M_{e}(k) \otimes_{\So_{k}} \left(\bigoplus_{i_{1}+ \cdots +i_{k} = n}
Ind_{\So_{i_{1}}\times \cdots \times \So_{i_{k}}}^{\So_{n}} (N_{g_{1}}(i_{1})\otimes
\cdots \otimes N_{g_{k}}(i_{k}))\right).$$
The differential on $M\circ N$ is given by $d_{M\circ N} := d_{M}\circ id_{N} + id_{M}\circ' d_{N}$.\\
The differential on $M \circ (N,\, N')$ is given by
$$d_{M \circ (N,\, N')} := d_{M} \circ (id_{N},\, id_{N'}) + id_{M} \circ' (d_{N},\, id_{N'}) + id_{M} \circ (id_{N},\, d_{N'}).$$
Moreover, for any dg $\So$-modules $M,\, N$, we denote by $M \circ_{(1)} N$ the dg $\So$-module $M \circ (I,\, N)$. When $f : M \rightarrow M'$ and $g : N \rightarrow N'$, the map $f \circ (id_I,\, g) : M \circ_{(1)} N  \rightarrow M' \circ_{(1)} N'$ is denoted by $f\circ_{(1)} g$.

\subsection{(Co)operad}
An \emph{operad} is a monoid in the monoidal category of dg $\So$-modules with respect to the monoidal product $\circ$. A \emph{morphism of operads} is a morphism of dg $\So$-modules commuting with the operads structure. The notion of \emph{cooperad} is the dual version, i.e. a comonoid in the category of dg $\So$-modules. However, we use the invariants for the diagonal actions in the definition of the monoidal product instead of the coinvariants, that is,
$$\bigoplus_{k \geq 0} \left( M(k) \otimes \left(\bigoplus_{i_1 + \cdots + i_k = n} (N(i_1) \otimes \cdots \otimes N(i_k)) \otimes \K[\So_{n}] \right)^{\So_{i_{1}}\times \cdots \times \So_{i_{k}}} \right) ^{\So_{k}}.$$
Since we work over a field of characteristic $0$, the invariants are in one-to-one correspondence with the coinvariants and both definitions are equivalent. The definition with the invariants allows to define properly the signs.\\

The \emph{unit} of an operad $\Po$ is denoted by $\iota_{\Po} : I \rightarrow \Po$ and the \emph{counit} of a cooperad $\C$ is denoted by $\eta_{\C} : \C \rightarrow I$. Moreover when $(\Po,\, \gamma)$ is an operad, we define the \emph{partial product} $\gamma_p$ by
$$\Po \circ_{(1)} \Po \mono \Po \circ \Po \xrightarrow{\gamma} \Po$$
and when $(\C,\, \Delta)$ is a cooperad, we define the \emph{partial coproduct} $\Delta_p$ by
$$\C \xrightarrow{\Delta} \C \circ \C \epi \C \circ_{(1)} \C.$$

\begin{ex1}
Let $V$ be a dg $\K$-module. The dg $\So$-module $End(V):= \{\textrm{Hom}(V^{\otimes n},\, V)\}_{n\geq 0}$, endowed with the composition of maps, is an operad.
\end{ex1}

\subsection{Module over an operad and relative composition product}

A \emph{right $\Po$-module $(\mathcal{L},\, \rho)$} is an dg $\So$-module endowed with a map $\rho : \mathcal{L} \circ \Po \rightarrow \mathcal{L}$ compatible with the product and the unit of the operad $\Po$. We define similarly the notion of \emph{left $\Po$-module}.

We define the \emph{relative composition product $\mathcal{L} \circ_{\Po} \R$} between a right $\Po$-module $(\mathcal{L},\, \rho)$ and a left $\Po$-module $(\R,\, \lambda)$ by the coequalizer diagram
$$\xymatrix{\mathcal{L} \circ \Po \circ \R  \ar@<0.5ex>[r]^{\hspace{0.2cm} \rho \circ id_{\R}} \ar@<-0.5ex>[r]_{\hspace{0.2cm} id_{\mathcal{L}} \circ \lambda} & \mathcal{L} \circ \R \ar@{->>}[r] & \mathcal{L} \circ_{\Po} \R}.$$

\subsection{Algebra over an operad}
Let  $\Po$ be an operad. An \emph{algebra over the operad $\Po$}, or a \emph{$\Po$-algebra}, is a dg $\K$-module $V$ endowed with a morphism of operads $\Po \rightarrow End(V)$.

Equivalently, a $\Po$-algebra structure is given by a map $\gamma_V : \Po (V) \rightarrow V$ which is compatible with the composition product and the unity, where
$$\Po (V) := \Po \circ (V,\, 0,\, 0,\, \cdots) = \bigoplus_{n\geq 0} \Po(n) \otimes_{\So_{n}} V^{\otimes n}.$$

\subsection{Coalgebra over a cooperad}
Dually, let  $\C$ be a cooperad. A \emph{coalgebra over the cooperad $\C$}, or a \emph{$\C$-coalgebra}, is a dg $\K$-module $V$ endowed with a map $\delta : V \rightarrow \C(V) = \oplus_{n\geq 0}(\C (n) \otimes V^{\otimes n})^{\So_{n}}$ which satisfies compatibility properties. The notation $(-)^{\So_{n}}$ stands for the space of invariant elements.

\subsection{Module over a $\Po$-algebra}
Let $\Po$ be a dg $\So$-module and let $A$ be a dg vector space.
For a dg vector space $M$, we define the vector space $\Po (A,\, M)$ by the formula
$$\Po (A,\, M) := \Po \circ (A,\, M) = \bigoplus_n \Po(n) \otimes_{\So_n} (\overset{n}{\underset{j=1}{\oplus}} A \otimes \cdots \otimes \underbrace{M}_{j^{th} \textrm{ position}} \otimes \cdots \otimes A).$$
Let $(\Po,\, \gamma)$ be an operad and let $(A,\, \gamma_{A})$ be a $\Po$-algebra. An \emph{$A$-module $(M,\, \gamma_M,\ \iota_M)$}, or \emph{$A$-module over $\Po$}, is a vector space $M$ endowed with two maps $\gamma_M : \Po (A,\, M) \rightarrow M$ and $\iota_M : M \rightarrow \Po(A,\, M)$ such that the following diagrams commute

$$\makebox{\xymatrix{\Po (\Po (A),\, \Po (A,\, M)) \ar[rr]^{\hspace{12pt} id_{\Po} (\gamma_A,\, \gamma_M)} \ar[dd]_{\cong} & & \Po (A,\, M) \ar[d]^{\gamma_M}\\
& & M & \textrm{and}\quad \\
(\Po \circ \Po) (A,\, M) \ar[rr]_{\hspace{4pt} \gamma (id_{A},\, id_{M})} & & \Po (A,\, M) \ar[u]_{\gamma_M}}}
\raisebox{-20pt}{\xymatrix{M \ar[r]^{\hspace{-12pt} \iota_M} \ar[dr]_= & \Po (A,\, M) \ar[d]^{\gamma_M}\\
& M.}}$$

$$\hspace{2cm} (Associativity) \hspace{4cm} (Unitarity)$$
The category of $A$-modules over the operad $\Po$ is denoted by $\M_A^{\Po}$. The objects in $\M_A^{\Po}$ are differential graded $A$-modules over $\Po$. However, the morphisms in $\M_A^{\Po}$ are only maps of graded $A$-modules over $\Po$.

\begin{ex}
\begin{itemize}
\item[]
\item The operad $\Po = \Ass$ encodes associative algebras (not necessarily with unit). Then the map $\gamma_{n} : \Ass (n)\otimes_{\So_{n}}A^{\otimes n} \rightarrow A$ stands for the associative product of $n$ elements, where $\Ass (n) = \K[\So_{n}]$. We represent an element in $\Ass (n)$ by a corolla with $n$ entries. Then, an element in $\Ass (A,\, M)$ can be represented by $\vcenter{
\xymatrix@M=0pt@R=5pt@C=4pt{& & & &\\ a_{1} & a_{2} & \cdots & m & \cdots & a_{n}\\
\ar@{-}[ddrrr] & \ar@{-}[ddrr] & \ar@{-}[ddr] & \ar@{-}[dd] & \ar@{-}[ddl] & \ar@{-}[ddll] \\
&&&&&\\
 & & & \ar@{-}[d] &  \\ & & & &\\ & & & &}}$. However,
 $$\vcenter{
\xymatrix@M=0pt@R=5pt@C=4pt{a_{1} & \cdots & a_{k} & m & a_{k+1} & \cdots & a_{n}\\
\ar@{-}[ddrrr] & \ar@{-}[ddrr] & \ar@{-}[ddr] & \ar@{-}[dd] & \ar@{-}[ddl] & \ar@{-}[ddll] \\
&&&&&\\
 & & & \ar@{-}[d] &  \\ & & & &}} = \gamma \circ \gamma \circ \gamma \left(\vcenter{
\xymatrix@M=0pt@R=4pt@C=3pt{a_{1} & \cdots & a_{k} && m && a_{k+1} & \cdots & a_{n}\\
\ar@{-}[drdrdrdr] && \ar@{-}[drdr] && \ar@{-}[dr] && \ar@{-}[dldldl] && \ar@{-}[dldldldl]\\
\ar@{.}[rrrrrrrr] &&&&&&&&\\
& \ar@{.}[rrrrrr] &&&&&&&\\
&& \ar@{.}[rrrr] &&&&&&\\
&&& \ar@{.}[rr] & \ar@{-}[d] &&&&\\
&&&&&&&&}}\right),$$
then by several uses of the associativity diagram of $\gamma_{M}$, we get that an $A$-module over the operad $\Ass$ is given by two morphisms $A\otimes M \rightarrow M$ and $M\otimes A \rightarrow M$. Finally, we get the classical notion of dg $A$-bimodule.
\item The operad $\Po = \Com$ encodes classical associative and commutative algebras. We have $\Com (n) = \K$ and an element in $\Com (A,\, M)$ can be represented by $\vcenter{
\xymatrix@M=0pt@R=5pt@C=4pt{ a_{1} & a_{2} & \cdots & m & \cdots & a_{n}\\
\ar@{-}[ddrrr] & \ar@{-}[ddrr] & \ar@{-}[ddr] & \ar@{-}[dd] & \ar@{-}[ddl] & \ar@{-}[ddll] \\
&&&&&\\
 & & & \ar@{-}[d] &  \\ & & & &}}$ where the corolla is non-planar. Like before, an $A$-module structure over the operad $\Com$ is given by a morphism $A\otimes M \rightarrow M$. Hence, we get the classical notion of dg $A$-module.
\item The operad $\Po = \Lie$ encodes the Lie algebras. In this case, an $A$-module over the operad $\Lie$ is actually a classical dg Lie module or equivalently a classical associative module over the universal enveloping algebra of the Lie algebra $A$.
\end{itemize}
\end{ex}

Goerss and Hopkins defined in \cite{GoerssHopkins} a free $A$-module. We recall here the definition.

\begin{prop}[Proposition $1.10$ of \cite{GoerssHopkins}]\label{adjunction1}
The forgetful functor $U : \M_A^{\Po} \rightarrow g\mathcal{M}od_{\K}$ has a left adjoint, denoted by
$$N \mapsto A\otimes^{\Po} N.$$
That is we have an isomorphism of dg modules $$(\emph{Hom}_{\M_{A}^{\Po}}(A\otimes^{\Po} N,\, M),\, \partial) \cong (\emph{Hom}_{g\mathcal{M}od_{\K}} (N,\, U M),\, \partial)$$ for all $N \in g\mathcal{M}od_{\K}$ and $M \in \M_{A}^{\Po}$.
\end{prop}

A description of $A\otimes^{\Po} N$ is given by the following coequalizer diagram in $dg\mathcal{M}od_{\K}$
$$\xymatrix{\Po (\Po (A),\, N) \ar@<0.5ex>[r]^{c_{0}} \ar@<-0.5ex>[r]_{c_{1}} & \Po (A,\, N) \ar@{->>}[r] & A\otimes^{\Po} N,}$$
where the two first maps are given by the operad product
$$\Po (\Po (A),\, N) \mono (\Po \circ \Po)(A,\, N) \xrightarrow{\gamma (id_{A},\, id_{N})} \Po (A,\, N)$$
and the $\Po$-algebra structure
$$\Po (\Po (A),\, N) \xrightarrow{id_{\Po} (\gamma_{A},\, id_{N})} \Po (A,\, N).$$

\begin{rem}
We have to make note of the fact that the symbol $\otimes^{\Po}$ is just a notation and not a classical tensor product (except in the case $\Po = \Com$), as we will see in the following examples.
\end{rem}

\begin{ex}
\begin{itemize}
\item[]
\item When $\Po = \Ass$, we can write
$$\begin{array}{ll}
\vcenter{
\xymatrix@M=0pt@R=6pt@C=4pt{a_{1} & \cdots & a_{l} & n & a_{l+1} & \cdots & a_{k}\\
\ar@{-}[drrr] & \ar@{-}[drr] & \ar@{-}[dr] & \ar@{-}[d] & \ar@{-}[dl] & \ar@{-}[dll] & \ar@{-}[dlll]\\
 & & & \ar@{-}[d] & && \\ & & & & & &}} & = c_{0} \left(\vcenter{
\xymatrix@M=0pt@R=6pt@C=4pt{a_{1} & \cdots & a_{l} & n & a_{l+1} & \cdots & a_{k}\\
\ar@{-}[dr] & \ar@{-}[d] & \ar@{-}[dl] & \ar@{-}[ddd] & \ar@{-}[dr] & \ar@{-}[d] & \ar@{-}[dl]\\
\ar@{.}[rrrrrr] & \ar@{-}[drr] &&&& \ar@{-}[dll] &\\
& \ar@{.}[rrrr] &&&&&&\\
&&&&&&}}\right) = c_{1} \left(\vcenter{
\xymatrix@M=0pt@R=6pt@C=4pt{a_{1} & \cdots & a_{l} & n & a_{l+1} & \cdots & a_{k}\\
\ar@{-}[dr] & \ar@{-}[d] & \ar@{-}[dl] & \ar@{-}[ddd] & \ar@{-}[dr] & \ar@{-}[d] & \ar@{-}[dl]\\
\ar@{.}[rrrrrr] & \ar@{-}[drr] &&&& \ar@{-}[dll] &\\
& \ar@{.}[rrrr] &&&&&&\\
&&&&&&}}\right)\\
& \\
& = \vcenter{
\xymatrix@M=0pt@R=6pt@C=4pt{a_{1} \cdots a_{l} & n & a_{l+1} \cdots a_{k}\\
\ar@{-}[dr] & \ar@{-}[dd] & \ar@{-}[dl]\\
& &  \\ & &\\ & &}},
\end{array}$$
then we get $A\otimes^{\Ass} N \cong \vcenter{
\xymatrix@M=0pt@R=6pt@C=4pt{N\\
\ar@{-}[dd]\\ &\\ &}} \oplus \vcenter{
\xymatrix@M=0pt@R=6pt@C=4pt{A && N\\
\ar@{-}[dr] && \ar@{-}[dl]\\
& \ar@{-}[d] &  \\ & &}} \oplus \vcenter{
\xymatrix@M=0pt@R=6pt@C=4pt{N && A\\
\ar@{-}[dr] && \ar@{-}[dl]\\
& \ar@{-}[d] &  \\ & &}} \oplus \vcenter{
\xymatrix@M=0pt@R=6pt@C=4pt{A & N & A\\
\ar@{-}[dr] & \ar@{-}[dd] & \ar@{-}[dl]\\
& &  \\ & &}} \cong N\oplus A\otimes N \oplus N\otimes A \oplus A\otimes N \otimes A \cong (\K\oplus A)\otimes N \otimes (\K\oplus A)$ as modules over $\K$.
\item When $\Po = \Com$, we get $A\otimes^{\Com} N \cong N\oplus A\otimes N \cong (\K\oplus A)\otimes N$ as modules over $\K$.
\item When $\Po = \Lie$, we get $A\otimes^{\Lie} N \cong U^{e}(A) \otimes N$ as modules over $\K$, where $U^{e}(A)$ is the enveloping algebra of the Lie algebra $A$.
\end{itemize}
\end{ex}

These examples lead to the study of the $A$-module $A\otimes^{\Po} \K$ which is the \emph{enveloping algebra of the $\Po$-algebra} $A$ (defined in \cite{HinichSchechtman, GetzlerJones}). It has a multiplication given by
$$(A\otimes^{\Po} \K) \otimes (A\otimes^{\Po} \K) \cong A\otimes^{\Po} (A \otimes^{\Po} \K) \rightarrow A\otimes^{\Po} \K,$$
where the arrow is induced by the composition $\gamma$ of the operad (indeed, the kernel of the map $\Po (A,\, \Po (A,\, \K)) \epi A\otimes^{\Po} (A\otimes^{\Po} \K)$ is sent to $0$ by the map $\Po (A,\, \Po (A,\, \K)) \mono (\Po \circ \Po)(A,\, \K) \xrightarrow{\gamma (id_{A},\, id_{\K})} \Po (A,\, \K) \epi A\otimes^{\Po} \K$). This multiplication is associative and has a unit $\K \rightarrow A\otimes^{\Po} \K$.

\begin{prop}[Proposition $1.14$ of \cite{GoerssHopkins}]\label{Freyd-Mitchell}
The category $\M_{A}^{\Po}$ of $A$-modules over $\Po$ is isomorphic to the category of left unitary $A\otimes^{\Po} \K$-modules $g\mathcal{M}od_{A\otimes^{\Po} \K}$.
\end{prop}

\begin{rem}
We work in a differential graded setting. The differential on $A\otimes^{\Po}\K$ is induced by the differential on $\Po (A,\, \K)$. It is easy to see that the isomorphism is compatible with the graded differential framework.
\end{rem}

Given a map of $\Po$-algebras $B \xrightarrow{f} A$, there exists a forgetful functor $f^* : \M_{A}^{\Po} \rightarrow \M_{B}^{\Po}$, whose left adjoint gives the notion of free $A$-module on a $B$-module.

\begin{prop}[Lemma $1.16$ of \cite{GoerssHopkins}]\label{adjunction2}
The forgetful functor $f^* : \M_A^{\Po} \rightarrow \M_{B}^{\Po}$ has a left adjoint denoted by
$$N \mapsto f_{!} (N) := A\otimes_{B}^{\Po} N.$$
That is we have an isomorphism of dg modules
$$\emph{Hom}_{\M_{A}^{\Po}}(f_{!}(N),\, M) \cong \emph{Hom}_{\M_{B}^{\Po}} (N,\, f^* (M))$$ for all $M \in \M_{A}^{\Po}$ and $N \in \M_{B}^{\Po}$.
\end{prop}

It is also possible to make explicit the $A$-module $A\otimes_{B}^{\Po} N$ as the following coequalizer
$$\xymatrix{A\otimes^{\Po} (B\otimes^{\Po} N) \ar@<0.5ex>[r] \ar@<-0.5ex>[r] & A\otimes^{\Po} N \ar@{->>}[r] & A\otimes_{B}^{\Po} N.}$$
The module $A\otimes^{\Po} (B\otimes^{\Po} N)$ is a quotient of $\Po (A,\, \Po (B,\, N))$, then we define on $\Po (A,\, \Po (B,\, N))$ the composite
{\small $$\Po (A,\, \Po (B,\, N)) \xrightarrow{id_{\Po} (id_{A},\, id_{\Po} (f,\, id_{N}))} \Po (A,\, \Po (A,\, N)) \mono (\Po \circ \Po)(A,\, N) \xrightarrow{\gamma (id_{A},\, id_{N})} \Po (A,\, N) \epi A\otimes^{\Po} N.$$}
\hspace{-0.25cm} This map induced the first arrow $A\otimes^{\Po} (B\otimes^{\Po} N) \rightarrow A\otimes^{\Po} N$.\\
Similarly, the second map is induced by the composite
$$\Po (A,\, \Po (B,\, N)) \xrightarrow{id_{\Po} (id_{A},\, \gamma_{N})} \Po (A,\, N) \epi A\otimes^{\Po} N,$$
where $\gamma_{N}$ encodes the $B$-module structure on $N$.

\begin{rem}
The $A$-module $A\otimes_{B}^{\Po} N$ is a quotient of the free $A$-module $A\otimes^{\Po} N$. As for the notation $\otimes^{\Po}$, we have to be careful about the notation $\otimes_{B}^{\Po}$ which is not a classical tensor product over $B$ (except for $\Po = \Com$), as we see in the following examples.
\end{rem}

\begin{ex}
Provided a morphism of algebras $B \xrightarrow{f} A$, we have the dg $\K$-modules isomorphisms
\begin{itemize}
\item $A\otimes_{B}^{\Ass} N \cong (\K \oplus A)\otimes_{B} N\otimes_{B} (\K \oplus A)$, where the map $B\rightarrow \K \oplus A$ is given by $f$,
\item $A\otimes_{B}^{\Com} N \cong (\K \oplus A)\otimes_{B} N$, where the map $B\rightarrow \K \oplus A$ is given by $f$,
\item $A \otimes_{B}^{\Lie} N \cong U^{e}(A)\otimes_{B} N$, where $U^{e}(A)$ is the enveloping algebra of the Lie algebra $A$.
\end{itemize}
In all these examples, the notation $\otimes_{B}$ stands for the usual tensor product over $B$.
\end{ex}

\section{André-Quillen cohomology of algebras over an operad}

First we recall the conceptual definition of André-Quillen cohomology with coefficients of an algebra over an operad from \cite{Hinich, GoerssHopkins}. Then we recall the constructions and theorems of Koszul duality theory of operads \cite{GinzburgKapranov}. Finally, we recall the definition of twisting morphism given by \cite{GetzlerJones}. This section contains no new result but we will use these three theories throughout the text. We only want to emphasize that operadic resolutions from Koszul duality theory define functorial cofibrant resolutions on the level of algebras and then provide explicit chain complexes which compute André-Quillen cohomology.

We work with the cofibrantly generated model category of algebras over an operad and of modules over an operad given in \cite{GetzlerJones}, \cite{Hinich} and \cite{BergerMoerdijk}.\\

\subsection{Derivation and cotangent complex}

To study the structure of the $\Po$-algebra $A$, we derive the functor of $\Po$-derivations from $A$ to $M$ in the Quillen sense (non-abelian setting).

\subsubsection{\bf Algebras over a $\Po$-algebra}
Let $A$ be a $\Po$-algebra. A $\Po$-algebra $B$ endowed with an augmentation, that is a map of $\Po$-algebras $B\xrightarrow{f} A$, is called a \emph{$\Po$-algebra over $A$}. We denote by $\Po$-$\mathcal{A}lg/A$ the category of dg $\Po$-algebras over $A$ (the morphisms are given by the morphisms of graded algebras which commute with the augmentation maps).

\subsubsection{\bf Derivation}
Let $B$ be a $\Po$-algebra over $A$ and let $M$ be an $A$-module. An \emph{$A$-derivation from $B$ to $M$} is a linear map $d : B \rightarrow M$ such that the following diagram commutes
$$\xymatrix@C=2cm{\Po (B) = \Po \circ B \ar[d]_{\gamma_{B}} \ar[r]^{id_{\Po} \circ' d} & \Po (B,\, M) \ar[r]^{id_{\Po} \circ (f,\, id_{M})} & \Po (A,\, M) \ar[d]^{\gamma_{M}}\\
B \ar[rr]_d && M,}$$
where the infinitesimal composite of morphisms $\circ'$ was defined in \ref{dgsm}.
We denote by Der$_{A}(B,\, M)$ the set of $A$-derivations from $B$ to $M$.

This functor is representable on the right by the abelian extension of $A$ by $M$ and on the left by the $B$-module $\Omega_{\Po}B$ of K\"ahler differential forms as follows.

\subsubsection{\bf Abelian extension}
Let A be a $\Po$-algebra and let $M$ be an $A$-module. The \emph{abelian extension of $A$ by $M$}, denoted by $A \ltimes M$, is the $\Po$-algebra over $A$ whose underlying space is $A\oplus M$ and whose algebra structure is given by
$$\Po (A\oplus M) \epi \Po (A) \oplus \Po (A,\, M) \xrightarrow{\gamma_{A} + \gamma_{M}} A\oplus M.$$
The morphism $A\ltimes M \rightarrow A$ is just the projection on the first summand.

\begin{lem}[Definition $2.1$ of \cite{GoerssHopkins}]
Let $A$ be a $\Po$-algebra and $M$ be an $A$-module. Then there is an isomorphism of dg modules
$$\emph{Der}_{A}(B,\, M) \cong \emph{Hom}_{\Po \textrm{-} \mathcal{A}lg/A}(B,\, A\ltimes M).$$
\end{lem}

\begin{pf}
Any morphism of $\Po$-algebras $g : B \rightarrow A\ltimes M$ is the sum of the augmentation $B \rightarrow A$ and a derivation $d : B \rightarrow M$ and vice versa.
$\cqfd$
\end{pf}

\begin{lem}[Lemma $2.3$ of \cite{GoerssHopkins}]\label{bij}
Let $B$ be a $\Po$-algebra over $A$ and $M$ be an $A$-module. There is a $B$-module $\Omega_{\Po}B$ and an isomorphism of dg modules
$$\emph{Der}_{A}(B,\, M) \cong \emph{Hom}_{\M_{B}^{\Po}}(\Omega_{\Po}B,\, f^*(M)),$$
where the forgetful functor $f^*$ endows $M$ with a $B$-module structure.
Moreover, when $B = \Po (V)$ is a free algebra, we get $\Omega_{\Po}B \cong B\otimes^{\Po} V$.
\end{lem}

The second part of the lemma is given by the fact that Der$_{A}(\Po (V),\, M) \cong$ Hom$_{g\M od_{\K}} (V,\, M)$, that is any derivation from a free $\Po$-algebra is characterized by the images of its generators.

The $B$-module $\Omega_{\Po}B$ is called the \emph{module of K\"ahler differential forms}. It can be made explicit by the coequalizer diagram
$$\xymatrix{B\otimes^{\Po} \Po (B) \ar@<0.5ex>[r] \ar@<-0.5ex>[r] & B\otimes^{\Po} B \ar@{->>}[r] & \Omega_{\Po}B,}$$
where the first arrow is $B\otimes^{\Po} \gamma_{B}$ and the map
$$\Po (B,\, \Po (B)) \mono (\Po \circ \Po)(B,\, B) \xrightarrow{\gamma (id_{B},\, id_{B})} \Po (B,\, B) \epi B\otimes^{\Po}B$$
factors through $B\otimes^{\Po}\Po (B)$ to give the second arrow.

\begin{cor}
Let $B$ be a $\Po$-algebra over $A$ and $M$ be an $A$-module. There is an isomorphism of dg modules
$$\emph{Der}_{A}(B,\, M) \cong \emph{Hom}_{\M_{A}^{\Po}}(A\otimes_{B}^{\Po} \Omega_{\Po}B,\, M).$$
\end{cor}

\begin{pf}
We use Lemma \ref{bij} and the fact that $A\otimes_{B}^{\Po} -$ is left adjoint to the forgetful functor $f^*$ (Proposition \ref{adjunction2}).
$\cqfd$ 
\end{pf}

Finally, we get a pair of adjoint functors
$$A\otimes_{-}^{\Po}\Omega_{\Po} -\quad :\quad \Po \textrm{-} \mathcal{A}lg/A \rightleftharpoons \M_{A}^{\Po}\quad :\quad A\ltimes -.$$

From now on, we work over a ground field $\K$ of characteristic $0$.

We recall the model category structures on $\Po\textrm{-}\mathcal{A}lg/A$ and $\M_{A}^{\Po}$ given in \cite{Hinich}. It is obtained by the following transfer principle (see also \cite{GetzlerJones} and \cite{BergerMoerdijk}). Let $\D$ be a cofibrantly generated model category and let $\E$ be a category with small colimits and finite limits. Assume that $F : \D \rightleftharpoons \E : G$ is an adjunction with left adjoint $F$. Then the category $\E$ inherits a cofibrantly generated model category structure from $\D$, provided that $G$ preserves filtered colimits and that Quillen's small object (or Quillen's path-object) argument is verified. In this model category structure, a map $f$ in $\E$ is a weak equivalence (resp. fibration) if and only if $G(f)$ is a weak equivalence (resp. fibration) in $\D$.

In \cite{Hinich}, Hinich transfers the model category structure of the category of chain complexes over $\K$ to the category of $\Po$-algebras (see Theorem $4.1.1$ of \cite{Hinich}, every operad is $\Sigma$-split since $\K$ is of characteristic $0$). Finally, we obtain a model category structure on $\Po\textrm{-}\mathcal{A}lg/A$ in which $g : B \rightarrow B'$ is a weak equivalence (resp. a fibration) when the underlying map between differential graded modules is a quasi-isomorphism (resp. surjection). The category $\M_{A}^{\Po}$ of $A$-modules is isomorphic to the category $g\mathcal{M}od_{A\otimes^{\Po} \K}$ of differential graded module over the enveloping algebra $A\otimes^{\Po} \K$ (Proposition \ref{Freyd-Mitchell}). Then the category $\M_{A}^{\Po}$ inherits a model category structure in which $g : M \rightarrow M'$ is a weak equivalence (resp. a fibration) when $g$ is a quasi-isomorphism (resp. surjection) of $A\otimes^{\Po} \K$-modules.

\begin{prop}
The pair of adjoint functors
$$A\otimes_{-}^{\Po}\Omega_{\Po} -\quad :\quad \Po\textrm{-}\mathcal{A}lg/A \rightleftharpoons \M_{A}^{\Po}\quad :\quad A\ltimes -$$
forms a Quillen adjunction.
\end{prop}

\begin{pf}
By Lemma $1.3.4$ of \cite{Hovey}, it is enough to prove that $A\ltimes -$ preserves fibrations and acyclic fibrations. Let $g : M \epi M'$ be a fibration (resp. acyclic fibration) between $A$-modules. Then $g$ is a surjection (resp. a surjective quasi-isomorphism). The image of the map $g$ under the functor $A\ltimes -$ is $id_{A}\oplus g : A\ltimes M \rightarrow A\ltimes M'$, denoted by $id_{A}\ltimes g$. It follows that $id_{A}\ltimes g$ is surjective (resp. surjective and a quasi-isomorphism), which completes the proof.
$\cqfd$
\end{pf}

Thus, we consider the derived functors and we get the following adjunction between the homotopy categories
$$\Lb (A\otimes_{-}^{\Po}\Omega_{\Po} -)\quad :\quad \textrm{Ho}(\Po\textrm{-}\mathcal{A}lg/A) \rightleftharpoons \textrm{Ho}(\M_{A}^{\Po})\quad :\quad \mathbb{R} (A\ltimes -).$$
It follows that the cohomology of
$$\textrm{Hom}_{\textrm{Ho}(\M_{A}^{\Po})}(A\otimes_{R}^{\Po} \Omega_{\Po}R,\, M) \cong \textrm{Der}_{A}(R,\, M) \cong \textrm{Hom}_{\textrm{Ho}(\Po\textrm{-}\mathcal{A}lg/A)}(R,\, A\ltimes M)$$
is independent of the choice of the cofibrant resolution $R$ of $A$ in the model category of $\Po$-algebras over $A$.

\subsubsection{\bf André-Quillen (co)homology and cotangent complex}

Let $R \qiso A$ be a cofibrant resolution of $A$. The \emph{cotangent complex} is the total (left) derived functor of the previous adjunction and a representation of it is given by
$$\Lb_{R/A} := A\otimes_{R}^{\Po}\Omega_{\Po}R \in \textrm{Ho}(\M_{A}^{\Po}).$$
The \emph{André-Quillen cohomology of the $\Po$-algebra $A$ with coefficients in an $A$-module $M$} is defined by
$$\textrm{H}_{\Po}^{\bullet}(A,\, M) := \textrm{H}^{\bullet}(\textrm{Hom}_{\textrm{Ho}(\M_{A}^{\Po})}(\Lb_{R/A},\, M)).$$
The \emph{André-Quillen homology of the $\Po$-algebra $A$ with coefficients in an $A$-module $M$} is defined by
$$\textrm{H}^{\Po}_{\bullet}(A,\, M) := \textrm{H}_{\bullet}(M\otimes_{A\otimes^{\Po}\K} \Lb_{R/A}).$$
The study of the André-Quillen homology with coefficients is analogous to the study of the André-Quillen cohomology with coefficients. In this paper, we only work with André-Quillen cohomology with coefficients.

\begin{rem}
We use the left derived functor of the adjunction to define the André-Quillen cohomology. It is equivalent to define the André-Quillen cohomology by means of the right derived functor. We make this choice here because we are interested in considering homomorphisms in a modules category.
\end{rem}

\subsection{Bar construction of an operad and Koszul operad}\label{operadicresolution}

To make this cohomology theory explicit, we need a cofibrant resolution for algebras over an operad. In the model category of algebras over an operad, a cofibrant object is a retract of a quasi-free algebra endowed with a good filtration (for example, a non-negatively graded algebra). So we look for quasi-free resolutions of algebras. Operadic resolutions provide such functorial cofibrant resolutions for algebras. There are mainly three operadic resolutions: the simplicial bar construction which induces a Godement type resolution for algebras, the (co)augmented (co)bar construction on the level of (co)operads and the Koszul complex for operads. This last one induces the bar-cobar resolution (or Boardman-Vogt resolution \cite{BoardmanVogt, BergerMoerdijk2}) on the level of algebras. The aim of the two next subsections is to recall the operadic resolutions.

Here, we briefly recall the (co)bar construction of a (co)operad and the notion of Koszul operad. We refer to \cite{GinzburgKapranov, GetzlerJones, Fresse} for a complete exposition.

\subsubsection{\bf Bar construction}\label{barcons}
Let $\Po$ be an augmented operad. We denote by $sV$ the suspension of $V$ (that-is-to-say $(sV)_{d} := V_{d-1}$). The \emph{bar construction of $\Po$} is the quasi-free cooperad
$$\B(\Po) := (\F ^c(s\overline{\Po}),\, d_{\B(\Po)} := d_{1} - d_{2}),$$
where the map $d_{1}$ is induced by the internal differential of the operad ($d_{s\overline{\Po}} := id_{\K s}\otimes d_{\Po}$) and the component $d_{2}$ is induced by the product of the operad by
$$\F _{(2)}(s\overline{\Po}) \cong \bigoplus_{2\textrm{-vertices trees}} \K s \otimes \overline{\Po} \otimes \K s \otimes \overline{\Po} \xrightarrow{id_{\K s} \otimes \tau \otimes id_{\overline{\Po}}}  \bigoplus_{2\textrm{-vertices trees}} \K s \otimes \K s \otimes \overline{\Po} \otimes \overline{\Po} \xrightarrow{\Pi_{s} \otimes \gamma_{\Po}} \K s \otimes \overline{\Po},$$
where $\tau : \overline{\Po} \otimes \K s \rightarrow \K s \otimes \overline{\Po}$ is the \emph{symmetry isomorphism} given explicitly by $\tau (o_{1}\otimes o_{2}) := (-1)^{|o_{1}| |o_{2}|} o_{2}\otimes o_{1}$ and $\Pi_{s} : \K s\otimes \K s \rightarrow \K s$ is the morphism of degree $-1$ induced by $\Pi_{s}(s\otimes s) := s$.

\begin{rem}
Assume that $\Po$ is weight graded. Then the bar construction is bigraded by the number $(w)$ of non-trivial indexed vertices and by the total weight $(\rho)$
$$\B_{(w)}(\Po) := \oplus_{\rho \in \mathbb{N}} \B_{(w)}(\Po)^{(\rho)}.$$
\end{rem}

Dually, we define the \emph{cobar construction of a coaugmented cooperad $\C$} by $$\Omega(\C) := (\F (s^{-1}\overline{\C}),\, d_{1} - d_{2}).$$

From now on, we assume that $\Po$ is an augmented operad and $\C$ is a coaugmented cooperad.

\subsubsection{\bf Quadratic operad}
A operad $\Po$ is \emph{quadratic} when $\Po = \F (V)/(R)$, where $V$ is the $\So$-module of generators, $\F (V)$ is the free operad and the space of relations $R$ lives in $\F _{(2)}(V)$, the set of trees with two vertices. We endow $\F (V)$ with a weight grading, which differs from the homological degree, given by the number of vertices, this induces a weight grading on each quadratic operad. In this paper, we consider only non-negatively weight graded operad and we say that a weight graded dg operad $\Po$ is connected when $\Po = \K \oplus \Po^{(1)} \oplus \Po^{(2)} \oplus \cdots$, where $\Po^{(0)} = \K$ is concentrated in homological degree $0$.

\subsubsection{\bf Koszul operad}
We define the \emph{Koszul dual cooperad of $\Po$} by the weight graded dg $\So$-module
$$\Po^{\ash}_{(\rho)} := \textrm{H}_{\rho}(\B_{(\bullet)}(\Po)^{(\rho)},\, d_{2}).$$
An operad is called a \emph{Koszul} operad when the injection $\Po^{\ash} \mono \B(\Po)$ is a quasi-isomorphism.

When $\Po$ is of finite type, that is $\Po(n)$ is finite dimensional for each $n$, we can dualize linearly the cooperad $\Poa$ to get the \emph{Koszul dual operad of $\Po$}, denoted by $\Po^{!}$. For any $\So_{n}$-module $V$, we denote by $V^{\vee}$ the $\So_{n}$-module $V^{*}\otimes (sgn_{n})$, where $(sgn_{n})$ is the one-dimensional signature representation of $\So_{n}$. We define $\Po^!(n) := \Poa(n)^{\vee}$. The product on $\Po^!$ is given by $^t\Delta_{\Poa} \circ \omega$ where $\omega : {\Poa}^{\vee} \circ {\Poa}^{\vee} \rightarrow (\Poa \circ \Poa)^{\vee}$.

\subsubsection{\bf Algebras up to homotopy}

Let $\Po$ be a Koszul operad. We define $\Po_{\infty} := \Omega (\Poa)$. A $\Po_{\infty}$-algebra is called an \emph{algebra up to homotopy} or \emph{homotopy $\Po$-algebra} (see \cite{GinzburgKapranov}). The notion of $\Po_{\infty}$-algebras is a lax version of the notion of $\Po$-algebras.

\begin{ex}
\begin{itemize}
\item[]
\item When $\Po = \Ass$, we get the notion of $A_{\infty}$-algebras;
\item when $\Po = \Lie$, we get the notion of $L_{\infty}$-algebras;
\item when $\Po = \Com$, we get the notion of $C_{\infty}$-algebras.
\end{itemize}
\end{ex}

\subsection{Operadic twisting morphism}\label{deftor}

We refer to \cite{GetzlerJones, MerkulovVallette} for a general and complete treatment. Let $\alpha$, $\beta : \C \rightarrow \Po$ be morphisms of $\So$-modules. We define the convolution product
$$\alpha \star \beta : \C \xrightarrow{\Delta_{p}} \C \circ_{(1)} \C \xrightarrow{\alpha \circ_{(1)} \beta} \Po \circ_{(1)} \Po \xrightarrow{\gamma_{p}} \Po.$$
The $\So$-module Hom$(\C,\, \Po)$ is endowed with an operad structure. Moreover, the convolution product is a pre-Lie product on Hom$(\C,\, \Po)$, that is, it satisfies the relation
$$(\alpha \star \beta) \star \gamma - \alpha \star (\beta \star \gamma) = (-1)^{|\beta||\gamma|}[(\alpha \star \gamma) \star \beta - \alpha \star (\gamma \star \beta)] \textrm{ for all $\alpha$, $\beta$ and $\gamma$ in $\Hom (\C,\, \Po)$}.$$

\subsubsection{\bf Definition}
An \emph{operadic twisting morphism} is a map $\alpha : \C \rightarrow \Po$ of degree $-1$ satisfying the \emph{Maurer-Cartan equation}
$$\partial(\alpha) + \alpha \star \alpha = 0.$$
We denote the set of operadic twisting morphisms from $\C$ to $\Po$ by \textrm{Tw}($\C$, $\Po$).

In the weight graded case, we assume that the twisting morphisms and the internal differentials preserve the weight.

\begin{thm}[Theorem $2.17$ of \cite{GetzlerJones}]\label{GJres}
The functors $\Omega$ and $\B$ form a pair of adjoint functors between the category of connected coaugmented cooperads and augmented operads. The natural bijections are given by the set of operadic twisting morphisms:
$$\emph{Hom}_{dg-Op}(\Omega(\C),\, \Po) \cong \emph{Tw($\C$, $\Po$)} \cong \emph{Hom}_{dg-Coop}(\C,\, \B(\Po)).$$
\end{thm}

\begin{ex}
We give examples of operadic twisting morphisms.
\begin{itemize}
\item When $\C = \B(\Po)$ is the bar construction on $\Po$, the previous theorem gives a natural operadic twisting morphism $\pi : \B(\Po) = \mathcal{F}^c(s\overline{\Po}) \epi s \overline{\Po} \xrightarrow{s^{-1}} \overline{\Po} \mono \Po$. This morphism is universal in the sense that each twisting morphism $\alpha : \C \rightarrow \Po$ factorizes uniquely through the map $\pi$
$$\xymatrix{\C \ar[rr]^{\alpha} \ar@{-->}[dr]_{f_{\alpha}} && \Po\\
& \B (\Po), \ar[ur]_{\pi}}$$
where $f_{\alpha}$ is a morphism of dg cooperads.
\item When $\C = \Po^{\ash}$ is the Koszul dual cooperad of a quadratic operad $\Po$, the map $\kappa : \Po^{\ash} \mono \B(\Po) \xrightarrow{\pi} \Po$ is an operadic twisting morphism (the precomposition of an operadic twisting morphism by a map of dg cooperads is an operadic twisting morphism). Actually we have $\Po^{\ash} \mono \F ^c(sV)$ and the map $\kappa$ is given by $\Po^{\ash} \epi \Poa_{(1)} \cong sV \xrightarrow{s^{-1}} V \mono \Po$.
\item When $\Po = \Omega(\C)$ is the cobar construction on $\C$, the previous theorem gives a natural operadic twisting morphism $\iota : \C \epi \overline{\C} \xrightarrow{s^{-1}} s^{-1}\overline{\C} \mono \Omega(\C) = \mathcal{F}(s^{-1}\overline{\C})$. This morphism is universal in the sense that each twisting morphism $\alpha : \C \rightarrow \Po$ factorizes uniquely through the map $\iota$
$$\xymatrix{& \Omega (\C) \ar@{-->}[dr]^{g_{\alpha}}\\
\C \ar[rr]^{\alpha} \ar[ur]^{\iota} && \Po,
}$$
where $g_{\alpha}$ is a morphism of dg operads.
\end{itemize}
\end{ex}

\subsubsection{\bf Twisted composition product}\label{twistedcompositionproduct}

Let $\Po$ be a dg operad and let $\C$ be a dg cooperad. Let $\alpha : \C \rightarrow \Po$ be an operadic twisting morphism. The \emph{twisted composition product} $\Po \circ_{\alpha} \C$ is the $\So$-module $\Po \circ \C$ endowed with a differential $d_{\alpha} := d_{\Po \circ \C} - \delta_{\alpha}^{l}$, where $\delta_{\alpha}^{l}$ is defined by the composite
$$\delta_{\alpha}^{l} : \Po \circ \C \xrightarrow{id_{\Po} \circ' \Delta_{\C}} \Po \circ (\C,\, \C \circ \C) \xrightarrow{id_{\Po} \circ (id_{\C},\, \alpha \circ id_{\C})} \Po \circ (\C,\, \Po \circ \C) \mono (\Po \circ \Po) \circ \C \xrightarrow{\gamma \circ id_{\C}} \Po \circ \C.$$
Since $\alpha$ is an operadic twisting morphism, $d_{\alpha}$ is a differential.

When $A$ is a $\Po$-algebra, we denote by $\C \circ_{\alpha} A$ the chain complex $(\C (A),\, d_{\alpha} := d_{\C (A)} + \delta_{\alpha}^{r})$, where $\delta_{\alpha}^{r}$ is the composite
$$\C (A) \xrightarrow{\Delta_{p}\circ id_{A}} (\C \circ_{(1)} \C) (A) \xrightarrow{id_{\C}\circ_{(1)}  \alpha \circ id_{A}} (\C \circ_{(1)} \Po) (A) \xrightarrow{id_{\C} \circ \gamma_{A}} C (A).$$

Finally, we denote by $\Po \circ_{\alpha} \C \circ_{\alpha} A$ the vector space $\Po \circ \C (A)$ endowed with the differential
$$d_{\alpha} := d_{\Po \circ \C (A)} - \delta_{\alpha}^{l} \circ id_{A} + id_{\Po} \circ' \delta_{\alpha}^r = d_{\Po \circ (\C \circ_{\alpha} A)} - \delta_{\alpha}^{l} \circ id_{A}.$$

The notation $d_{\alpha}$ stands for different differentials. The differential is given without ambiguity by the context.

\subsubsection{\bf Operadic resolutions}
In \cite{GetzlerJones}, Getzler and Jones produced functorial resolutions of algebras given by the following theorems.

\begin{thm}[Theorem $2.19$ of \cite{GetzlerJones}]\label{barcobarres}
The augmented bar construction gives a resolution
$$\xymatrix{\Po \circ_{\pi} \B(\Po) \circ_{\pi} A \ar@{->>}[r]^{\hspace{0.9cm}\sim} & A.}$$
\end{thm}

\begin{thm}[Theorem $2.25$ of \cite{GetzlerJones}]
When the operad $\Po$ is Koszul, there is a smaller resolution of $A$ given by the Koszul complex
$$\xymatrix{\Po \circ_{\kappa} \Po^{\ash} \circ_{\kappa} A \ar@{->>}[r]^{\hspace{0.7cm}\sim} & A.}$$
\end{thm}

The augmented bar resolution admits a dual version.

\begin{thm}[Theorem $4.18$ of \cite{Vallette2}]
For every weight graded coaugmented cooperad $\C$, there is an isomorphism
$$\Omega (\C) \circ_{\iota} \C \qiso I.$$
\end{thm}

This gives, for all $\Omega (\C)$-algebra $A$, a quasi-isomorphism $\Omega (\C) \circ_{\iota} \C \circ_{\iota} A \qiso A$.

\subsection{Description of the cotangent complex}

Thanks to these resolutions, we can describe the underlying vector space of the cotangent complex.

\subsubsection{\bf Quasi-free resolution}\label{quasifreeresolution}

Let $A$ be a $\Po$-algebra, let $C$ be a $\C$-coalgebra endowed with a filtration $F_pC$ such that $F_{-1}C = \{0\}$ and let $\alpha : \C \rightarrow \Po$ be an operadic twisting morphism. We denote by $\Po \circ_{\alpha} C$ the complex  $(\Po(C),\, d_{\alpha} := d_{\Po \circ C} - \delta_{\alpha}^{l})$. The differential $\delta_{\alpha}^{l}$ on $\Po (C)$ is given by
 $$\delta_{\alpha}^{l} : \Po (C) \xrightarrow{id_{\Po} \circ' \Delta} \Po \circ (C,\, \C (C)) \xrightarrow{id_{\Po} \circ (id_{C},\, \alpha \circ id_{C})} \Po \circ (C,\, \Po (C)) \mono \Po \circ \Po (C) \xrightarrow{\gamma \circ id_{c}} \Po (C).$$
A \emph{quasi-free resolution of $A$} is a complex $\Po \circ_{\alpha} C$ such that $\Po \circ_{\alpha} C \qiso A$ and ${\delta_{\alpha}^l}_{|F_pC} \subset \Po (F_{p-1}C)$.

Except the normalized cotriple construction, all the previous resolutions are of this form when $A$ is non-negatively graded. With this resolution, we make the cotangent complex explicit.

\begin{thm}\label{cotangentcomplex}
Let $\Po (C)$ be a quasi-free resolution of the $\Po $-algebra $A$. With this resolution, the cotangent complex has the form
$$\Lb_{\Po (C)/A} \cong A\otimes^{\Po}C.$$
\end{thm}

\begin{pf}
The cotangent complex is isomorphic to
$$\begin{array}{llll}
A\otimes_{R}^{\Po}\Omega_{\Po}R & = & A\otimes_{\Po (C)}^{\Po}\Omega_{\Po}(\Po (C))&\\
& \cong & A\otimes_{\Po (C)}^{\Po}(\Po (C)\otimes^{\Po}C) & (\textrm{Lemma \ref{bij}})\\
& \cong & A\otimes^{\Po}C & (\textrm{Propositions \ref{adjunction1} and \ref{adjunction2}}).
\end{array}$$
$\cqfd$
\end{pf}

When we use the augmented bar construction, we get the cotangent complex for any algebra over any operad. However this complex may be huge and it can be useful to work with smaller resolutions. When we use the Koszul resolution, we can use the Koszul complex and we get the cotangent complex of an algebra over a Koszul operad. For homotopy algebras, we use the coaugmented cobar construction. In this paper, we consider only resolutions coming from operadic resolutions. In \cite{Milles2}, we work with even smaller resolutions, but which are not functorial with respect to the algebra.\\

To describe completely the cotangent complex, we have to make its differential explicit. In the next section, we will trace the boundary map on Der$_{A}(R,\, M)$ through the various isomorphisms.

\section{Lie theoretic description}\label{twist}

We endow the chain complex defining the André-Quillen cohomology with a structure of Lie algebra. The notion of twisting morphism (or twisting cochain) first appeared in \cite{Brown} and in \cite{Moore} (see also \cite{HusemollerMooreStasheff}). It is a particular kind of maps between a coassociative coalgebra and an associative algebra. Getzler and Jones extend this definition to (co)algebras over (co)operads (see $2.3$ of \cite{GetzlerJones}). We show that the differential on the cotangent complex $A\otimes^{\Po} C$ is obtained by twisting the internal differential by a twisting morphism.\\

In the sequel, let $(\Po,\, \gamma)$ denote an operad, $(\C,\, \Delta)$ denote a cooperad and $(C,\, \Delta_{C})$ denote a $\C$-coalgebra.

\subsection{A Lie algebra structure}

Let $\alpha : \C \rightarrow \Po$ be an operadic twisting morphism. Let $C$ be a $\C$-coalgebra and let $A$ be a $\Po$-algebra. Let $M$ be an $A$-module. For all $\varphi$ in Hom$_{g\mathcal{M}od_{\K}}(C,\, A)$ and $g$ in Hom$_{g\mathcal{M}od_{\K}}(C,\, M)$, we define $\alpha[\varphi,\, g] := \sum_{n\geq 1} \alpha[\varphi,\, g]_n$, where $\alpha[\varphi,\, g]_n$ is the composite
$$C \xrightarrow{\Delta_{C}} \C (C) \epi (\C(n) \otimes C^{\otimes n})^{\So_n} \xrightarrow{\alpha \otimes \varphi^{\otimes n-1} \otimes g} \Po(n) \otimes A^{\otimes n-1} \otimes M \epi \Po (A,\, M) \xrightarrow{\gamma_M} M.$$
The notation $\otimes_{H}$ stands for the Hadamard product: for any $\So$-modules $M$ and $N$, $(M \otimes_{H} N) (n) := M(n)\otimes N(n)$. Let $End_{s^{-1}\K}$ be the cooperad defined by
$$End_{s^{-1}\K}(n) := \textrm{Hom}((s^{-1}\K)^{\otimes n},\, s^{-1}\K)$$
endowed with the natural action of $\So_{n}$. When $(C,\, \Delta_C)$ is a $\C$-coalgebra, we endow $s^{-1}C := s^{-1}\K \otimes C$ with a structure of $End_{s^{-1}\K} \otimes_{H} \C$-coalgebra given by
$$\Delta_{s^{-1}C} : s^{-1}C \xrightarrow{\Delta_{C}(n)} s^{n-1}s^{-n}(\C(n)\otimes C^{\otimes n})^{\So_n} \xrightarrow{\tau_n} \left((End_{s^{-1}\K}(n) \otimes \C(n))\otimes (s^{-1}C)^{\otimes n}\right)^{\So_n},$$
where $\Delta_{C}(n)$ is the composite $C \xrightarrow{\Delta_C} \C (C) \epi (\C(n) \otimes C^{\otimes n})^{\So_n}$ and $\tau_n$ is a map which permutes components and is induced by compositions of $\tau$ (seen in Section $1.2.1$). The differential on $s^{-1}C$ is given by $d_{s^{-1}C} := id_{s^{-1}\K}\otimes d_C$.\\

In the following results, the operad $\Po$ is quadratic and binary and the cooperad $\C = \Poa$ is the Koszul dual cooperad of $\Po$. The twisting morphism $\kappa : \Poa \rightarrow \Po$ is defined in the examples after Section \ref{GJres}.

\begin{thm}\label{algLie}
Let $\Po$ be a quadratic binary operad and let $\C = \Po^{\ash}$ be the Koszul dual cooperad of $\Po$. Let $A$ be a $\Po$-algebra and $C$ be a $\Po^{\ash}$-coalgebra. The chain complex $$(\emph{Hom}_{g\mathcal{M}od_{\K}}(C,\, A),\, \kappa [-,\, -],\, \partial)$$
forms a dg Lie algebra whose bracket $\kappa [-,\, -]$ is of degree $-1$, that is
$$\kappa [\varphi,\, \psi] = - (-1)^{(|\varphi|-1)(|\psi|-1)}\kappa [\psi,\, \varphi].$$
\end{thm}

\begin{pf}
There is an isomorphism of chain complexes
$$\begin{array}{ccc}
\textrm{Hom}_{g\mathcal{M}od_{\K}}^{\bullet}(C,\, A) & \iso & \textrm{Hom}_{g\mathcal{M}od_{\K}}^{\bullet +1}(s^{-1}C,\, A)\\
\varphi & \mapsto & (\bar{\varphi} : s^{-1}c \mapsto \varphi(c)),
\end{array}$$
since $\overline{\partial(\varphi)} = \overline{d_A \circ \varphi - (-1)^{|\varphi|}\varphi \circ d_{C}} = d_A \circ \bar{\varphi} - (-1)^{|\varphi|-1}\bar{\varphi} \circ d_{s^{-1}C} = \partial(\bar{\varphi})$. Moreover, we have the equality $\overline{\kappa[\varphi,\, \psi]} = (-1)^{|\bar{\varphi}|}\bar{\kappa}[\bar{\varphi},\, \bar{\psi}]$, where $\bar{\kappa}(s^{n-1}\mu^{c}) := \kappa(\mu^{c})$ is not a map of $\So_{n}$-modules.

We show now that the dg module
$$(\textrm{Hom}_{g\mathcal{M}od_{\K}}^{\bullet}(s^{-1}C,\, A),\, (-1)^{|\bar{\varphi}|}\bar{\kappa}[\bar{\varphi},\, \bar{\psi}],\, \partial)$$
forms a Lie algebra. Since $C$ is a $\Poa$-coalgebra, we get that $(s^{-1}C)^* \cong sC^*$ is a $\Po^{!}$-algebra. That is, there is a morphism of operads $\Po^{!} \rightarrow End(sC^*)$. Hence, we obtain a morphism $\Po^{!}\otimes_{H}\Po \rightarrow End(sC^*)\otimes_{H}End(A) \cong End(sC^*\otimes A)$. We apply Theorem $29$ of \cite{Vallette} and we get that Hom$_{g\mathcal{M}od_{\K}}(s^{-1}C,\, A) \cong sC^{*}\otimes A$ is a Lie algebra. The Lie algebra structure is given by $(-1)^{|\bar{\varphi}|}\bar{\kappa}[\bar{\varphi},\, \bar{\psi}]$, which is of degree $0$ since $\kappa$ is non-zero only on $\Poa(2)$. Therefore $\textrm{Hom}_{g\mathcal{M}od_{\K}}^{\bullet +1}(s^{-1}C,\, A)$ is a Lie algebra with bracket of degree $0$.
$\cqfd$
\end{pf}

\begin{thm}
Let $\Po$ be a quadratic binary operad and take $\C = \Po^{\ash}$. Let $A$ be a $\Po$-algebra, let $C$ be a $\C$-coalgebra and let $M$ be an $A$-module. Then the dg module $$(\emph{Hom}_{g\mathcal{M}od_{\K}}(C,\, M),\, \kappa[-,\, -],\, \partial)$$
is a dg Lie module over $(\emph{Hom}_{g\mathcal{M}od_{\K}}(C,\, A),\, \kappa[-,\, -],\, \partial)$.
\end{thm}

\begin{pf}
The proof is analoguous to the proof of Theorem \ref{algLie} in the following way. A $A$-module structure over the operad $\Po$ is equivalent to a map of operads $\Po \rightarrow End_{A}(M)$, where $End_{A}(M) := End(A)\oplus End(A,\, M)$ with
$$End(A,\, M)(n) := \bigoplus_{j=1}^{n} \textrm{Hom}(\underbrace{A\otimes \cdots \otimes A}_{\textrm{$j-1$ times}} \otimes M \otimes \underbrace{A\otimes \cdots \otimes A}_{\textrm{$n-j$ times}},\, M).$$
The composition product is given by the composition of maps when possible and zero otherwise. We get Hom$_{g\mathcal{M}od_{\K}}(s^{-1}C,\, M) \cong sC^{*}\otimes M$ and there is a map of operads $\Lie \rightarrow \Po^{!}\otimes_{H} \Po \rightarrow End(sC^{*})\otimes End_{A} (M) \cong End_{sC^{*}\otimes A}(sC^{*}\otimes M)$. Therefore, Hom$_{g\mathcal{M}od_{\K}}(C,\, M)$ is a dg Lie module over Hom$_{g\mathcal{M}od_{\K}}(C,\, A)$.
$\cqfd$
\end{pf}

\subsection{Algebraic twisting morphism}\label{twistexamples}

In this section, we define the notion of twisting morphism on the level of (co)algebras introduced in $2.3$ of \cite{GetzlerJones}. Assume now that $\alpha : \C \rightarrow \Po$ is an operadic twisting morphism. Let $A$ be a $\Po$-algebra and let $C$ be a $\C$-coalgebra. For all $\varphi$ in Hom$_{\mathcal{M}od_{\K}}(C,\, A)$, we define the maps
$$\star_{\alpha} (\varphi) : C \xrightarrow{\Delta_{C}} \C (C) \xrightarrow{\alpha \circ \varphi} \Po (A) \xrightarrow{\gamma_{A}} A.$$

An \emph{algebraic twisting morphism with respect to $\alpha$} is a map $\varphi : C \rightarrow A$ of degree $0$ satisfying the Maurer-Cartan equation
$$\partial(\varphi) + \star_{\alpha} (\varphi) = 0.$$
We denote by Tw$_{\alpha}(C, A)$ the set of algebraic twisting morphisms with respect to $\alpha$.

\begin{ex}
We consider the two examples of Section \ref{deftor} once again.
\begin{itemize}
\item The map $\eta_{\B(\Po)}(A) := \eta_{\B(\Po)}\circ id_A : \B(\Po) (A) \epi I \circ A \cong A$ is an algebraic twisting morphism with respect to $\pi$. For simplicity, assume $d_{A} = 0$. We get
$$\begin{array}{ll}
\partial(\eta_{\B(\Po)}(A)) & = d_{A}\circ \eta_{\B(\Po)}(A) - \eta_{\B(\Po)}(A) \circ d_{\pi}^{r}\\
& =  -\eta_{\B(\Po)}(A)\circ (d_{\B(\Po)}\circ id_{A} + \delta_{\pi}^r)\\
& =  -\eta_{\B(\Po)}(A)\circ \delta_{\pi}^r
\end{array}$$
since $d_{\B(\Po)} = 0$ on $\F _{(0)}(s\overline{\Po})$. Then $\partial (\eta_{\B(\Po)}(A)) (e)$ is non-zero if and only if $e = s\mu \otimes (a_{1}\otimes \cdots \otimes a_{n})\in \F _{(1)}(s\overline{\Po})(A)$ and is equal to $-\mu (a_{1},\cdots , a_{n})$ in this case. Moreover, $\star_{\pi}(\eta_{\B(\Po)}(A))$ satisfies the same properties. So the assertion is proved.
\item The map $\eta_{\Po^{\ash}}(A) : \Po^{\ash} (A) \mono \B(\Po)(A) \epi A$ is an algebraic twisting morphism with respect to $\kappa$.
\end{itemize}
Let us now make explicit the maps $\kappa$ and $\eta_{\Po^{\ash}}(A)$ in the cases $\Po = \Ass$, $\Com$ and $\Lie$. We refer to \cite{Vallette} for the categorical definition of the Koszul dual cooperad.
\begin{itemize}
\item When $\Po = \Ass$, the Koszul dual $\Ass^{\ash}$ is a cooperad cogenerated by the elements $s\vcenter{\xymatrix@M=0pt@R=4pt@C=4pt{
\ar@{-}[dr] &  &\ar@{-}[dl]  \\
 &\ar@{-}[d] &  \\  & &}}$, that is the elements $\vcenter{\xymatrix@M=0pt@R=4pt@C=4pt{
\ar@{-}[dr] &  &\ar@{-}[dl]  \\
 &\ar@{-}[d] &  \\  & &}} \in \Ass(2)$ suspended by an $s$ of degree $1$, with corelations $s\vcenter{\xymatrix@M=0pt@R=4pt@C=4pt{
\ar@{-}[dr] &  &\ar@{-}[dl]  \\
 &\ar@{-}[d] &  \\  & &}} \otimes (s\vcenter{\xymatrix@M=0pt@R=4pt@C=4pt{
\ar@{-}[dr] &  &\ar@{-}[dl]  \\
 &\ar@{-}[d] &  \\  & &}}\otimes \vcenter{\xymatrix@M=0pt@R=4pt@C=4pt{
\ar@{-}[dd]&  \\
 & \\  &}}) - s\vcenter{\xymatrix@M=0pt@R=4pt@C=4pt{
\ar@{-}[dr] &  &\ar@{-}[dl]  \\
 &\ar@{-}[d] &  \\  & &}} \otimes (\, \vcenter{\xymatrix@M=0pt@R=4pt@C=4pt{
\ar@{-}[dd]&  \\
 & \\  &}}\otimes s\vcenter{\xymatrix@M=0pt@R=4pt@C=4pt{
\ar@{-}[dr] &  &\ar@{-}[dl]  \\
 &\ar@{-}[d] &  \\  & &}})$, that we can represent by $s^2(\vcenter{\xymatrix@M=0pt@R=4pt@C=4pt{
& & & & \\ \ar@{-}[ddrr] & & \ar@{-}[dl] & & \ar@{-}[ddll]\\
&  & &  & \\
& & \ar@{-}[d] & &\\
& & }} - \vcenter{\xymatrix@M=0pt@R=4pt@C=4pt{& & & & \\
\ar@{-}[drdr] & &\ar@{-}[dr] & & \ar@{-}[dldl]  \\
& & & & \\
& &\ar@{-}[d] & & \\
& & }})$. The map $\kappa : \Ass^{\ash} \rightarrow \Ass$ sends $s\vcenter{\xymatrix@M=0pt@R=4pt@C=4pt{
\ar@{-}[dr] & &\ar@{-}[dl]  \\
 &\ar@{-}[d] &  \\  & &}}$ onto $\vcenter{\xymatrix@M=0pt@R=4pt@C=4pt{
\ar@{-}[dr] &  &\ar@{-}[dl]  \\
 &\ar@{-}[d] &  \\  & &}}$ and is zero elsewhere. The map $\eta_{\Ass^{\ash}}(A)$ sends $A$ onto $A$ and is zero elsewhere.
\item When $\Po = \Com$, the map $\kappa$ sends the cogenerator of $\Com^{\ash}$ on the generator of $\Com$ and is zero outside $\Com^{\ash}(2)$. The map $\eta_{\Com^{\ash}}(A)$ is just the projection onto $A$.
\item When $\Po = \Lie$, the map $\kappa$ sends the cogenerator of $\Lie^{\ash}$ on the generator of $\Lie$ and is zero outside $\Lie^{\ash}(2)$ and the map $\eta_{\Lie^{\ash}}(A)$ is just the projection onto $A$.
\end{itemize}
\end{ex}

When $\Po$ is a binary quadratic operad, $\C = \Po^{\ash}$ is its Koszul dual cooperad and $\alpha = \kappa$, then algebraic twisting morphisms with respect to $\kappa$ are in one-to-one correspondence with solutions of the Maurer-Cartan equation in the dg Lie algebra introduced in Theorem \ref{algLie}.

\subsection{Twisted differential}

Let $\alpha : \C \rightarrow \Po$ be an operadic twisting morphism and let $\varphi : C \rightarrow A$ be an algebraic twisting morphism with respect to $\alpha$. We associate to $\alpha$ and $\varphi$ a \emph{twisted differential} $\partial_{\alpha,\, \varphi}$, denoted simply by $\partial_{\varphi}$, on Hom$_{g\mathcal{M}od_{\K}}(C,\, M)$ by the formula
$$\partial_{\varphi} (g) := \partial (g) + \alpha[\varphi,\, g].$$

\begin{lem}
If $\alpha \in$ \emph{Tw($\C$, $\Po$)} and $\varphi \in$ \emph{Tw$_{\alpha}$(C, A)}, then ${\partial_{\varphi}}^2 = 0$.
\end{lem}

\begin{pf}
We recall that $|\alpha| = -1$ and $|\varphi| = 0$. Let us modify a little bit the operator $\alpha[\varphi,\, g]_{n}$. We define for all $\psi$ in Hom$_{g\mathcal{M}od_{\K}}(C,\, A)$ and $g$ in Hom$_{g\mathcal{M}od_{\K}}(C,\, M)$ the operator $\alpha[\varphi,\, (\psi,\, g)]_{n}$ to be the composite:
$$C \xrightarrow{\Delta_{C,\, n}} (\C(n) \otimes C^{\otimes n})^{\So_n} \xrightarrow{\sum_{j} \alpha \otimes \varphi^{j-1} \otimes \psi \otimes \varphi^{n-j-1} \otimes g} \Po(n) \otimes A^{\otimes n-1} \otimes M \epi \Po (A,\, M) \xrightarrow{\gamma_M} M.$$
We define $\alpha[\varphi,\, (\psi,\, g)] := \sum_{n\geq 2}\alpha[\varphi,\, (\psi,\, g)]_{n}$.

(We have to pay attention to the fact that sign $(-1)^{|\psi||g|}$ may appear. The elements of $\C(C)$ are invariant under the action of the symmetric groups, so they are of the form $\sum_{\sigma \in \So_{n}} \varepsilon_{\sigma} \mu^{c}\cdot \sigma \otimes c_{\sigma^{-1}(1)} \otimes \cdots \otimes c_{\sigma^{-1}(n)}$, where $\varepsilon_{\sigma}$ depends on $(-1)^{|c_{i}||c_{j}|}$. For example, $\varepsilon_{(12)} = (-1)^{|c_{1}||c_{2}|}$, $\varepsilon_{(123)} = (-1)^{|c_{1}||c_{3}|+|c_{2}||c_{3}|}$ and $\varepsilon_{(132)} = (-1)^{|c_{1}||c_{2}|+|c_{1}||c_{3}|}$. Moreover, the coinvariant elements in $\Po (A,\, M)$ satisfy $\mu \otimes_{\So_{n}} (a_{1} \otimes \cdots \otimes a_{n}) = \varepsilon_{\sigma} \mu \cdot \sigma \otimes_{\So_{n}} (a_{\sigma^{-1}(1)} \otimes \cdots \otimes a_{\sigma^{-1}(n)})$. The image of $(-1)^{|c_{1}||c_{2}|} \mu^{c} \cdot (12) \otimes c_{2} \otimes c_{1}$ under $\gamma_{M} \circ (\alpha \otimes \psi \otimes g)$ in $M$ is
$$\hspace{-6cm} (-1)^{|c_{1}||c_{2}|+|\mu^{c}|(|\psi|+|g|)+|g||c_{2}|}\alpha(\mu^{c})(\psi(c_{2}),\, g(c_{1}))$$
$$\hspace{3.9cm} = (-1)^{|c_{1}||c_{2}|+|\mu^{c}|(|\psi|+|g|)+|g||c_{2}|}(-1)^{|\psi(c_{2})||g(c_{1})|} \mu(g(c_{1}),\, \psi(c_{2}))$$
$$\hspace{2cm} = (-1)^{|\psi||g|}(-1)^{|\mu^{c}|(|\psi|+|g|)+|\psi||c_{1}|} \mu(g(c_{1}),\, \psi(c_{2})).$$
Therefore, the operator $\alpha[\varphi,\, (\psi,\, g)]$ can be understood as follows
$$\alpha[\varphi,\, (\psi,\, g)] = \sum \left( \vcenter{\xymatrix@R=6pt@C=6pt{\varphi \ar@{-}[drrr] & \varphi \ar@{-}[drr] & \psi \ar@{-}[dr] & \varphi \ar@{-}[d] & \varphi \ar@{-}[dl] & g \ar@{-}[dll] & \varphi \ar@{-}[dlll]\\
&&& \alpha \ar@{-}[d] &&&\\
&&&&&&}} + (-1)^{|\psi||g|}\vcenter{\xymatrix@R=6pt@C=6pt{\varphi \ar@{-}[drrr] & \varphi \ar@{-}[drr] & g \ar@{-}[dr] & \varphi \ar@{-}[d] & \varphi \ar@{-}[dl] & \psi \ar@{-}[dll] & \varphi \ar@{-}[dlll]\\
&&& \alpha \ar@{-}[d] &&&\\
&&&&&&}}\right).)$$

The maps $\Delta_{C}$ and $\gamma_{M}$ are maps of dg modules and we have the equality
$$\partial(\alpha \otimes \varphi^{\otimes n-1} \otimes \psi) = \partial(\alpha) \otimes \varphi^{\otimes n-1} \otimes \psi + (-1)^{|\alpha|}\alpha \otimes \partial(\varphi^{\otimes n-1}) \otimes \psi + (-1)^{|\alpha|} \alpha \otimes \varphi^{\otimes n-1} \otimes \partial(\psi),$$
where $\partial(\varphi^{\otimes n-1}) = \sum_{j} \varphi^{j-1} \otimes \partial(\varphi) \otimes \varphi^{n-j-1}$. Therefore we get
$$\partial (\alpha[\varphi,\, g]) = \partial(\alpha)[\varphi,\, g] + (-1)^{|\alpha|} \alpha[\varphi,\, (\partial(\varphi),\, g)] + (-1)^{|\alpha|} \alpha[\varphi,\, \partial(g)].$$
It follows that
$$\begin{array}{lcl}
{\partial_{\varphi}}^2 (g) & = & \partial_{\varphi} (\partial(g) + \alpha[\varphi,\, g])\\
& = & \partial^2 (g) + \partial(\alpha[\varphi,\, g]) + \alpha[\varphi,\, \partial(g)] + \alpha[\varphi,\, \alpha[\varphi,\, g]]\\
& = &  \partial(\alpha)[\varphi,\, g] + (-1)^{|\alpha|} \alpha[\varphi,\, (\partial(\varphi),\, g)] + (-1)^{|\alpha|} \alpha[\varphi,\, \partial(g)]\\
&& + \alpha[\varphi,\, \partial(g)] + \alpha[\varphi,\, \alpha[\varphi,\, g]]\\
& = & \partial(\alpha)[\varphi,\, g] - \alpha[\varphi,\, (\partial(\varphi),\, g)] + \alpha[\varphi,\, \alpha[\varphi,\, g]].
\end{array}$$
The following picture
$$\begin{array}{lcl}
\sum \vcenter{\xymatrix@M=2pt@R=3pt@C=2pt{&& \varphi \ar@{-}[dr] & g \ar@{-}[d] & \varphi \ar@{-}[dl] &&\\
\varphi \ar@{-}[drrr] & \varphi \ar@{-}[drr] && \alpha \ar@{-}[d] && \varphi \ar@{-}[dll] & \varphi \ar@{-}[dlll]\\
&&& \alpha \ar@{-}[d] &&&\\
&&&&&&}} & \hspace{-0.3cm} = & \hspace{-0.2cm} \sum \left( \vcenter{\xymatrix@M=2pt@R=4pt@C=2pt{\varphi & g & \varphi \ar@{-}[dr] & \varphi \ar@{-}[d] & \varphi \ar@{-}[dl] & \varphi & \varphi\\
\ar@{-}[drrr] & \ar@{-}[drr] && \alpha \ar@{-}[d] && \ar@{-}[dll] & \ar@{-}[dlll]\\
&&& \alpha \ar@{-}[d] &&&\\
&&&&&&}} + \vcenter{\xymatrix@M=2pt@R=4pt@C=2pt{\varphi & \varphi & \varphi \ar@{-}[dr] & g \ar@{-}[d] & \varphi \ar@{-}[dl] & \varphi & \varphi\\
\ar@{-}[drrr] & \ar@{-}[drr] && \alpha \ar@{-}[d] && \ar@{-}[dll] & \ar@{-}[dlll]\\
&&& \alpha \ar@{-}[d] &&&\\
&&&&&&}} + \vcenter{\xymatrix@M=2pt@R=4pt@C=2pt{\varphi & \varphi & \varphi \ar@{-}[dr] & \varphi \ar@{-}[d] & \varphi \ar@{-}[dl] & g & \varphi\\
\ar@{-}[drrr] & \ar@{-}[drr] && \alpha \ar@{-}[d] && \ar@{-}[dll] & \ar@{-}[dlll]\\
&&& \alpha \ar@{-}[d] &&&\\
&&&&&&}}\right)\\
& \hspace{-0.3cm} - & \hspace{-0.2cm} \sum \left((-1)^{|g|}\vcenter{\xymatrix@M=2pt@R=3pt@C=2pt{&& \varphi \ar@{-}[dr] & \varphi \ar@{-}[d] & \varphi \ar@{-}[dl] &&\\
\varphi \ar@{-}[drrr] & g \ar@{-}[drr] && \alpha \ar@{-}[d] && \varphi \ar@{-}[dll] & \varphi \ar@{-}[dlll]\\
&&& \alpha \ar@{-}[d] &&&\\
&&&&&&}} + \vcenter{\xymatrix@M=2pt@R=3pt@C=2pt{&& \varphi \ar@{-}[dr] & \varphi \ar@{-}[d] & \varphi \ar@{-}[dl] &&\\
\varphi \ar@{-}[drrr] & \varphi \ar@{-}[drr] && \alpha \ar@{-}[d] && g \ar@{-}[dll] & \varphi \ar@{-}[dlll]\\
&&& \alpha \ar@{-}[d] &&&\\
&&&&&&}} \right)
\end{array}$$
models the equation
$$\alpha[\varphi,\, \alpha[\varphi,\, g]] = (\alpha \star \alpha)[\varphi,\, g] - \alpha[\varphi,\, (\star_{\alpha}(\varphi),\, g)]$$
(the sign $(-1)^{|g|}$ appears when we permute $\alpha$ and $g$). Thus
$${\partial_{\varphi}}^2 (g) = (\partial(\alpha) + \alpha \star \alpha)[\varphi,\, g] - \alpha[\varphi,\, (\partial(\varphi) + \star_{\alpha} (\varphi),\, g)].$$
Since $\alpha$ is an operadic twisting morphism and $\varphi$ is an algebraic twisting morphism with respect to $\alpha$, this concludes the proof.
$\cqfd$
\end{pf}

\subsection{The cotangent complex of an algebra over an operad}\label{generalcase}

From now on, we trace through the isomorphisms of Theorem \ref{cotangentcomplex} in order to make the differential on the cotangent complex explicit. Finally, for appropriate differentials, we obtain the isomorphism of differential graded modules
$$\textrm{Der}_{A}(\Po (C),\, M) \cong \textrm{Hom}_{\M_{A}^{\Po}}(A\otimes^{\Po} C,\, M),$$
where $\Po (C)$ is a quasi-free resolution of $A$.

We have in mind the resolutions obtained by means of the augmented bar construction on the level of operad, applied to an algebra, or the Koszul complex on an algebra or the coaugmented cobar construction on the level of cooperads, applied to a homotopy algebra.

The space Der$_{A} (\Po (C),\, M)$ is endowed with the following differential
$$\partial(f) = d_{M} \circ f - (-1)^{|f|} f \circ d_{\alpha},$$
where $d_{\alpha}$ was defined in Section \ref{quasifreeresolution}.

\begin{prop}\label{iso1}
With the above notations, we have the following isomorphism of dg modules
$$(\emph{Der}_{A}(\Po \circ_{\alpha} C,\, M),\, \partial) \cong (\emph{Hom}_{g\mathcal{M}od_{\K}}(C,\, M),\, \partial_{\varphi} = \partial + \alpha[\varphi,\, -]), \textrm{ where } C = \C(A).$$
\end{prop}

\begin{pf}
First, the isomorphism of $\K$-modules between Der$_{A}(\Po (C),\, M)$ and Hom$_{g\mathcal{M}od_{\K}}(C,\, M)$ is given by the restriction on the generators $C$.

We verify that this isomorphism commutes with the respective differentials. We fix the notations $\bar{f} := f_{|C}$ and $n := |\bar{f}| = |f|$. On the one hand, we have
$$\begin{array}{ll}
\partial(f)_{|C} & = (d_{M} \circ f)_{|C} - (-1)^{|f|} (f\circ d_{\alpha})_{|C}\\
& = d_{M} \circ \bar{f} - (-1)^{n} f\circ (d_{\Po} \circ id_{C} + id_{\Po} \circ' d_C - \delta_{\alpha}^l)_{|C}.
\end{array}$$
Moreover, $(d_{\Po} \circ id_{C})_{|C} = 0$ since $(d_{\Po})_{|\Po(1)} = 0$ and $f\circ (id_{\Po} \circ' d_C)_{|C} = \bar{f} \circ d_C$. Thus
$$\partial(f)_{|C} = d_{M} \circ \bar{f} - (-1)^{n} \bar{f} \circ d_C +(-1)^n f\circ {\delta_{\alpha}^l}_{|C}.$$
On the other hand,
$$\partial_{\varphi} (\bar{f}) = d_{M} \circ \bar{f} - (-1)^{n} \bar{f} \circ d_C + \alpha[\varphi,\, \bar{f}].$$
With the signs $\alpha \otimes \bar{f} = (-1)^{|\alpha| |\bar{f}|} (id\otimes \bar{f})\otimes (\alpha \otimes id)$ and using the fact that $f$ is a derivation, we verify that $(-1)^{n}f\circ {\delta_{\alpha}^l}_{|C} = \alpha[\varphi,\, \bar{f}]$.
$\cqfd$
\end{pf}

Let us construct a twisted differential on the free $A$-module $A\otimes^{\Po}C$ as follows. Since $A\otimes^{\Po}C$ is a quotient of $\Po (A,\, C)$, we define a map
$$\delta_{1}^{l}(n) : \Po (A,\, C) \xrightarrow{id_{\Po} (id_{A},\, \Delta_{C}(n))} \Po (A,\, (\C(n)\otimes C^{\otimes n})^{\So_{n}}) \xrightarrow{id_{\Po} (id_{A},\, \alpha \otimes \varphi^{\otimes n-1} \otimes id_{C})}$$
$$\hspace{2cm}\Po (A,\, \Po(n) \otimes A^{\otimes n-1} \otimes C) \rightarrow (\Po \circ \Po) (A,\, C) \xrightarrow{\gamma (id_{A},\, id_{C})} \Po (A,\, C).$$
This map sends the elements $\mu \otimes \gamma_{A}(\nu_{1}\otimes a_{1}\otimes \cdots \otimes a_{i_{1}})\otimes \cdots \otimes c \otimes \cdots \otimes \gamma_{A}(\nu_{k}\otimes \cdots \otimes a_{n})$ and $\gamma_{\Po}(\mu \otimes \nu_{1} \otimes \cdots \otimes \nu_{k})\otimes a_{1}\otimes \cdots \otimes a_{i_{1}}\otimes \cdots \otimes c \otimes \cdots \otimes a_{n}$ to the same image, for $c\in C$ and $a_{j} \in A$ and $\mu,\, \nu_{j} \in \Po$. So $\delta_{1}^{l}(n)$ induces a map on the quotient
$$\delta_{\alpha,\, \varphi}^{l}(n) : A\otimes^{\Po} C \rightarrow A\otimes^{\Po} C.$$
We write $\delta_{1}^{l} := \sum \delta_{1}^{l}(n)$ and $\delta_{\alpha,\, \varphi}^{l} := \sum \delta_{\alpha,\, \varphi}^{l}(n)$, or simply $\delta_{\varphi}^{l}$.

We define the twisted differential $\partial_{\alpha,\, \varphi}$, or simply $\partial_{\varphi}$ on Hom$_{\M_{A}^{\Po}}(A\otimes^{\Po}C,\, M)$ by
$$\begin{array}{rl}
\partial_{\varphi} (f) := & \partial(f) + (-1)^{|f|} f \circ \delta_{\varphi}^{l}\\
= & d_{M} \circ f - (-1)^{|f|}f\circ (d_{A\otimes^{\Po}C} - \delta_{\varphi}^{l}),
\end{array}$$
where the differential $d_{A\otimes^{\Po} C}$ is induced by the natural differential on $\Po (A,\, C)$. 
So we consider the twisted differential $d_{\varphi} := d_{A\otimes^{\Po}C} - \delta_{\varphi}^{l}$ on $A\otimes^{\Po}C$. Once again, the notation $\partial_{\varphi}$ stands for several differentials and the relevant one is given without ambiguity by the context.

\begin{thm}\label{iso2}
With the above notations, the following three dg modules are isomorphic
$$(\emph{Der}_{A}(\Po \circ_{\alpha} C,\, M),\, \partial) \cong (\emph{Hom}_{g\mathcal{M}od_{\K}}(C,\, M),\, \partial_{\varphi}) \cong (\emph{Hom}_{\M_{A}^{\Po}}(A\otimes^{\Po} C,\, M),\, \partial_{\varphi}).$$
\end{thm}

\begin{pf}
We already know the isomorphism of $\K$-modules given by the restriction
$$(\textrm{Hom}_{\M_{A}^{\Po}}(A\otimes^{\Po} C,\, M),\, \partial) \cong (\textrm{Hom}_{g\mathcal{M}od_{\K}}(C,\, M),\, \partial)$$
from the preliminaries. We now verify that this isomorphism commutes with the differentials. With the notation $\bar{f} := f_{|C}$, we have
$$\partial_{\varphi} (\bar{f}) = d_{M} \circ \bar{f} - (-1)^{|\bar{f}|} \bar{f} \circ d_C + \alpha[\varphi,\, \bar{f}]$$
and
$$\partial_{\varphi}(f)_{|C} = (d_{M}\circ f - (-1)^{|f|}f\circ d_{A\otimes^{\Po}C} + (-1)^{|f|}f\circ \delta_{\varphi}^{l})_{|C}.$$
Since $(f\circ d_{A\otimes^{\Po}C})_{|C} = \bar{f} \circ d_C$, we just need to show the equality $\alpha[\varphi,\, \bar{f}] = (-1)^{|f|}(f \circ \delta_{\varphi}^{l})_{|C}$. This holds since $M \in \M_{A}^{\Po}$ and $f$ is a morphism of $A$-modules over $\Po$ and the structure of $A$-module on $C$ into $A\otimes^{\Po} C$ is just the projection $\Po (A,\, C) \epi A\otimes^{\Po} C$.
$\cqfd$
\end{pf}

Finally, when $\Po \circ_{\alpha} C \qiso A$ is a quasi-free resolution of $A$, the chain complex
$$(A\otimes^{\Po}C,\, d_{\varphi} = d_{A\otimes^{\Po}C} - \delta_{\varphi}^{l})$$
is a representation of the cotangent complex. In our cases, we have $C = \C(A)$. Then a representation of the cotangent complex is given by
$$(A\otimes^{\Po}\C(A),\, d_{\varphi} = d_{A\otimes^{\Po}\C(A)} - \delta_{\varphi}^{l} + \delta_{\varphi}^{r}),$$
where $\delta_{\varphi}^{l}$ is induced by
$$\hspace{-1cm} \Po (A,\, \C(A)) \xrightarrow{id_{\Po} \circ (id_{A},\, \Delta_{p} \circ id_{A})} \Po (A,\, (\C \circ_{(1)} \C) (A)) \xrightarrow{id_{\Po} \circ (id_{A},\, \alpha \circ_{(1)} id_{\C} \circ id_{A})}$$
$$\hspace{1cm} \Po (A,\, (\Po \circ_{(1)} \C) (A)) \mono (\Po \circ \Po)(A,\, \C(A)) \xrightarrow{\gamma \circ (id_{A},\, id_{\C(A)})} \Po (A,\, \C(A))$$
and $\delta_{\varphi}^{r}$ is induced by
$$\hspace{-1cm} \Po (A,\, \C(A)) \xrightarrow{id_{\Po} \circ (id_{A},\, \Delta_{p} \circ id_{A})} \Po (A,\, (\C \circ_{(1)} \C) (A)) \xrightarrow{id_{\Po} \circ (id_{A},\, id_{\C} \circ_{(1)} \alpha \circ id_{A})}$$
$$\hspace{2.5cm} \Po (A,\, (\C \circ_{(1)} \Po) (A)) \mono \Po (A,\, \C (A,\, \Po (A))) \xrightarrow{id_{\Po} \circ (id_{A},\, id_{\C} \circ (id_{A},\, \gamma_{A}))} \Po (A,\, \C(A)).$$

\begin{rem}
Applying this description to the resolutions of algebras obtained by means of the augmented bar construction or by means of the Koszul complex, we obtain two different chain complexes which allow us to compute the André-Quillen cohomology. The one using the Koszul resolution is smaller since $\Po^{\ash} \mono \B(\Po)$. However the differential on the one using the augmented bar construction is simpler as the differential strongly depends on the coproduct. The cooperad $\Po^{\ash}$ is often given up to isomorphism, therefore it is difficult to make it explicit.
\end{rem}

\section{Applications and new examples of cohomology theories}

We apply the previous general definitions to nearly all the operads we know. We explain which resolution can be used each time. Sometimes, it corresponds to known chain complexes. We also show that the cotriple cohomology corresponds to André-Quillen cohomology. Among the new examples, we make the André-Quillen cohomology for algebras over the operad $\Perm$ explicit. We do the same for homotopy $\Po$-algebras. From now on, we assume that the algebras are non-negatively graded.

\subsection{Applications}

For some operads, an explicit chain complex computing the cohomology theory for the associated algebras has already been proposed by various authors.

\begin{itemize}
\item When $\Po = \Ass$ is the operad of associative algebras, $A\otimes^{\Po}\Po^{\ash}(A) \cong (\K \oplus A) \otimes B(A)\otimes (\K \oplus A)$ (by \ref{adjunction1}) is the normalized Hochschild complex (see Section $1.1.14$ of \cite{Loday}). The André-Quillen cohomology of associative algebras is the Hochschild cohomology (see also Chapter IX, Section $6$ of \cite{CartanEilenberg}).
\item When $\Po = \Lie$ is the operad of Lie algebras, $A\otimes^{\Po}\Po^{\ash}(A) \cong U^{e}(A)\otimes \Lambda (A)$ since $\Lie^{\ash}(A) \cong \Lambda (A)$. The André-Quillen cohomology of Lie algebras is Chevalley-Eilenberg cohomology (see Chapter XIII of \cite{CartanEilenberg}).
\item When $\Po = \Com$ is the operad of commutative algebras, the complex $A\otimes^{\Po}\Po^{\ash}(A) \cong (\K \oplus A)\otimes \Com^{\ash}(A)$, only valid in characteristic $0$, gives the cohomology theory of commutative algebras defined by Quillen in \cite{Quillen}. It corresponds to Harrison cohomology defined in \cite{Harrison}. We refer to \cite{Loday} for the relationship between the different definitions.
\item When $\Po = \Dias$ is the operad of diassociative algebras and with the Koszul resolution, we get the chain complex and the associated cohomology defined by Frabetti in \cite{Frabetti}.
\item When $\Po = \Leib$ is the operad of Leibniz algebras and with the Koszul resolution, the André-Quillen cohomology of Leibniz algebras is the cohomology defined by Loday and Pirashvili in \cite{LodayPirashvili}.
\item For the operad $\Po = \Poiss$ encoding Poisson algebras, Fresse followed, as in this paper, the ideas of Quillen to make a cohomology of Poisson algebras explicit with the Koszul resolution \cite{Fresse2}.
\item When $\Po = \Prelie$ and with the Koszul resolution, the André-Quillen cohomology of pre-Lie algebras is the one defined by Dzhumadil'daev in \cite{Dzhumadil}.
\item When $\Po = \Zinb$, or equivalently $\Leib^{!}$, and with the Koszul resolution, the André-Quillen cohomology of Zinbiel algebras is the one given in \cite{Balavoine}.
\end{itemize}

More generally,
\begin{itemize}
\item Balvoine introduces a chain complex in \cite{Balavoine}. When the operad $\Po$ is a binary Koszul operad, the chain complex computing the André-Quillen cohomology obtained with the Koszul resolution corresponds to the one defined by Balavoine. Thus, the cohomology theories are the same in this case.
\end{itemize}

\subsection{The case of Perm algebras}\label{Perm}

We denote by $\Perm$ the operad corresponding to Perm algebras defined in \cite{Chapoton}.\\

Let us recall that a basis for $\Perm(n)$ is given by corollas in space with $n$ leaves labelled by $1$ to $n$ with one leaf underlined. So $\Perm(n)$ is of dimension $n$. The composition product in $\Perm$ is given by the path traced through the upper underlined leaf from the root. For example, $\def\objectstyle{\scriptstyle}
\def\labelstyle{\scriptstyle}
\gamma \left(\vcenter{\xymatrix@M=0pt@R=4pt@C=4pt{& \underline{1} & 2 & 3 &&&\\
& \ar@{-}[dr] & \ar@{-}[d] & \ar@{-}[dl] &&&\\
\ar@{.}[rrrrrr] &&&&&&\\
&& \underline{1} && 2 & 3 &\\
&& \ar@{-}[drr] && \ar@{-}[d] & \ar@{-}[dl] &\\
&&&& \ar@{-}[d] &&\\
&&&&&&}}\right) = \vcenter{\xymatrix@M=0pt@R=4pt@C=4pt{& \underline{1} & 2 & 3 & 4 & 5 &\\
& \ar@{-}[drr] & \ar@{-}[dr] & \ar@{-}[d] & \ar@{-}[dl] & \ar@{-}[dll] &\\
&&& \ar@{-}[d] &&&\\
&&&&&&}}$ and $\def\objectstyle{\scriptstyle}
\def\labelstyle{\scriptstyle}
\gamma \left(\vcenter{\xymatrix@M=0pt@R=4pt@C=4pt{&& \underline{1} & 2 & 3 &&\\
&& \ar@{-}[dr] & \ar@{-}[d] & \ar@{-}[dl] &&\\
\ar@{.}[rrrrrr] &&&&&&\\
& \underline{1} && 2 && 3 &\\
& \ar@{-}[drr] && \ar@{-}[d] && \ar@{-}[dll] &\\
&&& \ar@{-}[d] &&&\\
&&&&&&}}\right) = \vcenter{\xymatrix@M=0pt@R=4pt@C=4pt{& \underline{1} & 2 & 3 & 4 & 5 &\\
& \ar@{-}[drr] & \ar@{-}[dr] & \ar@{-}[d] & \ar@{-}[dl] & \ar@{-}[dll] &\\
&&& \ar@{-}[d] &&&\\
&&&&&&}}$.

 In \cite{ChapotonLivernet}, the authors show that the Koszul dual operad of the operad $\Perm$ is the operad $\Prelie$ and that the operad $\Prelie$ is Koszul. It follows that the operad $\Perm$ is Koszul (see \cite{GinzburgKapranov} for  general facts about Koszul duality of operads). Since $\Perm^{\ash} \cong \Prelie^{\vee}$, it is possible to understand the coproduct on $\Perm^{\ash}$ if we know the product on $\Prelie$. Chapoton and Livernet gave an explicit basis for $\Prelie$ and made explicit the product. This basis of $\Prelie$ is given by the rooted trees of degree $n$, that is with $n$ vertices, denoted $\RT(n)$. Then we need to understand the coproduct on $\Prelie^*$ which is given by
$$\Delta : \Prelie^* \xrightarrow{^t\gamma} (\Prelie \circ \Prelie)^* \iso \Prelie^* \circ \Prelie^*,$$
where $\Prelie^*(n) := \Prelie(n)^*$ and $^t\gamma (f) := f\circ \gamma$. A rooted tree is represented as in \cite{ChapotonLivernet}, with its root at the bottom. We make explicit the coproduct on a particular element
$$\def\objectstyle{\scriptstyle}
\def\labelstyle{\scriptstyle}
\Delta\left(\vcenter{\xymatrix@R=4pt@C=2pt{*+[o][F-]{1} \ar@{-}[dr] && 
*+[o][F-]{3} \ar@{-}[dl]\\
& *+[o][F-]{2}}}\right) = \vcenter{\xymatrix@R=4pt@C=4pt{*+[o][F-]{1}}} \circ_{1} \vcenter{\xymatrix@R=4pt@C=2pt{*+[o][F-]{1} \ar@{-}[dr] && 
*+[o][F-]{3} \ar@{-}[dl]\\
& *+[o][F-]{2}}} + \vcenter{\xymatrix@R=4pt@C=4pt{*+[o][F-]{1} \ar@{-}[d]\\
*+[o][F-]{2}}} \circ_{2} \vcenter{\xymatrix@R=4pt@C=4pt{*+[o][F-]{2} \ar@{-}[d]\\
*+[o][F-]{1}}} + \vcenter{\xymatrix@R=4pt@C=4pt{*+[o][F-]{2} \ar@{-}[d]\\
*+[o][F-]{1}}} \circ_{1} \vcenter{\xymatrix@R=4pt@C=4pt{*+[o][F-]{1} \ar@{-}[d]\\
*+[o][F-]{2}}} + \vcenter{\xymatrix@R=4pt@C=2pt{*+[o][F-]{1} \ar@{-}[dr] && 
*+[o][F-]{3} \ar@{-}[dl]\\
& *+[o][F-]{2}}} \circ_{1} \vcenter{\xymatrix@R=4pt@C=4pt{*+[o][F-]{1}}} + \vcenter{\xymatrix@R=4pt@C=2pt{*+[o][F-]{1} \ar@{-}[dr] && 
*+[o][F-]{3} \ar@{-}[dl]\\
& *+[o][F-]{2}}} \circ_{2} \vcenter{\xymatrix@R=4pt@C=4pt{*+[o][F-]{1}}} + \vcenter{\xymatrix@R=4pt@C=2pt{*+[o][F-]{1} \ar@{-}[dr] && 
*+[o][F-]{3} \ar@{-}[dl]\\
& *+[o][F-]{2}}} \circ_{3} \vcenter{\xymatrix@R=4pt@C=4pt{*+[o][F-]{1}}}.$$

Let $A$ be a $\Perm$-algebra. The cotangent complex has the following form
$$\begin{array}{lcl}
A\otimes^{\Po} \Poa (A) & = & A\otimes^{\Po} \RT (A) \cong \vcenter{
\xymatrix@M=0pt@R=6pt@C=4pt{\underline{\RT(A)}\\
\ar@{-}[dd]\\ &\\ &}} \oplus \vcenter{\xy (3,0)*{} ; (3,2)*{} **\dir{-},
(3,2)*{} ; (0,4)*{} **\dir{-},
(3,2)*{} ; (6,4)*{} **\dir{-},
(-1,7.4)*{\underline{A}} ; (8,7)*{\RT(A)}
\endxy} \oplus \ \vcenter{\xy (5,0)*{} ; (5,2)*{} **\dir{-},
(5,2)*{} ; (2,4)*{} **\dir{-},
(5,2)*{} ; (8,4)*{} **\dir{-},
(-1,7)*{\underline{\RT(A)}} ; (8,7.4)*{A}
\endxy}\\
& \cong & \RT(A) \oplus A\otimes \RT(A) \oplus \RT(A) \otimes A,
\end{array}$$
where $\RT(A) = \oplus_n \RT(n) \otimes_{\So_n} A^{\otimes n}$.

\subsubsection{\bf When the algebra is trivial}

We assume first that $A$ is a trivial algebra, that is $\gamma_A \equiv 0$. To make the differential on the cotangent complex explicit, we just need to describe the restriction $\RT(A) \rightarrow \RT(A) \otimes A \oplus A\otimes \RT(A)$ since the differential is zero on $A \otimes \RT(A) \oplus \RT(A) \otimes A$. Let $T$ be in $\RT(n)$. There are several possibilities:
\begin{itemize}
 \item[i)] the rooted tree $T$ has the form $\def\objectstyle{\scriptstyle}
\def\labelstyle{\scriptstyle}
\vcenter{\xymatrix@R=4pt@C=4pt{*++[o][F-]{T_1} \ar@{-}[d]\\
*++[o][F-]{1}}}$, where $T_1$ is in $\RT(n-1)$. In that case, the term $\def\objectstyle{\scriptstyle}
\def\labelstyle{\scriptstyle}
\vcenter{\xymatrix@R=4pt@C=4pt{*+[o][F-]{2} \ar@{-}[d]\\
*+[o][F-]{1}}} \circ_2 T_1$ appears in $\Delta_{\RT} (T)$, so the image of $T\otimes a_1 \otimes \cdots \otimes a_n$ under $d_{\varphi}$ in $A\otimes \RT (A)$ contains $-a_1 \otimes (T_1 \otimes a_2 \otimes \cdots \otimes a_n)$;
\item[ii)] there exists $T_2$ in $\RT(n-1)$ such that the rooted tree $T$ can be written $\def\objectstyle{\scriptstyle}
\def\labelstyle{\scriptstyle}
\vcenter{\xymatrix@R=4pt@C=4pt{*++[o][F-]{1} \ar@{-}[d]\\
*++[o][F-]{T_2}}}$. In that case, the term $\def\objectstyle{\scriptstyle}
\def\labelstyle{\scriptstyle}
\vcenter{\xymatrix@R=4pt@C=4pt{*+[o][F-]{1} \ar@{-}[d]\\
*+[o][F-]{2}}} \circ_2 T_2$ appears in $\Delta_{\RT} (T)$, so the image of $T\otimes a_1 \otimes \cdots \otimes a_n$ under $d_{\varphi}$ in $\RT (A) \otimes A$ contains $-(T_2 \otimes a_2 \otimes \cdots \otimes a_n) \otimes a_1$;
\item[iii)] the rooted tree has the form $\def\objectstyle{\scriptstyle}
\def\labelstyle{\scriptstyle}
\vcenter{\xymatrix@R=4pt@C=4pt{*++[o][F-]{T_3} \ar@{-}[d]\\
*++[o][F-]{n}}}$, where $T_3$ is in $\RT(n-1)$. In that case, the term $\def\objectstyle{\scriptstyle}
\def\labelstyle{\scriptstyle}
\vcenter{\xymatrix@R=4pt@C=4pt{*+[o][F-]{1} \ar@{-}[d]\\
*+[o][F-]{2}}} \circ_1 T_3$ appears in $\Delta_{\RT} (T)$, so the image of $T\otimes a_1 \otimes \cdots \otimes a_n$ under $d_{\varphi}$ in $A\otimes \RT (A)$ contains $-a_n \otimes (T_3 \otimes a_1 \otimes \cdots \otimes a_{n-1})$;
\item[iv)] there exists $T_4$ in $\RT(n-1)$ such that the rooted tree can be written \(\def\objectstyle{\scriptstyle}
\def\labelstyle{\scriptstyle}
\vcenter{\xymatrix@R=4pt@C=4pt{*++[o][F-]{n} \ar@{-}[d]\\
*++[o][F-]{T_4}}}\). In that case, the term $\def\objectstyle{\scriptstyle}
\def\labelstyle{\scriptstyle}
\vcenter{\xymatrix@R=4pt@C=4pt{*+[o][F-]{2} \ar@{-}[d]\\
*+[o][F-]{1}}} \circ_1 T_4$ appears in $\Delta_{\RT} (T)$, so the image of $T\otimes a_1 \otimes \cdots \otimes a_n$ under $d_{\varphi}$ in $\RT (A) \otimes A$ contains $-(T_4 \otimes a_1 \otimes \cdots \otimes a_{n-1}) \otimes a_n$;
\end{itemize}
A rooted tree $T$ has the shape i) and iv), or ii) and iii), or ii) and iv), or i) only, or ii) only, or iii) only, or iv) only, or finally a shape not described in i) to iv). In this last case, the differential is $0$. Otherwise, the image under the differential of an element $T\otimes a_1 \otimes \cdots \otimes a_n$ in $\RT(A)$ is given by the sum of the corresponding terms in i) to iv). For example, if $T$ can be written $\def\objectstyle{\scriptstyle}
\def\labelstyle{\scriptstyle}
\vcenter{\xymatrix@R=4pt@C=4pt{*++[o][F-]{T_1} \ar@{-}[d]\\
*++[o][F-]{1}}}$ and \(\def\objectstyle{\scriptstyle}
\def\labelstyle{\scriptstyle}
\vcenter{\xymatrix@R=4pt@C=4pt{*++[o][F-]{n} \ar@{-}[d]\\
*++[o][F-]{T_4}}}\), we get $d_{\varphi}(T\otimes a_1 \otimes \cdots \otimes a_n) = -a_1 \otimes (T_1 \otimes a_2 \otimes \cdots \otimes a_n) -(T_4 \otimes a_1 \otimes \cdots \otimes a_{n-1}) \otimes a_n$.

\subsubsection{\bf For any Perm algebra}

For a general $\Perm$-algebra $A$, we no longer assume a priori that the restriction of the differential $d_{A\otimes^{\Po}\Poa (A)}$ to $\Poa(A)$, that is $d_{\alpha}$, is zero. For a rooted tree $T$ in $\RT(n)$, we define the function $f$ by $f(T, i, j) = 1$ if $T = \def\objectstyle{\scriptstyle}
\def\labelstyle{\scriptstyle}
\vcenter{\xymatrix@R=4pt@C=4pt{*++[o][F-]{T_{2}} \ar@{-}[d] & *++[o][F-]{T_{3}} \ar@{-}[d]\\
*++[o][F-]{i} \ar@{-}[d] & *++[o][F-]{j} \ar@{-}[l]\\
*++[o][F-]{T_1} &}}$ for some rooted tree $T_{1}$ and some families of rooted trees $T_{2}$ and $T_{3}$, and $f(T, i, j) = 0$ otherwise. There exists a rooted tree $T_{i}$ in $\RT(n-1)$ such that $T$ appears in the product $T_{i} \circ_{i} \def\objectstyle{\scriptstyle}
\def\labelstyle{\scriptstyle}
\vcenter{\xymatrix@R=4pt@C=4pt{*+[o][F-]{2} \ar@{-}[d]\\
*+[o][F-]{1}}}$ if and only if $f(T, i, i+1) = 1$ (take $T_{i} = \def\objectstyle{\scriptstyle}
\def\labelstyle{\scriptstyle}
\vcenter{\xymatrix@R=4pt@C=4pt{*+[o][F-]{T'_{2}} \ar@{-}[d] & *+[o][F-]{T'_{3}} \ar@{-}[dl]\\
*++[o][F-]{i} \ar@{-}[d] &\\
*+[o][F-]{T'_1} &}}$ where $T'_{j}$ is the family of trees $T_{j}$ with vertices $k > i$ replaced by $k+1$). Similarly there exists a rooted tree $T_{i}$ in $\RT(n-1)$ such that $T$ appears in the product $T_{i} \circ_{i} \def\objectstyle{\scriptstyle}
\def\labelstyle{\scriptstyle}
\vcenter{\xymatrix@R=4pt@C=4pt{*+[o][F-]{1} \ar@{-}[d]\\
*+[o][F-]{2}}}$ if and only if $f(T, i+1, i) = 1$. We define $E^1(T) := \{i\, |\, f(T, i, i+1) = 1\}$ and $E^2(T) := \{i\, |\, f(T, i+1, i) = 1\}$. We obtain
$$\begin{array}{lcl}
d_{\alpha}(T\otimes a_{1} \otimes \cdots \otimes a_{n}) & = & \sum_{i\in E^1(T)} T_{i} \otimes a_{1} \otimes \cdots \otimes \gamma_{A}(\vcenter{\xy (2,0)*{} ; (2,2)*{} **\dir{-},
(2,2)*{} ; (0,4)*{} **\dir{-},
(2,2)*{} ; (4,4)*{} **\dir{-},
(0,5)*{-}
\endxy} \otimes a_{i} \otimes a_{i+1}) \otimes \cdots \otimes a_{n}\\
& + & \sum_{i\in E^2(T)} T_{i} \otimes a_{1} \otimes \cdots \otimes \gamma_{A}(\vcenter{\xy (2,0)*{} ; (2,2)*{} **\dir{-},
(2,2)*{} ; (0,4)*{} **\dir{-},
(2,2)*{} ; (4,4)*{} **\dir{-},
(4,5)*{-}
\endxy} \otimes a_{i} \otimes a_{i+1}) \otimes \cdots \otimes a_{n},
\end{array}$$
where $T_{i}$ is the rooted tree such that $T$ appears in the product $T_{i} \circ_{i} \def\objectstyle{\scriptstyle}
\def\labelstyle{\scriptstyle}
\vcenter{\xymatrix@R=4pt@C=4pt{*+[o][F-]{2} \ar@{-}[d]\\
*+[o][F-]{1}}}$ or $T_{i} \circ_{i} \def\objectstyle{\scriptstyle}
\def\labelstyle{\scriptstyle}
\vcenter{\xymatrix@R=4pt@C=4pt{*+[o][F-]{1} \ar@{-}[d]\\
*+[o][F-]{2}}}$. Finally, on $\RT(A)$, the differential on the cotangent complex is given by $d_{\varphi} = d_{\alpha} - \delta_{\varphi}^l$.

We describe now the differential $\delta_{\varphi}^l$ on $A\otimes \RT(A)$ thanks to the description i) - iv) of the previous section.
\begin{itemize}
\item[i)-ii)] The term $\gamma_{A}(\vcenter{\xy (2,0)*{} ; (2,2)*{} **\dir{-},
(2,2)*{} ; (0,4)*{} **\dir{-},
(2,2)*{} ; (4,4)*{} **\dir{-},
(0,5)*{-}
\endxy} \otimes a_{0} \otimes a_{1}) \otimes (T_{i} \otimes a_{2} \otimes \cdots \otimes a_{n})$ appears in $\delta_{\varphi}^l(a_{0} \otimes (T \otimes a_{1} \otimes \cdots \otimes a_{n}))$ (with $i = 1$ or $2$);
\item[iii)-iv)] the term $\gamma_{A}(\vcenter{\xy (2,0)*{} ; (2,2)*{} **\dir{-},
(2,2)*{} ; (0,4)*{} **\dir{-},
(2,2)*{} ; (4,4)*{} **\dir{-},
(0,5)*{-}
\endxy} \otimes a_{0} \otimes a_{n}) \otimes (T_{i} \otimes a_{1} \otimes \cdots \otimes a_{n-1})$ appears in $\delta_{\varphi}^l(a_{0} \otimes (T \otimes a_{1} \otimes \cdots \otimes a_{n}))$ (with $i = 3$ or $4$).
\end{itemize}
Similarly, we describe the differential $\delta_{\varphi}^l$ on $\RT(A) \otimes A$.
\begin{itemize}
\item[i)] The term $\gamma_{A}(\vcenter{\xy (2,0)*{} ; (2,2)*{} **\dir{-},
(2,2)*{} ; (0,4)*{} **\dir{-},
(2,2)*{} ; (4,4)*{} **\dir{-},
(0,5)*{-}
\endxy} \otimes a_{1} \otimes a_{n+1}) \otimes (T_{1} \otimes a_{2} \otimes \cdots \otimes a_{n})$ appears in $\delta_{\varphi}^l((T \otimes a_{1} \otimes \cdots \otimes a_{n}) \otimes a_{n+1})$;
\item[ii)] the term $(T_{2} \otimes a_{2} \otimes \cdots \otimes a_{n}) \otimes \gamma_{A}(\vcenter{\xy (2,0)*{} ; (2,2)*{} **\dir{-},
(2,2)*{} ; (0,4)*{} **\dir{-},
(2,2)*{} ; (4,4)*{} **\dir{-},
(0,5)*{-}
\endxy} \otimes a_{1} \otimes a_{n+1})$ appears in $\delta_{\varphi}^l((T \otimes a_{1} \otimes \cdots \otimes a_{n}) \otimes a_{n+1})$;
\item[iii)] the term $\gamma_{A}(\vcenter{\xy (2,0)*{} ; (2,2)*{} **\dir{-},
(2,2)*{} ; (0,4)*{} **\dir{-},
(2,2)*{} ; (4,4)*{} **\dir{-},
(0,5)*{-}
\endxy} \otimes a_{n} \otimes a_{n+1}) \otimes (T_{3} \otimes a_{1} \otimes \cdots \otimes a_{n-1})$ appears in $\delta_{\varphi}^l((T \otimes a_{1} \otimes \cdots \otimes a_{n}) \otimes a_{n+1})$;
\item[iv)] the term $(T_{4} \otimes a_{1} \otimes \cdots \otimes a_{n-1}) \otimes \gamma_{A}(\vcenter{\xy (2,0)*{} ; (2,2)*{} **\dir{-},
(2,2)*{} ; (0,4)*{} **\dir{-},
(2,2)*{} ; (4,4)*{} **\dir{-},
(0,5)*{-}
\endxy} \otimes a_{n} \otimes a_{n+1})$ appears in $\delta_{\varphi}^l((T \otimes a_{1} \otimes \cdots \otimes a_{n}) \otimes a_{n+1})$.
\end{itemize}
Finally, the differential on the cotangent complex $\RT(A) \oplus A\otimes \RT(A) \oplus \RT(A) \otimes A$ is given by $d_{\alpha} + id_{A}\otimes d_{\alpha} + d_{\alpha}\otimes id_{A} - \delta_{\varphi}^l$.

\subsection{The case of $A_{\infty}$-algebras}\label{Ainfinity}

Markl gave in \cite{Markl} a definition for a cohomology theory for homotopy associative algebras. In this section, we make explicit the André-Quillen cohomology for homotopy associative algebras and we recover the complex defined by Markl.\\

The operad $A_{\infty} = \Omega(\Ass^{\ash}) = \F (\vcenter{\xymatrix@M=0pt@R=4pt@C=4pt{\ar@{-}[dr] && \ar@{-}[dl]\\
& \ar@{-}[d] &\\
&&}}, \vcenter{\xymatrix@M=0pt@R=4pt@C=4pt{\ar@{-}[dr] & \ar@{-}[dd] & \ar@{-}[dl]\\
&&\\
&&}}, \vcenter{\xymatrix@M=0pt@R=4pt@C=4pt{\ar@{-}[drr] & \ar@{-}[dr] && \ar@{-}[dl] & \ar@{-}[dll]\\
&& \ar@{-}[d] &&\\
&&&&}}, \ldots)$ is the free operad on one generator in each degree greater than $1$. We have the resolution $R := A_{\infty} \circ \Ass^{\ash} (A) \qiso A$ and we get
$$\Lb_{R|A} = \bigoplus_{\genfrac{}{}{0cm}{1}{l, l_1, l_2 \geq 0}{k \geq 0}} \bigoplus_{\genfrac{}{}{0cm}{1}{i_1 + \cdots + i_k = l_1}{j_1 + \cdots + j_k = l_2}} A^{\otimes i_1} | \cdots | A^{\otimes i_k} | A^{\otimes l} | A^{\otimes j_k} | \cdots | A^{\otimes j_1}.$$
Actually, an element in $\Lb_{R|A}$ should be seen as a planar tree
$$\def\objectstyle{\scriptstyle}
\def\labelstyle{\scriptstyle}
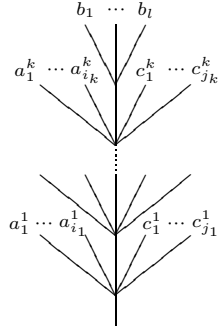
\xy <0.4cm, 0cm>:
(2,0)*{} ; (2,1)*{} **\dir{-},
(2,1)*{} ; (-0.5,3)*{} **\dir{-},
(2,1)*{} ; (1,3)*{} **\dir{-},
(2,1)*{} ; (2,3)*{} **\dir{-},
(2,1)*{} ; (3,3)*{} **\dir{-},
(2,1)*{} ; (4.5,3)*{} **\dir{-},
(2,3)*{} ; (-0.5,5)*{} **\dir{-},
(2,3)*{} ; (1,5)*{} **\dir{-},
(2,3)*{} ; (2,5)*{} **\dir{-},
(2,3)*{} ; (3,5)*{} **\dir{-},
(2,3)*{} ; (4.5,5)*{} **\dir{-},
(2,5)*{} ; (2,6)*{} **\dir{.},
(2,6)*{} ; (-0.5,8)*{} **\dir{-},
(2,6)*{} ; (1,8)*{} **\dir{-},
(2,6)*{} ; (2,8)*{} **\dir{-},
(2,6)*{} ; (3,8)*{} **\dir{-},
(2,6)*{} ; (4.5,8)*{} **\dir{-},
(2,8)*{} ; (1,10)*{} **\dir{-},
(2,8)*{} ; (2,10)*{} **\dir{-},
(2,8)*{} ; (3,10)*{} **\dir{-},
(-1.2,3.4)*{a^1_1}, (-0.4,3.4)*{\cdots}, (0.6,3.4)*{a^1_{i_1}}, (3.2,3.4)*{c^1_1}, (4,3.4)*{\cdots}, (5,3.4)*{c^1_{j_1}},
(-1,8.5)*{a^k_1}, (0,8.5)*{\cdots}, (1,8.5)*{a^k_{i_k}}, (3,8.5)*{c^k_1}, (4,8.5)*{\cdots}, (5,8.5)*{c^k_{j_k}},
(1,10.5)*{b_1}, (2,10.5)*{\cdots}, (3,10.5)*{b_{l}}
\endxy$$
where some $i_t$ or $j_t$ may be $0$.

An element in $\Lb_{R|A}$ is written $a^1_1 \cdots a^1_{i_1}| \cdots |a^k_1 \cdots a^k_{i_k} [b_1 \cdots b_l] c^k_1 \cdots c_{j_k}^k| \cdots |c_1^1 \cdots c^1_{j_1}$.

A structure of $A_{\infty}$-algebra on $A$ is given by maps $\mu_n : A^{\otimes n} \rightarrow A$ satisfying compatibility relations and a structure of $A$-module over the operad $A_{\infty}$ on $M$ is given by maps $\mu_{n, i} : A^{\otimes i-1}\otimes M \otimes A^{\otimes n-i} \rightarrow M$ for $n \geq 2$ and $1\leq i \leq n$ satisfying some compatibility relations.

In this case, the twisting morphism $\alpha$ is the injection $\Ass^{\ash} \mono \Omega(\Ass^{\ash})$ and the twisting morphism on the level of (co)algebras $\varphi$ is the projection $\Ass^{\ash}(A) \epi A$.

When $d_A = 0$, the differential on the cotangent complex is the sum of three terms that we will make explicit. Otherwise, we have to add a term induced by $d_A$. The first part of the differential is $d_{A\otimes^{A_{\infty}} \Ass^{\ash}(A)}$ given by $d_{\alpha}$ and $d_{A_{\infty}}$.

We use the fact that $\Delta_{p} : \Ass^{\ash} \rightarrow \Ass^{\ash} \circ_{{(1)}} \Ass^{\ash}$ is given by the formula
$$\Delta_{p}(\mu^c_n) = \sum_{\lambda,\, k} (-1)^{\lambda + k(l-\lambda+k)} \mu^c_{l+1-k} \otimes (\underbrace{id\otimes \cdots \otimes id}_{\lambda}\otimes \mu_k^c \otimes \underbrace{id \otimes \cdots \otimes id}_{l-\lambda-k})$$
to give on $\Ass^{\ash}(A)$ the differential
$$d_{\alpha}([b_1 \cdots b_l]) = \sum_{\lambda,\, k} (-1)^{\lambda + k(l-\lambda-k) + (|b_1|+ \cdots + |b_{\lambda}|)(k-1)} [b_1 \cdots b_{\lambda}\, \mu_k (b_{\lambda+1} \cdots b_{\lambda + k}) b_{\lambda + k + 1} \cdots b_l].$$

Contrary to $\Ass$ and $\Ass^{\ash}$, $A_{\infty}$ has a non-zero differential which induces a non-zero differential on $\Lb_{R|A}$ (also denoted $d_{A_{\infty}}$ by abuse of notations). We get\\
$d_{A_{\infty}}(a^1_1 \cdots a^1_{i_1}| \cdots |a^k_1 \cdots a^k_{i_k} [b_1 \cdots b_l] c^k_1 \cdots c_{j_{k}}^{k}| \cdots |c_{1}^{1} \cdots c^1_{j_1}) =$
$$\begin{array}{l}
- \sum \varepsilon_{\lambda,\, k,\, t} a^1_1 \cdots |a^t_1 \cdots \mu_k(a_{\lambda+1}^t \cdots a^t_{\lambda +k}) \cdots a^t_{i_t} | \cdots [\cdots ] c^k_1 \cdots | \cdots | \cdots c^1_{j_1}\\
- \sum \varepsilon_{\lambda,\, k,\, t} a^1_1 \cdots |a^t_1 \cdots a_{\lambda}^t | a_{\lambda+1}^t \cdots a^t_{i_t} | \cdots [\cdots ] \cdots | c^t_1 \cdots c_{k-i_t+\lambda-1}^t | c_{k-i_t+\lambda}^t \cdots c^t_{j_t}  | \cdots c^1_{j_1}\\
- \sum \varepsilon_{\lambda,\, k,\, t} a^1_1 \cdots | \cdots | \cdots a^k_{i_k} [\cdots ] c^k_1 \cdots | c^t_1 \cdots \mu_k(c_{\lambda+1}^t \cdots c^t_{\lambda +k}) \cdots c^t_{j_t}  | \cdots c^1_{j_1},
  \end{array}$$
where $\varepsilon_{\lambda,\, k,\, t} = (-1)^{i_1 + j_1 + \cdots + i_{t-1} + j_{t-1} + \lambda + k(i_t + j_t +1-\lambda+k)}$.

The second part of the differential is the twisted one induced by $\delta_{\varphi}^l$. We get\\
$\delta_{\varphi}^l(a^1_1 \cdots a^1_{i_1}| \cdots |a^k_1 \cdots a^k_{i_k} [b_1 \cdots b_l] c^k_1 \cdots c_{j_{k}}^{k}| \cdots |c_{1}^{1} \cdots c^1_{j_1}) =$
$$\sum_{\lambda,\, k} \epsilon \cdot a^1_1 \cdots |b_1 \cdots b_{\lambda}[b_{\lambda+1} \cdots b_{\lambda+k}] b_{\lambda + k +1} \cdots b_l | \cdots c^1_{j_1},$$
where $\epsilon := (-1)^{i_1 + j_1 + \cdots + i_k + j_k + (|a^1_1| + \cdots + |a^k_{i_k}|)(l-k+1) + (|b_1| + \cdots + |b_{\lambda}|)(k-1) + \lambda + k(l-\lambda + k)}$.

\subsection{The case of $L_{\infty}$-algebras}

The case of $L_{\infty}$-algebras can be made explicit in the same way, with trees in space instead of planar trees. We recover then the definitions given by Hinich and Schechtman in \cite{HinichSchechtman}.

\subsection{The case of $\Po_{\infty}$-algebras}

The general case of homotopy $\Po$-algebras can be treated similarly as follows. Let $\Po$ be a Koszul operad and let $\Po_{\infty} := \Omega (\Poa)$ be its Koszul resolution. Any $\Po_{\infty}$-algebra $A$ admits a resolution $\Po_{\infty} \circ_{\iota} \Poa \circ_{\iota} A \qiso A$, where $\iota : \Poa \rightarrow \Po_{\infty} = \Omega (\Poa)$ is the universal twisting morphism. The cotangent complex has the same form as in the previous cases.

\section{The cotangent complex and the module of K\"ahler differential forms}

In this section, we show that the André-Quillen cohomology of a $\Po$-algebra $A$ is an Ext-functor over the enveloping algebra of $A$ if and only if the cotangent complex of $A$ is a resolution of the module of K\"ahler differential forms. Moreover, we prove that the André-Quillen cohomology theory of an operad is an Ext-functor over its enveloping algebra. We recall that we consider only non-negatively graded $\Po$-algebras in order to have cofibrant resolutions.

\subsection{André-Quillen cohomology as an Ext-functor}\label{AQaaEf}

Let $R$ be a cofibrant resolution of a $\Po$-algebra $A$. Then there is a map $$\Lb_{R/A} = A\otimes^{\Po}_{R}\Omega_{\Po}(R) \rightarrow A\otimes^{\Po}_A \Omega_{\Po}(A) \cong \Omega_{\Po}(A).$$
If the functor $A\otimes^{\Po}_{-}\Omega_{\Po}(-)$ sends cofibrant resolutions to cofibrant resolutions, then the André-Quillen cohomology is the following Ext-functor
$$\textrm{H}_{\Po}^{\bullet}(A,\, M) \cong \textrm{Ext}_{A\otimes^{\Po}\K}^{\bullet}(\Omega_{\Po}(A),\, M).$$

Moreover, we will see in this subsection that the reverse implication is true. Let $X_{\bullet} \qiso \Omega_{\Po}(A)$ be a cofibrant resolution in $\M_{A}^{\Po}$ and consider a quasi-free resolution $R = \Po \circ \C(A)$ of $A$. The cotangent complex $\Lb_{R/A} \cong A\otimes^{\Po}\C(A)$ is a quasi-free $A$-module over $\Po$ since $R$ is quasi-free, so this realization of the cotangent complex is a cofibrant $A$-module over $\Po$. The model category structure on $\M_{A}^{\Po}$ and the commutative diagram
$$\xymatrix@H=16pt{0 \ar[r] \ar@{>->}[d] & X_{\bullet} \ar@{->>}[d]^{\sim}\\
A\otimes^{\Po}\C(A) \ar[r] \ar@{-->}[ur] & \Omega_{\Po}(A)}$$
give a map $A\otimes^{\Po}\C(A) \rightarrow X_{\bullet}$. This last map induces a map
$$\mathrm{H}^{\bullet}_{\Po}(A,\, M) \leftarrow \mathrm{H}^{\bullet}_{\Po}(\mathrm{Hom}_{A\otimes^{\Po}\K \textrm{-}mod}(X_{\bullet},\, M)) \cong \textrm{Ext}_{A\otimes^{\Po}\K}^{\bullet}(\Omega_{\Po}(A),\, M).$$
When this map is an isomorphism, the \emph{André-Quillen cohomology is an Ext-functor over the $\Po$-enveloping algebra}.\\

We prove the following homological lemmas.

\begin{lem}\label{homology1}
Let $\varphi : V \rightarrow W$ be a map of dg vector spaces. If $\varphi^* : V^* \leftarrow W^*$ is an isomorphism then $\varphi : V \rightarrow W$ is an isomorphism, where $V^* := \Hom_{\K}(V,\, \K)$.
\end{lem}

\begin{pf}
Let $x \in V$ non zero and $H$ be a supplementary of $\K x$ in $V = \K x \oplus H$. Since $\varphi^*$ is surjective, there exists $g \in W^*$ such that $x^* = \varphi^* (g) = g \circ \varphi$, where $x^*$ is the map in $V^*$ which is $1$ on $x$ and $0$ on $H$. Thus $1 = x^* (x) = g \circ \varphi (x)$, so $\varphi(x) \neq 0$ and $\varphi$ is injective. Dually we show that $\varphi$ is surjective.
$\cqfd$
\end{pf}

\begin{lem}\label{homology2}
Let $S$ be a dg unitary associative algebra over $\K$ and let $\varphi : M \rightarrow N$ be a map of dg left $S$-modules. If $\varphi^{*} : \Hom_{S\textrm{-}mod}(M,\, M') \xleftarrow{\sim} \Hom_{S\textrm{-}mod}(N,\, M')$ is a quasi-isomorphism for all dg left $S$-module $M'$, then $\varphi : M \xrightarrow{\sim} N$.
\end{lem}

\begin{pf}
We endow $\Hom_{\K}(S,\, \K)$ with a structure of dg left $S$-module by $s\cdot f (x) := f(s^{-1}\cdot x)$ for $s \in S$ and $f \in \Hom_{\K}(S,\, \K)$ and $x \in S$. We have the adjunction
$$\Hom_{S\textrm{-}mod}(M,\, \Hom_{\K}(S,\, \K)) \cong \Hom_{\K}(M \otimes_S S,\, \K) \cong \Hom_{\K}(M,\, \K),$$
which is an isomorphism of dg left $S$-modules (where $\K$ is endowed with a trivial structure). Thus $\varphi^*$ induces a quasi-isomorphism $\Hom_{\K}(M,\, \K) \xleftarrow{\sim} \Hom_{\K}(N,\, \K)$. Since the differential on $\K$ is $0$, we get $\mathrm{H}_{\bullet}(\Hom_{\K}(M,\, \K)) \cong \Hom_{\K}(\mathrm{H}_{\bullet}(M),\, \K)$. We conclude using Lemma \ref{homology1}.
$\cqfd$
\end{pf}

\begin{lem}
Let $\Po$ be a dg operad and let $A$ be a $\Po$-algebra. Let $\C$ be a cooperad and let $\alpha : \C \rightarrow \Po$ be an operadic twisting morphism such that $\Po \circ \C (A)$ is a quasi-free resolution of $A$. There exists a spectral sequence which converges to the cohomology of $A$ with coefficients in $M$, such that
$$E_{2}^{p,\, q} \cong \emph{Ext}_{A\otimes^{\Po}\K}^{p}(\emph{H}_{q}(A\otimes^{\Po}\C(A)),\, M) \Rightarrow \emph{H}_{\Po}^{p+q}(A,\, M).$$
\end{lem}

\begin{pf}
The arguments of Section $5.3.1$ of \cite{Balavoine} are still valid here and give the convergence of the spectral sequence.
$\cqfd$
\end{pf}

\begin{thm}\label{extfunctor}
Let $\Po$ be a dg operad and let $A$ be a $\Po$-algebra. Let $R$ be a cofibrant resolution of $A$. The following properties are equivalent:
\begin{itemize}
\item[$(P_{0})$] the André-Quillen cohomology of $A$ is an Ext-functor over the enveloping algebra $A\otimes^{\Po}\K$, that is $\emph{H}^{\bullet}_{\Po}(A,\, M) \cong \emph{Ext}_{A\otimes^{\Po}\K}^{\bullet}(\Omega_{\Po}(A),\, M)$;
\item[$(P_{1})$] the cotangent complex is quasi-isomorphic to the module of Kähler differential forms, that is $\Lb_{R/A} \qiso \Omega_{\Po}(A)$.
\end{itemize}
\end{thm}

\begin{pf}
A representation of the cotangent complex is given by $A\otimes^{\Po}\C(A)$, where $\C$ is a cooperad and $\alpha : \C \rightarrow \Po$ is a Koszul morphism, e.g. $\C = \B (\Po)$ and $\alpha = \pi$. When $A\otimes^{\Po}\C(A) \qiso \Omega_{\Po}(A)$, as $A\otimes^{\Po}\C(A)$ is a quasi-free $A\otimes^{\Po}\K$-module, the André-Quillen cohomology is by definition an Ext-functor and the property $(P_1)$ implies the property $(P_0)$. Conversely, we assume that H$_{\Po}^{\bullet}(-,\, A)$ is an Ext-functor. We apply Lemma \ref{homology2} to $S = A\otimes^{\Po}\K$, to $M = A\otimes^{\Po}\C(A)$ and to $N = X_{\bullet}$ a cofibrant resolution of $\Omega_{\Po}(A)$. This gives that the property $(P_0)$ implies the property $(P_1)$.
$\cqfd$
\end{pf}

\subsection{André-Quillen cohomology of operads as an Ext-functor}

Rezk defined a cohomology theory for operads following the ideas of Quillen in \cite{Rezk}. Baues, Jibladze and Tonks proposed in \cite{BauesJibladzeTonks} a cohomology theory for monoids in particular monoidal categories, which includes the case of operads. Later Merkulov and Vallette gave in \cite{MerkulovVallette} the cohomology theory ``à la Quillen'' for properads, and so for operads. Merkulov and Vallette define the cotangent complex associated to the resolution of an operad. Let $\Omega(\C) \qiso \Po$ be a cofibrant resolution of the operad $\Po$. We get
$$\mathbb{L}_{\Omega(\C)/\Po} \cong \Po \circ_{(1)}(s^{-1}\overline{\C} \circ \Po) \rightarrow \Omega_{\textrm{op.}}(\Po) \cong \Po \circ_{(1)} (\overline{\Po} \circ \Po)/\! \sim \ \cong \Po \circ_{(1)} \overline{\Po},$$
where $\Omega_{\textrm{op.}}(\Po)$ is the left $\Po$-module of Kähler differential forms (we can see $\Omega_{\textrm{op.}}(\Po)$ as $\Omega_{S}(\Po)$ with $S$ the coloured operad whose algebras are operads). The differential on $\mathbb{L}_{\Omega(\C)/\Po}$ is made explicit as a truncation of the functorial cotangent complex defined in Section \ref{deffunccc}. This enables to define the \emph{André-Quillen cohomology of an operad with coefficients in an infinitesimal $\Po$-bimodule}.

\subsubsection{\bf Infinitesimal bimodule}

An infinitesimal $\Po$-bimodule is an $\So$-module $M$ endowed with two degree $0$ maps $\Po \circ (\Po,\, M) \rightarrow M$ and $M \circ \Po \rightarrow M$ satisfying the commutativity of certain diagrams. We refer to Section $3$ of \cite{MerkulovVallette} for an explicit definition.\\

The notion of operad is a generalization of the notion of associative algebra. Thus, the following lemma can be seen as a generalization of the one in the case of associative algebra.

\begin{lem}\label{acycomplcot}
Let $\Po$ be an augmented dg operad and $\Omega(\C) \qiso \Po$ be a cofibrant resolution. The map $\Po \circ_{(1)}(s^{-1}\overline{\C} \circ \Po) \rightarrow \Omega_{\textrm{op.}}(\Po) \cong \Po \circ_{(1)} \overline{\Po}$ is a quasi-isomorphism.
\end{lem}

\begin{pf}
Since the result does not depend on the cofibrant resolution, we show it in the underlying case $\C = \B(\Po)$. We filter the complex $\Po \circ_{(1)}(s^{-1}\overline{\B(\Po)}\circ \Po)$ by the total number of elements of $\overline{\Po}$ in $\overline{\B(\Po)}\circ \Po$
$$F_p \Po \circ_{(1)}(s^{-1}\overline{\B(\Po)}\circ \Po) := \bigoplus_{w + k \leq p} \Po \circ_{(1)} (s^{-1}\overline{\B_{(w)}(\Po)} \circ (I \oplus \underbrace{\overline{\Po}}_{k \textrm{ times}})).$$
The differential in $\Po \circ_{(1)}(s^{-1}\overline{\B(\Po)}\circ \Po)$ is given by $d_{\Po \circ_{(1)}(s^{-1}\overline{\B(\Po)} \circ \Po)} - \delta^l + \delta^r$. The term $-\delta^l$ decreases $w$ and possibly $k$. The part of $d_{\Po \circ_{(1)}(s^{-1}\overline{\B(\Po)} \circ \Po)}$ induced by $d_{\Po}$ keeps $w+k$ constant and the part induced by $d_2$ of $\B(\Po)$ keeps $w+k$ constant when the application of $\gamma$ is given by $\Po \circ I \cong \Po \cong I \circ \Po$ and decreases $w+k$ by one otherwise. The term $\delta^r$ behaves as the part of the differential induced by $d_2$. Then, the differential respects the filtration. The filtration is bounded below and exhaustive so we can apply the classical theorem of convergence of spectral sequence (cf. Theorem $5.5.1$ of \cite{Weibel}) to obtain that the spectral sequence associated to the filtration converges to the homology of $\Po \circ_{(1)}(s^{-1}\overline{\B(\Po)}\circ \Po)$. The differential $d^0$ on the $E^0_{p,\, \bullet}$ page is given by $d_{\Po} \circ_{(1)} id_{s^{-1}\overline{\B(\Po^{triv})} \circ \Po^{triv}} + id_{\Po} \circ_{(1)} d_{s^{-1}\overline{\B(\Po^{triv})} \circ \Po^{triv}}$, where $\Po^{triv}$ is the underlying dg $\So$-module of $\Po$ endowed with a trivial composition structure, that is $\gamma_{\Po^{triv}} \equiv 0$. By Maschke's theorem, since $\K$ is a field of characteristic $0$, every $\K[\So_{n}]$-module is projective. Then, by the K\"unneth formula, we get
$$\mathrm{H}_{\bullet}(E^{0}_{p,\, \bullet}) = \mathrm{H}_{\bullet}(\Po \circ_{(1)} (s^{-1}\overline{\B(\Po^{triv})} \circ \Po^{triv})) = \mathrm{H}_{\bullet}(\Po) \circ_{(1)} \mathrm{H}_{\bullet}(s^{-1}\overline{\B(\Po^{triv})} \circ \Po^{triv}).$$
Similarly to the proof of $I \qiso B(\Po) \circ \Po$ (see Theorem $2.19$ in \cite{GetzlerJones}), we see that $\overline{\Po} \qiso s^{-1}\overline{B(\Po)}\circ \Po$. Then, for $\Po^{triv}$, we have
$$\mathrm{H}_{\bullet}(E^{0}_{p,\, \bullet}) = \mathrm{H}_{\bullet}(\Po) \circ_{(1)} \mathrm{H}_{\bullet}(s^{-1}\overline{\B(\Po^{triv})} \circ \Po^{triv}) = \mathrm{H}_{\bullet}(\Po) \circ_{(1)} \mathrm{H}_{\bullet}(\overline{\Po}) = \mathrm{H}_{\bullet}(\Po \circ_{(1)} \overline{\Po}).$$
Finally, the spectral sequence collapses at rank $1$ and the Lemma is true.
$\cqfd$
\end{pf}

As a corollary of the previous Lemma, we get

\begin{thm}
The André-Quillen cohomology of operads with coefficients in an infinitesimal $\Po$-bimodule is the Ext-functor
$$\emph{H}^{\bullet}(\Po,\, \mathcal{M}) \cong \emph{Ext}_{\Po \circ_{(1)}(I \circ \Po)}^{\bullet}(\Omega_{\textrm{op.}}(\Po),\, \mathcal{M}).$$
\end{thm}

\begin{pf}
We combine Theorem \ref{extfunctor} and Lemma \ref{acycomplcot}.
$\cqfd$
\end{pf}

\section{The functorial cotangent complex}

In this section, we introduce a \emph{functorial cotangent complex} and a \emph{functorial module of K\"ahler differential forms}, depending only on the operad. 
We prove that when the André-Quillen cohomology is an Ext-functor then the map between these two complexes is a quasi-isomorphism. 
We define the module of obstructions and we show that it is acyclic when the André-Quillen cohomology is an Ext-functor, giving a way, by contraposition, to show that the André-Quillen cohomology isn't an Ext-functor. 

Under some \emph{PBW conditions}, we prove that the modules of obstructions is acyclic if and only if the André-Quillen cohomology is an Ext-functor.

\subsection{Definition of the functorial cotangent complex}\label{deffunccc}

As we explain in Section \ref{operadicresolution}, the resolutions of algebras we use in this paper come from operadic resolutions. They all have the form $\Po \circ_{\alpha} \C \qiso I$, where $\alpha : \C \rightarrow \Po$ is an operadic twisting morphism. We call such a twisting morphism a \emph{Koszul morphism}. We define (a representation of) the \emph{functorial cotangent complex} based on such type of resolutions as follows.

We consider the dg infinitesimal $\Po$-bimodule $L_{\Po} := \Po (I,\, \C \circ \Po) = \Po \circ_{(1)} (\C \circ \Po)$ endowed with the differential $d_{L_{\Po}} := d_{\Po (I,\, \C \circ \Po)} - \delta_{L_{\Po}}^{l} + \delta_{L_{\Po}}^{r}$, where $\delta_{L_{\Po}}^{l}$ is defined by the composite
$$\hspace{-1cm} \Po \circ_{(1)} (\C \circ \Po) \xrightarrow{id_{\Po} \circ_{(1)} (\Delta_{p} \circ id_{\Po})} \Po \circ_{(1)} ((\C \circ_{(1)} \C) \circ \Po) \xrightarrow{id_{\Po} \circ_{(1)} (\alpha \circ id_{\C} \circ id_{\Po})}$$
$$\hspace{1cm} \Po \circ_{(1)} ((\Po \circ_{(1)} \C) \circ \Po) \mono (\Po \circ \Po \circ \Po) \circ_{(1)} (\C \circ \Po) \xrightarrow{\gamma \circ \gamma \circ_{(1)} id_{\C \circ \Po}} \Po \circ_{(1)} (\C \circ \Po)$$
and $\delta_{L_{\Po}}^{r}$ is defined by the composite
$$\hspace{-1cm} \Po \circ_{(1)} (\C \circ \Po) \xrightarrow{id_{\Po} \circ_{(1)} (\Delta_{p} \circ id_{\Po})} \Po \circ_{(1)} ((\C \circ_{(1)} \C) \circ \Po) \xrightarrow{id_{\Po} \circ_{(1)} (id_{\C} \circ \alpha \circ id_{\Po})}$$
$$\hspace{1cm} \Po \circ_{(1)} ((\C \circ_{(1)} \Po) \circ \Po) \mono \Po \circ_{(1)} (\C \circ \Po \circ \Po) \xrightarrow{id_{\Po} \circ_{(1)} (id_{\C} \circ \gamma)} \Po \circ_{(1)} (\C \circ \Po).$$
The right action is given by $\Po \circ_{(1)} (\C \circ \Po) \circ \Po \mono (\Po \circ \Po) \circ_{(1)} (\C \circ \Po \circ \Po) \xrightarrow{\gamma \circ_{(1)} id_{\C} \circ \gamma} \Po \circ_{(1)} (\C \circ \Po).$

\begin{prop}
Let $A$ be a $\Po$-algebra. With the above notations, there is an isomorphism of chain complexes
$$L_{\Po} \circ_{\Po} A \cong A\otimes^{\Po}\C (A).$$
\end{prop}

\begin{pf}
We write $L_{\Po} \circ A \cong \Po (A,\, \C \circ \Po (A))$. We use the description of the relative composition product $\circ_{\Po}$ and of the description $A\otimes^{\Po}N$ to get $L_{\Po} \circ_{\Po} A \cong A\otimes^{\Po}\C (A)$. The equality of the differentials comes from the same descriptions.
$\cqfd$
\end{pf}

\begin{cor}
Let $V$ be a dg trivial $\Po$-algebra, that is $\gamma_V \equiv 0$. There is an isomorphism of chain complexes
$$(L_{\Po} \circ_{\Po} I) \circ V \cong V \otimes^{\Po}\C(V).$$
\end{cor}

\begin{pf}
When the $\Po$-algebra $V$ is trivial, we get the isomorphism of underlying dg modules $(L_{\Po} \circ_{\Po} I) \circ V \cong L_{\Po} \circ_{\Po} V$, where $I$ can be seen as a left $\Po$-module with a trivial structure. The equality of the differentials follows from their definitions.
$\cqfd$
\end{pf}

We denote $\overline{L}_{\Po} := L_{\Po} \circ_{\Po} I$.

\subsection{Definition of the functorial module of Kähler differential forms}

Let $\Po$ be a dg operad. We define the \emph{functorial module of Kähler differential forms} as the following coequalizer diagram in the category of infinitesimal $\Po$-bimodules (see $10.3$ of \cite{Fresse3} for an equivalent definition)
$$\xymatrix@C=1.5cm{\Po \circ_{(1)} (\Po \circ \Po) \ar@<0.5ex>[r]^{\ \ id_{\Po} \circ_{(1)} \gamma} \ar@<-0.5ex>[r]_{\ \ c_{2}} & \Po \circ_{(1)} \Po \ar@{->>}[r] & \Omega_{\Po},}$$
where $c_{2}$ is given by the composite
$$\Po \circ_{(1)} (\Po \circ \Po) \xrightarrow{id_{\Po} \circ_{(1)} (id_{\Po} \circ' id_{\Po})} \Po \circ_{(1)} (\Po \circ (\Po,\, \Po)) \mono (\Po \circ \Po \circ \Po) \circ_{(1)} \Po \xrightarrow{\gamma \circ \gamma \circ_{(1)} id_{\Po}} \Po \circ_{(1)} \Po.$$
The right $\Po$-module action on $\Omega_{\Po}$ is induced by the right $\Po$-module action on $\Po \circ_{(1)} \Po$ given by
$$(\Po \circ_{(1)} \Po) \circ \Po \mono (\Po \circ \Po) \circ_{(1)} (\Po \circ \Po) \xrightarrow{\gamma \circ_{(1)} \gamma} \Po \circ_{(1)} \Po.$$

\begin{prop}
Let $A$ be a $\Po$-algebra. There is an isomorphism of chain complexes
$$\Omega_{\Po} \circ_{\Po} A \cong \Omega_{\Po} (A).$$
\end{prop}

\begin{pf}
We write $A \otimes^{\Po} A \cong (\Po \circ_{(1)} \Po) \circ_{\Po} A$ and $A \otimes^{\Po} \Po(A) \cong (\Po \circ_{(1)} \Po \circ \Po) \circ_{\Po} A$. Thanks to the description of $\Omega_{\Po}(A)$ given at the end of Lemma \ref{bij}, we get the result.
$\cqfd$
\end{pf}

\begin{cor}
Let $V$ be a dg trivial $\Po$-algebra. There is an isomorphism of chain complexes
$$(\Omega_{\Po} \circ_{\Po} I) \circ V \cong \Omega_{\Po}(V).$$
\end{cor}

\begin{pf}
When the $\Po$-algebra $V$ is trivial, we get $(\Omega_{\Po} \circ_{\Po} I) \circ V \cong \Omega_{\Po} \circ_{\Po} \circ V$.
$\cqfd$
\end{pf}

We denote $\overline{\Omega}_{\Po} := \Omega_{\Po} \circ_{\Po} I$.

\subsection{Homotopy category}

Let $\Po$ be an augmented dg operad and $\C$ be a coaugmented dg cooperad such that $\Po \circ_{\alpha} \C \qiso I$. We define the following surjective map of infinitesimal $\Po$-bimodules
$$L_{\Po} = \Po (I,\, \C \circ \Po) \epi \Po (I,\, I \circ \Po)/\! \sim \ \cong \Po \circ_{(1)} \Po /\! \sim \ \cong \Omega_{\Po}.$$
This map induces a map
$$A\otimes^{\Po} \C(A) \cong L_{\Po} \circ_{\Po} A \epi \Omega_{\Po} \circ_{\Po} A \cong \Omega_{\Po}(A)$$
which coincides with the map given in Section \ref{AQaaEf}.

The differential on $\mathbb{L}_{\Omega(\C)/\Po}$ and the augmentation $\mathbb{L}_{\Omega(\C)/\Po} \rightarrow \Po \circ_{(1)} (I \circ \Po)$ induce a differential on the cone $s\mathbb{L}_{\Omega(\C)/\Po} \oplus \Po \circ_{(1)} (I \circ \Po)$. With this differential, we have $L_{\Po} \cong s\mathbb{L}_{\Omega(\C)/\Po} \oplus \Po \circ_{(1)} (I \circ \Po)$. Then, $L_{\Po}$ is well-defined in the homotopy category of infinitesimal $\Po$-bimodules. The same is true for $\overline{L}_{\Po} = L_{\Po}\circ_{\Po}I$ and we call its image in the homotopy category of infinitesimal left $\Po$-modules the \emph{functorial cotangent complex}, that we denote by $\LbP$. We denote by $\mathbf{\Omega}_{\Po}$ the image of $\overline{\Omega}_{\Po} := \Omega_{\Po} \circ_{\Po} I$ in the homotopy category.

\subsection{Filtration on the cotangent complex}

Let $\alpha : \C \rightarrow \Po$ be a Koszul morphism between a connected weight graded dg cooperad and a weight graded dg operad. Let $A$ be a dg $\Po$-algebra. We filter $L_{\Po} \circ A \cong \Po (A,\, \C \circ \Po(A))$ by the weight in the first $\Po$  and the weight in $\C$:
$$F_{p} (L_{\Po} \circ A) := \bigoplus_{m+n\leq p} \Po^{(n)}\circ_{(1)} (\C^{(m)} \circ \Po) \circ A.$$

With the projection $L_{\Po} \circ A \epi L_{\Po} \circ_{\Po} A \cong A\otimes^{\Po} \C(A)$, it induces a filtration on $A\otimes^{\Po} \C(A)$ by $F_{p} (A\otimes^{\Po} \C(A)) := \textrm{im}(F_p (L_{\Po} \circ A) \epi A\otimes^{\Po} \C(A))$.

The differential on $A\otimes^{\Po} \C(A)$ is given by $id_{L_{\Po}} \circ_{\Po}' d_A + d_{\Po (I,\, \C \circ \Po)} \circ_{\Po} id_A - \delta_{L_{\Po}}^{l} \circ_{\Po} id_A + \delta_{L_{\Po}}^{r} \circ_{\Po} id_A$. 
The term depending on $d_\C$ keeps the sum $n+m$ constant (since it preserves the weight), the part $id_{L_{\Po}} \circ_{\Po}' d_A
- \delta_{L_{\Po}}^{l} \circ_{\Po} id_A$ and the term depending on $d_\Po$ may decrease the sum $n+m$ and the part $\delta_{L_{\Po}}^{r} \circ_{\Po} id_A$ decreases the sum $n+m$ (at least by 1 since $\C$ is connected weight graded
). It follows that the differential on the cotangent complex respects this filtration.

\begin{lem}
For any $\Po$-algebra $A$, the spectral sequence associated to the filtration $F_{p}$ converges to the homology of the cotangent complex
$$E_{p,\, q}^1 = \emph{H}_{p+q}(F_{p} (A\otimes^{\Po} \C(A))/F_{p-1} (A\otimes^{\Po} \C(A))) \Rightarrow \emph{H}_{p+q} (A\otimes^{\Po} \C(A)).$$
\end{lem}

\begin{pf}
This filtration is exhaustive and bounded below so we can apply the classical theorem of convergence of spectral sequences (cf. Theorem $5.5.1$ of \cite{Weibel}) to obtain the result.
$\cqfd$
\end{pf}

We denote by $d^{0}$ the differential on $E_{p,\, \bullet}^0$, which depends on $id_{L_{\Po}} \circ_{\Po}' d_A$, on $d_{\Po (I,\, \C \circ \Po)} \circ_{\Po} id_A$ and on $\delta_{L_{\Po}}^{l} \circ_{\Po} id_A$. 
We denote by $d^{r}$ the differential on $E_{p,\, \bullet}^r$
.

\subsection{Filtration of the module of K\"ahler differential forms}

By its definition as a coequalizer, $\Omega_{\Po}(A)$ 
is a quotient of $A\otimes^{\Po} A$. 
We filter $A \otimes^{\Po}A$ by the weight in $\Po$ and this induces a filtration $F_p'$ on $\Omega_{\Po}(A)$. The differential on $\Omega_{\Po}(A)$ respects the filtration
. We denote by ${d'}^0$ the differential on ${E'}_{p,\, \bullet}^0$, which is induced by the differential on $\Omega_{\Po}(A)$. 

Then, for any $\Po$-algebra $A$, the filtration $F'_p$ is exhaustive and bounded below, so the spectral sequence associated to $F'_{p}$ converges to the homology of the module of K\"ahler differential forms
$${E'}_{p,\, q}^1 = \emph{H}_{p+q}(F'_{p} (\Omega_{\Po}(A))/F'_{p-1} (\Omega_{\Po}(A))) \Rightarrow \emph{H}_{p+q} (\Omega_{\Po}(A)).$$

The map $A\otimes^{\Po} \C(A) \to \Omega_{\Po}(A)$ considered above is compatible with the filtrations $F_\bullet$ and $F'_\bullet$.

\subsection{The cotangent complex and the module of K\"ahler differential forms}

As in the previous sections, we endow the enveloping algebra $A\otimes^{\Po} \K$ with a filtration given by the weight in $\Po$. Following the notations given by Pirashvili in the review of \cite{Frabetti}, we say that $\Po$ is an operad satisfying the PBW property if for any $\Po$-algebra $A$, there is an isomorphism $gr(A\otimes^{\Po} \K) \cong A^{\mathrm{tr}} \otimes^{\Po}\K$, where $A^{\mathrm{tr}}$ is the underlying space of $A$ endowed with the trivial $\Po$-algebra structure $\gamma_{A^{\mathrm{tr}}} \equiv 0$. 
The study of the differential on the cotangent complex $A\otimes^{\Po} \C (A)$ shows that $\Po$ is an operad satisfying the PBW property if and only if for any $\Po$-algebra $A$, we have the isomorphism $gr(A \otimes^{\Po}\C (A)) \cong A^{\mathrm{tr}} \otimes^{\Po} \C (A^{\mathrm{tr}})$. 

We similarly say that the $\Po$-K\"ahler differentials satisfy the PBW property when for any $\Po$-algebra $A$ we have the isomorphism $gr(\Omega_{\Po}(A)) \cong \Omega_{\Po}(A^{\mathrm{tr}})$. 
These notions are different from the notion of PBW-operad defined in \cite{Hoffbeck}.

To shorten the following theorem, we say that $\Po$ is a \emph{PBW-operad} when
\begin{itemize}
\item
it is a connected weight graded operad (in this case, we can consider the cooperad $\C = \B\Po$ which is connected weight graded),
\item
it satisfies the PBW property, and
\item
$\Po$-K\"ahler differentials satisfy the PBW property.
\end{itemize}

We complete Theorem \ref{extfunctor} as follows.

\begin{thm}\label{acyclicity trivial}
The assertion
\begin{itemize}
\item[$(P_{0})$] the André-Quillen cohomology is an Ext-functor over the enveloping algebra $A\otimes^{\Po}\K$ for any $\Po$-algebra $A$;
\end{itemize}
implies the following properties equivalent assertions:
\begin{itemize}
\item[$(P_{1}')$] the cotangent complex is quasi-isomorphic to the module of Kähler differential forms for any dg vector space $V$, seen as an algebra with trivial structure, that is $\Lb_{R/V} \qiso \Omega_{\Po}(V)$;
\item[$(P_{2})$] the functorial cotangent complex $\LbP$ is quasi-isomorphic to the functorial module of K\"ahler differential forms $\mathbf{\Omega}_{\Po}$, that is $\LbP \qiso \mathbf{\Omega}_{\Po}$.
\end{itemize}
Under the assumption that $\Po$ is a PBW-operad, we have the equivalences $(P_0) \Leftrightarrow (P_1') \Leftrightarrow (P_2)$.
\end{thm}

\begin{pf}
The implication $(P_0) \Rightarrow (P_1')$ is clear by Theorem \ref{extfunctor}. 
The equivalence $(P_1') \Leftrightarrow (P_2)$ follows from the equalities $(V\otimes^{\Po} \C(V), d_{\varphi}) = (\overline{L}_{\Po} \circ V,\, d_{\overline{L}_{\Po} \circ V}) = (\overline{L}_{\Po} \circ V,\, d_{\overline{L}_{\Po}} \circ id_V)$ and $(\Omega_{\Po}(V),\, d_{\Omega_{\Po}(V)}) = (\overline{\Omega}_{\Po} \circ V,\, d_{\overline{\Omega}_{\Po}} \circ id_{V})$.

We assume now that $\Po$ be a connected weight graded operad satisfying the PBW property such that the $\Po$-Kähler differentials also satisfy the PBW property.
We prove the implication $(P_0) \Leftarrow (P_1')$. 
Suppose that the cotangent complex to be quasi-isomorphic to the module of K\"ahler differential forms for any dg vector space. 
Let $A$ be a $\Po$-algebra and denote by $V$ the underlying dg vector space of $A$ considered as a trivial algebra. We use the filtrations and the spectral sequences of the two previous sections. 
We have
$$
E_{p,\, q}^0 = \left(gr(A \otimes^{\Po}\C (A)),\, d^0_A\right) \cong \left(A^{\mathrm{tr}} \otimes^{\Po} \C (A^{\mathrm{tr}}),\, d^0_{A^{\mathrm{tr}}}\right)
$$
by the fact that $\Po$ satisfies the PBW property. 
Similarly,
$$
{E'}_{p,\, q}^0 = \left(gr(\Omega_{\Po}(A)),\, {d'}^0_A\right) \cong \left(\Omega_{\Po}(A^{\mathrm{tr}}),\, {d'}^0_{A^{\mathrm{tr}}}\right)
$$
because the $\Po$-Kähler differentials also satisfy the PBW property. 
By $(P'_1)$, we have a quasi-isomorphism $(E^0,\, d^0) \qiso ({E'}^0,\,  {d'}^0)$. 
This map is induced by the counit $\C \to I$ and it follows that $E^0$ is a subquotient of $A \otimes^{\Po}I (A)$ ($\C$ is coaugmented). 
Then the higher differentials $d^r_A$ and ${d'}^r_A$ coincide through the (quasi-)isomorphism because they only depend on $d_A$ and $d_\Po$. 
So we obtain $(P_1)$, and therefore $(P_0)$ by Theorem \ref{extfunctor}.
$\cqfd$
\end{pf}

\begin{rem}
When $\Po$ is an operad concentrated in homological degree $0$, $\overline{\Omega}_{\Po}$ is an $\So$-module concentrated in degree $0$. In this case, we say that $\LbP$ is acyclic when its homology is concentrated in degree $0$ and equal to $\mathbf{\Omega}_{\Po}$.
\end{rem}

\subsubsection{\bf First applications}\label{acyclicity cotangent complex}
The operads encoding associative algebras and Lie algebras are PBW-operads. 
Indeed they are connected weight graded since they admit a homogeneous quadratic presentation. Then the operad $\Ass$ satisfy the PBW property by means of the computation made in the Examples after Proposition \ref{adjunction1}. We have
\begin{align*}
gr(A \otimes^{\Ass} \K) & \cong gr^0 \oplus gr^1 \oplus gr^2 = \K \oplus \left(A\otimes \K \oplus \K \otimes A\right) \oplus (A \otimes A \otimes A)\\
& \cong (\K \oplus A^{\mathrm{tr}}) \otimes A^{\mathrm{tr}} \otimes (\K \oplus A^{\mathrm{tr}}) \cong A^{\textrm{tr}} \otimes^\Po \K.
\end{align*}
For the operad $\Lie$, the corresponding statement is the meaning of the Poincaré-Birkhoff-Witt theorem. 
It remains to show that the $\Ass$ (resp. $\Lie$)-Kähler differentials satisfy the PBW property. 
In the case of the operad $\Ass$, a computation similar as the computation for the enveloping algebra gives
$$
gr(\Omega_\Po (A)) \cong gr^0 \oplus gr^1 \cong dA \oplus dA\otimes A \cong dA^{\mathrm{tr}} \oplus dA^{\mathrm{tr}}\otimes A^{\mathrm{tr}} \cong \Omega_\Po(A^{\mathrm{tr}}),
$$
where we have noted $dA$ the ``linear'' copie of $A$ in $\Omega_\Po (A) \cong A \otimes^\Po dA /\sim$ and where we have made used of the relation $d(a \cdot b) = da \otimes b \pm a \otimes db$ and $a \otimes db \otimes c = \pm (d(a\cdot b) \otimes c - a \otimes d(b\cdot c))$. 
The case of the operad $\Lie$ can be proven analogously (we find $gr(\Omega_{\Lie}(A)) \cong S(A) \cong \Omega_{\Lie}(A^{\mathrm{tr}})$).

We now prove the acyclicity of the functorial cotangent complex in these cases. This gives a conceptual proof of the fact that for these operads the André-Quillen cohomology is an Ext-functor over the enveloping algebra $A\otimes^{\Po}\K$.
\begin{itemize}
 \item We have $L_{\Ass} = \vcenter{
\xymatrix@M=0pt@R=6pt@C=4pt{\Ass^{\ash}\\
\ar@{-}[dd]\\ &\\ &}} \oplus \vcenter{\xy (2,0)*{} ; (2,2)*{} **\dir{-},
(2,2)*{} ; (0,4)*{} **\dir{-},
(2,2)*{} ; (4,4)*{} **\dir{-},
(6,7)*{\Ass^{\ash}}
\endxy} \oplus \vcenter{\xy (4,0)*{} ; (4,2)*{} **\dir{-},
(4,2)*{} ; (2,4)*{} **\dir{-},
(4,2)*{} ; (6,4)*{} **\dir{-},
(0,7)*{\Ass^{\ash}}
\endxy} \oplus \vcenter{
\xymatrix@M=0pt@R=6pt@C=4pt{ & \Ass^{\ash} & \\
\ar@{-}[dr] & \ar@{-}[dd] & \ar@{-}[dl]\\
& &  \\ & &}}$. Then $L_{\Ass} (n)$ is generated by the elements $u_{n} := \vcenter{\xymatrix@M=0pt@R=6pt@C=4pt{1 && \cdots && n\\
\ar@{-}[drr] & \ar@{-}[dr] & \ar@{-}[dd] & \ar@{-}[dl] & \ar@{-}[dll]\\
\ar@{.}[rrrr] &&&&\\
&&&&&}}$, $r_{n} := \vcenter{\xymatrix@M=0pt@R=6pt@C=4pt{&&& 1 && \cdots && n\\
&&& \ar@{-}[drr] & \ar@{-}[dr] & \ar@{-}[d] & \ar@{-}[dl] & \ar@{-}[dll]\\
\ar@{.}[rrrrrrr] & \ar@{-}[drr] &&&& \ar@{-}[dll] &&\\
&&& \ar@{-}[d] &&&&&\\
&&&&&&&&}}$, $l_{n} := \vcenter{\xymatrix@M=0pt@R=6pt@C=4pt{1 && \cdots && n &&&\\
\ar@{-}[drr] & \ar@{-}[dr] & \ar@{-}[d] & \ar@{-}[dl] & \ar@{-}[dll] &&&\\
\ar@{.}[rrrrrrr] && \ar@{-}[drr] &&&& \ar@{-}[dll] &\\
&&&& \ar@{-}[d] &&&\\
&&&&&&&}}$ and $v_{n} := \vcenter{\xymatrix@M=0pt@R=6pt@C=4pt{&&& 1 && \cdots && n &&&\\
&&& \ar@{-}[drr] & \ar@{-}[dr] & \ar@{-}[d] & \ar@{-}[dl] & \ar@{-}[dll] &&&\\
\ar@{.}[rrrrrrrrrr] & \ar@{-}[drrrr] &&&& \ar@{-}[dd] &&&& \ar@{-}[dllll] &\\
&&&&&&&&&&\\
&&&&&&&&&&}}$. Since $d(u_{n}) = - l_{n-1} - (-1)^{n-1}r_{n-1}$, $d(r_{n}) = - v_{n-1} = (-1)^{n-1}d(l_{n})$ and $d(v_{n}) = 0$, we define a homotopy $h$ for $d$ by $h(u_{n}) := 0$, $h(l_{n}) = h(r_{n})(-1)^{n} := -\frac{1}{2} u_{n+1}$ and $h(v_{n}) := -\frac{1}{2} ((-1)^{n}l_{n+1} + r_{n+1})$.
\item We have
$$L_{\Lie} = \vcenter{
\xy (0,0)*{} ; (0,4)*{} **\dir{-},
(0,7)*{\Lie^{\ash}}, (-0.6,4)*{\textrm{\tiny $1$ \normalsize}}
\endxy} \oplus \vcenter{\xy (2,0)*{} ; (2,2)*{} **\dir{-},
(2,2)*{} ; (0,4)*{} **\dir{-},
(2,2)*{} ; (4,4)*{} **\dir{-},
(5,7)*{\Lie^{\ash}}, (-0.6,4)*{\textrm{\tiny $1$ \normalsize}}, (5.4,4)*{\textrm{\tiny $2$ \normalsize}}
\endxy} \oplus \vcenter{\xy (3,0)*{} ; (3,1)*{} **\dir{-},
(3,1)*{} ; (0,4)*{} **\dir{-},
(3,1)*{} ; (6,4)*{} **\dir{-},
(4.5,2.5)*{} ; (3,4)*{} **\dir{-},
(7,7)*{\Lie^{\ash}}, (-0.6,4)*{\textrm{\tiny $1$ \normalsize}}, (2.6,4)*{\textrm{\tiny $2$ \normalsize}}, (7.4,4)*{\textrm{\tiny $3$ \normalsize}}
\endxy} \oplus \cdots \oplus \vcenter{\xy (3,0)*{} ; (3,1)*{} **\dir{-},
(3,1)*{} ; (-1,6)*{} **\dir{-},
(3,1)*{} ; (8,6)*{} **\dir{-},
(6,4)*{} ; (5,6)*{} **\dir{-},
(10,9.5)*{\Lie^{\ash}}, (-1,7)*{\textrm{\tiny $1$ \normalsize}}, (4,4)*{\textrm{\tiny $\cdots$ \normalsize}}, (5,7)*{\textrm{\tiny $n$-$1$ \normalsize}}, (9.4,7)*{\textrm{\tiny $n$ \normalsize}}
\endxy} \oplus \cdots.$$
Then we can define the same homotopy as in \cite{CartanEilenberg}, Theorem $7.1$, Chap. XIII.
\item Following Frabetti in \cite{Frabetti}, we show the acyclicity of $\mathbb{L}_{\Dias}$.
\end{itemize}

\begin{rem}
We recall the following results.
\begin{itemize}
\item Loday and Pirashvili showed in \cite{LodayPirashvili} that the cohomology of Leibniz algebras can be written as an Ext-functor.
\item Dzhumadil'daev showed in \cite{Dzhumadil} that the cohomology of pre-Lie algebras can be written as an Ext-functor.
\end{itemize}
\end{rem}

\subsubsection{\bf The module of obstructions}\label{moduleobstruction}

Let $\Po$ be an augmented dg operad and let $\C$ be a coaugmented dg cooperad and let $\alpha : \C \rightarrow \Po$ be a twisting morphism. The map $\overline{L}_{\Po} \epi \overline{\Omega}_{\Po}$ is surjective and we defined
$$O_{\Po} := \textrm{ker} (\overline{L}_{\Po} \epi \overline{\Omega}_{\Po})$$
to get the following short exact sequence of dg $\So$-modules
$$O_{\Po} \mono \overline{L}_{\Po} \epi \overline{\Omega}_{\Po}.$$
Since $\overline{L}_{\Po}$ and $\overline{\Omega}_{\Po}$ are well-defined in the homotopy category of infinitesimal left $\Po$-modules, the same is true for $O_{\Po}$. Thus we define the \emph{module of obstructions} $\Kc$ by its image in the homotopy category of infinitesimal $\Po$-modules. We get the following short exact sequence
$$\Kc \mono \LbP \epi {\bf \Omega_{\Po}}.$$

We compute
$$O_{\Po} = Rel \oplus (\Po (I,\, \C \circ \Po)) \circ_{\Po} I,$$
where $Rel$ is the image of the set of relations
$$\{ \def\objectstyle{\scriptstyle}
\def\labelstyle{\scriptstyle}
\vcenter{\xymatrix@M=0pt@R=4pt@C=4pt{
& \ar@{-}[d] &&&\\
\ar@{-}[drr] & \ar@{-}[dr] & \ar@{-}[d] & \ar@{-}[dl] & \ar@{-}[dll]\\
&& \ar@{-}[d] &&\\
&&&&}} + \def\objectstyle{\scriptstyle}
\def\labelstyle{\scriptstyle}
\vcenter{\xymatrix@M=0pt@R=4pt@C=4pt{
&& \ar@{-}[d] &&\\
\ar@{-}[drr] & \ar@{-}[dr] & \ar@{-}[d] & \ar@{-}[dl] & \ar@{-}[dll]\\
&& \ar@{-}[d] &&\\
&&&&}} + \def\objectstyle{\scriptstyle}
\def\labelstyle{\scriptstyle}
\vcenter{\xymatrix@M=0pt@R=4pt@C=4pt{
&&& \ar@{-}[d] &\\
\ar@{-}[drr] & \ar@{-}[dr] & \ar@{-}[d] & \ar@{-}[dl] & \ar@{-}[dll]\\
&& \ar@{-}[d] &&\\
&&&&}}, \textrm{ where } \def\objectstyle{\scriptstyle}
\def\labelstyle{\scriptstyle}
\vcenter{\xymatrix@M=0pt@R=4pt@C=4pt{
\ar@{-}[drr] & \ar@{-}[dr] & \ar@{-}[d] & \ar@{-}[dl] & \ar@{-}[dll]\\
&& \ar@{-}[d] &&\\
&&&&}} = \gamma \left(\def\objectstyle{\scriptstyle}
\def\labelstyle{\scriptstyle}
\vcenter{\xymatrix@M=0pt@R=4pt@C=4pt{
& \ar@{-}[rd] & \ar@{-}[d] & \ar@{-}[ld] &\\
\ar@{.}[rrrr] \ar@{-}[drr] && \ar@{-}[d] && \ar@{-}[dll] \\
&& \ar@{-}[d] &&\\
&&&&}} \right) \textrm{ and } \def\objectstyle{\scriptstyle}
\def\labelstyle{\scriptstyle}
\vcenter{\xymatrix@M=0pt@R=4pt@C=4pt{
\ar@{-}[drr] && \ar@{-}[d] && \ar@{-}[dll] \\
&& \ar@{-}[d] &&\\
&&&&}},\, \def\objectstyle{\scriptstyle}
\def\labelstyle{\scriptstyle}
\vcenter{\xymatrix@M=0pt@R=4pt@C=4pt{
& \ar@{-}[rd] & \ar@{-}[d] & \ar@{-}[ld] &\\
&& \ar@{-}[d] &&\\
&&&&}}\, \in \Po \}$$
in $(\Po \circ_{(1)} I) \circ_{\Po} I$.

We deduce the following theorem.

\begin{thm}\label{obstruction}
Let $\Po$ be an augmented dg operad. The property
\begin{itemize}
\item[$(P_{0})$] The André-Quillen cohomology is an Ext-functor over the enveloping algebra $A\otimes^{\Po}\K$ for any $\Po$-algebra $A$;
\end{itemize}
implies the following property
\begin{itemize}
\item[$(P_{3})$] the homology of the module of obstructions $\Kc$ is acyclic.
\end{itemize}
When $\Po$ be a PBW-operad, we obtain the equivalence $(P_{0}) \Leftrightarrow (P_{3})$.
\end{thm}

\begin{pf}
The short exact sequence $O_{\Po} \mono \overline{L}_{\Po} \epi \overline{\Omega}_{\Po}$ induces a long exact sequence in homology which gives the equivalence $(P_3) \Leftrightarrow (P_2)$. 
Then the theorem follows from Theorem \ref{acyclicity trivial}.
$\cqfd$
\end{pf}

\subsection{Another approach}

In the parallel work \cite{Fresse3}, Fresse studied the homotopy properties of modules over operads. His method applied to the present question provides the following sufficient condition for the André-Quillen cohomology to be an Ext-functor. In this section, we show the relationship between the two approaches.\\

Let $\Po[1]$ be the $\So$-module defined in \cite{Fresse3} given by $\Po[1](n) := \Po(1+n)$. The $\So_n$-action is given by the action of $\So_n$ on $\{2,\ldots ,\, n+1 \} \subset \{1,\ldots ,\, n+1 \}$. Similarly to this definition, we define the $\So$-module $\Po[1]_j$ by $\Po[1]_j(n) := \Po(n+1)$ where the $\So_n$-action is given by the action $\So_n$ on $\{1,\ldots ,\, \hat{j},\, \ldots ,\, n+1 \} \subset \{1,\ldots ,\, n+1 \}$. Thus $\Po[1] = \Po[1]_1$. We have $(\Po \circ_{(1)} I) (n) \cong \underbrace{\Po(n) \oplus \cdots \oplus \Po(n)}_{n \textrm{ times}}$. As a right $\Po$-module, we have
$$(\Po \circ_{(1)} I)(n) \cong \underbrace{\Po[1]_1(n-1)\oplus \cdots \oplus \Po[1]_n(n-1)}_{n \textrm{ times}}.$$
When $\Po[1]$ is a quasi-free right $\Po$-module, that is $\Po[1] \cong (M \circ \Po,\, d)$, we get that $\Po[1]_j$ is a quasi-free right $\Po$-module $(M \circ \Po,\, d)$ thanks to the isomorphism
$$\Po[1] \rightarrow \Po[1]_j,\, \mu \mapsto \mu \cdot (1 \cdots j).$$
We define
$$M'(n) := \oplus_{k \geq 1} \underbrace{M(n) \oplus \cdots \oplus M(n)}_{k \textrm{ times}}.$$ Then $\Po \circ_{(1)} I$ is a retract of $(M'\circ \Po,\, d')$, which is quasi-free. When $\Po[1]$ is only a retract of a quasi-free right $\Po$-module, we get by the same argument that $\Po \circ_{(1)} I$ is a retract of a quasi-free right $\Po$-module. Then, $\Po[1]$ cofibrant as a right $\Po$-module implies that $\Po \circ_{(1)} I$ cofibrant as a right $\Po$-module. Finaly, when $\Po[1]$ is a cofibrant right $\Po$-module, $L_{\Po} \cong (\Po \circ_{(1)} I) \otimes (\C \circ \Po)$ is also cofibrant. Thus, when we assume moreover that $\Omega_{\Po}$ is a cofibrant right $\Po$-module, the quasi-isomorphism $L_{\Po} \qiso \Omega_{\Po}$ between cofibrant right $\Po$-modules gives a quasi-isomorphism $A \otimes^{\Po} \C(A) \cong L_{\Po} \circ_{\Po}A \qiso \Omega_{\Po} \circ_{\Po}A \cong \Omega_{\Po}(A)$ (since $A$ is cofibrant). Therefore, we have the following sufficient condition for the André-Quillen cohomology to be an Ext-functor.

\begin{thm}[Theorem $17.3.4$ in \cite{Fresse3}]
If $\Po[1]$ and $\Omega_{\Po}$ form cofibrant right $\Po$-modules, then we have
$$\emph{H}_{\Po}^{\bullet}(A,\, M) \cong \emph{Ext}_{A\otimes^{\Po}\K}^{\bullet}(\Omega_{\Po}(A),\, M).$$
\end{thm}

\section{Is André-Quillen cohomology an Ext-functor ?}

In the previous section, we showed that when the André-Quillen cohomology is an Ext-functor, the module of obstructions $\Kc$ is acyclic. It follows, by contraposition, that when the module of obstruction $\Kc$ is not acyclic, then the André-Quillen cohomology is not an Ext-functor. 
Moreover, when $\Po$ is a PBW-operad, the module of obstructions $\Kc$ is acyclic if and only if the André-Quillen cohomology is an Ext-functor. 

In this section, we apply these criteria to the operads $\Com$, $\Perm$ and to the minimal models of Koszul operads. In the case of the operads $\Com$ and $\Perm$, we provide universal obstructions for the André-Quillen cohomology to be an Ext-functor. In the case of an operad which is the cobar construction on a cooperad, we show that the obstructions always vanish. We apply this to the case of homotopy algebras.

\subsection{The case of commutative algebras}

We exhibit a non-trivial element in the homology of the module of obstructions. This gives a universal obstruction for the André-Quillen cohomology of commutative algebras to be an Ext-functor over the enveloping algebra $A\otimes^{\Com}\K$.

\begin{prop}\label{cascom}
The module of obstructions $\mathbb{O}_{\Com}$ is not acyclic. More precisely, we have $$\textrm{\emph{H}}_{1}(\mathbb{O}_{\Com}) \neq 0.$$
\end{prop}

\begin{pf}
Consider the element $\nu := \def\objectstyle{\scriptstyle}
\def\labelstyle{\scriptstyle}
\vcenter{\xymatrix@M=0pt@R=4pt@C=2pt{
1 && 2\\
\ar@{-}[dr] && \ar@{-}[dl]\\
& \ar@{-}[d] &\\
&&}}$ in $\Com^{\ash} \mono \B(\Com)$ and $\mu := \def\objectstyle{\scriptstyle}
\def\labelstyle{\scriptstyle}
\vcenter{\xymatrix@M=0pt@R=4pt@C=2pt{
1 && 2\\
\ar@{-}[dr] && \ar@{-}[dl]\\
& \ar@{-}[d] &\\
&&}}$ in $\Com$. The element $\mu \otimes (\nu \otimes id) = \def\objectstyle{\scriptstyle}
\def\labelstyle{\scriptstyle}
\vcenter{\xymatrix@M=0pt@R=4pt@C=4pt{
1 && 2 &&&\\
\ar@{-}[dr] && \ar@{-}[dl] & 3\\
\ar@{.}[rrrr] & \ar@{-}[rd] && \ar@{-}[ld] &\\
&& \ar@{-}[d] &\\
&&&&&}}$ lives in $O_{\Com}$. We compute
$$d_{O_{\Com}}(\mu \otimes (\nu \otimes id)) = \def\objectstyle{\scriptstyle}
\def\labelstyle{\scriptstyle}
\vcenter{\xymatrix@M=0pt@R=4pt@C=4pt{
& 1 &&&\\
& \ar@{-}[d] & 2 & 3\\
\ar@{.}[rrrr] & \ar@{-}[rd] & \ar@{-}[d] & \ar@{-}[ld] &&\\
&& \ar@{-}[d] &&\\
&&&&&}} + \def\objectstyle{\scriptstyle}
\def\labelstyle{\scriptstyle}
\vcenter{\xymatrix@M=0pt@R=4pt@C=4pt{
&& 2 &&\\
& 1 & \ar@{-}[d] & 3\\
\ar@{.}[rrrr] & \ar@{-}[rd] & \ar@{-}[d] & \ar@{-}[ld] &&\\
&& \ar@{-}[d] &&\\
&&&&&}} = 0.$$
Then $\mu \otimes (\nu \otimes id)$ is a cycle in $O_{\Com}$. However,
$$d_{O_{\Com}}\left(\mbox{\makebox[0pt][l]{
$\def\objectstyle{\scriptstyle}
\def\labelstyle{\scriptstyle}
\vcenter{\xymatrix@M=0pt@R=8pt@C=12pt{
&&&&&\\
&&&&&\\
\ar@{.}[rrrrr] && \ar@{-}[d] &&&\\
&&&&&}}$}
\hspace{-0.1cm}\raisebox{12pt}{$\def\objectstyle{\scriptstyle}
\def\labelstyle{\scriptstyle}
\left(\vcenter{\xymatrix@M=0pt@R=4pt@C=2pt{
1 && 2 && 3\\
\ar@{-}[ddrr] && \ar@{-}[dl] && \ar@{-}[ddll]\\
&&&&\\
&& \ar@{-}[d] &&\\
&&&&}} - \vcenter{\xymatrix@M=0pt@R=4pt@C=2pt{
1 && 2 && 3\\
\ar@{-}[ddrr] && \ar@{-}[dr] && \ar@{-}[ddll]\\
&&&&\\
&& \ar@{-}[d] &&\\
&&&&}}\right)$}}\right) = \def\objectstyle{\scriptstyle}
\def\labelstyle{\scriptstyle}
\vcenter{\xymatrix@M=0pt@R=4pt@C=4pt{
1 && 2 &&&\\
\ar@{-}[dr] && \ar@{-}[dl] & 3\\
\ar@{.}[rrrr] & \ar@{-}[rd] && \ar@{-}[ld] &\\
&& \ar@{-}[d] &\\
&&&&&}} - \def\objectstyle{\scriptstyle}
\def\labelstyle{\scriptstyle}
\vcenter{\xymatrix@M=0pt@R=4pt@C=4pt{
&& 2 && 3 &\\
& 1 & \ar@{-}[dr] && \ar@{-}[dl]\\
\ar@{.}[rrrr] & \ar@{-}[rd] && \ar@{-}[ld] &&\\
&& \ar@{-}[d] &&\\
&&&&&}}$$
and it is impossible to obtain $\mu \otimes (\nu \otimes id)$ as a boundary of an element in $O_{\Com}$. Therefore, this shows that H$_{1}(O_{\Com}) \neq 0$.
$\cqfd$
\end{pf}

\begin{rem}
The short exact sequence $O_{\Com} \mono \mathbb{L}_{\Com} \epi \Omega_{\Com}$ gives a long exact sequence in homology and, since H$_n(\Omega_{\Com}) = 0$ for all $n \geq 1$, we get also H$_{1}(\mathbb{L}_{\Com}) \neq 0$. It follows that there exists a commutative algebra such that the cotangent complex is not acyclic.
\end{rem}

Thanks to Theorem \ref{obstruction}, by contraposition, this gives a conceptual explanation to the fact that the André-Quillen cohomology of commutative algebras cannot always be written as an Ext-functor over the enveloping algebra $A\otimes^{\Com}\K$.

\subsection{The case of Perm-algebras}

The same argument applied to Perm algebras gives a conceptual explanation to the fact that the André-Quillen cohomology of Perm algebras cannot always be written as an Ext-functor over the enveloping algebra $A\otimes^{\Perm}\K$.

\begin{prop}
We have $$\textrm{\emph{H}}_{1}(\mathbb{O}_{\Perm}) \neq 0.$$
\end{prop}

\begin{pf}
The proof is similar to the proof of Proposition \ref{cascom}.
$\cqfd$
\end{pf}

\subsection{The case of algebras up to homotopy}

We show a new homotopy property for algebras over certain cofibrant operads. We apply this in the case of $\Po$-algebras up to homotopy to prove that the André-Quillen cohomology is always an Ext-functor over the enveloping algebra $A\otimes^{\Po_{\infty}}\K$.

\begin{thm}\label{exthomotopy}
Let $\C$ be a coaugmented weight graded dg cooperad and $\Po = \Omega(\C)$ the cobar construction on it. The André-Quillen cohomology of $\Omega(\C)$-algebras is an Ext-functor over the enveloping algebra $A\otimes^{\Omega(\C)}\K$. Explicitly, for any $\Omega(\C)$-algebra $A$ and any $A$-module $M$, we have
$$\mathrm{H}^{\bullet}_{\Omega(\C)}(A,\, M) \cong \emph{Ext}_{A\otimes^{\Omega(\C)}\K}^{\bullet}(\Omega_{\Omega(\C)}(A),\, M).$$
\end{thm}

\begin{pf}
As in this case of $A_{\infty}$-algebras, the twisting morphism $\alpha$ is the map $\C \rightarrow \Omega(\C)$ given in the examples after Theorem \ref{GJres} and the twisting morphism on the level of (co)algebras $\varphi$ is the projection $\C(A) \epi A$. As a dg $\So$-module, the module of obstructions has the following form
$$O_{\Omega(\C)} \cong Rel \oplus \bigoplus_{n\geq 0} \underbrace{(s^{-1}\overline{\C}) \circ_{(1)} \left((s^{-1}\overline{\C}) \circ_{(1)} \cdots \left((s^{-1}\overline{\C}) \right. \right.}_{n \textrm{ times}} \left. \left. \circ_{(1)} \overline{\C} \right) \cdots \right),$$
where $Rel \subset  \bigoplus_{n\geq 0} \underbrace{(s^{-1}\overline{\C}) \circ_{(1)} \left((s^{-1}\overline{\C}) \circ_{(1)} \cdots \left((s^{-1}\overline{\C}) \right. \right.}_{n \textrm{ times}} \left. \left. \circ_{(1)} I \right) \cdots \right)$ is defined in Section \ref{moduleobstruction}.  For any $1 \leq j \leq n$, let $\mu_{j,\, i_{j}}^{c} \in \overline{\C} \circ_{(1)} I$, where $i_{j}$ is the emphasized entry and $\mu_{j}^{c} \in \overline{\C}(m_{j})$. Let $\nu^{c} \in \overline{\C}(m)$. For $\sigma \in \So_{m_{1} + \cdots + m_{n} + m - n}$, we define the map $h$ by
$$\begin{array}{l}
h(s^{-1}\mu_{1,\, i_{1}}^{c} \otimes s^{-1}\mu_{2,\, i_{2}}^{c} \otimes \cdots \otimes s^{-1}\mu_{n,\, i_{n}}^{c} \otimes \nu^{c} \otimes \sigma) = 0 \textrm{ and on } Rel,\\
h\left(\sum_{i_{n} = 1}^{m_{n}} s^{-1}\mu_{1,\, i_{1}}^{c} \otimes \cdots \otimes s^{-1}\mu_{n,\, i_{n}}^{c} \otimes 1^{c} \otimes \sigma \right) = \varepsilon_{n-1} s^{-1}\mu_{1,\, i_{1}}^{c} \otimes \cdots \otimes s^{-1}\mu_{n-1,\, i_{n-1}}^{c} \otimes \mu_{n}^{c} \otimes \sigma,
\end{array}$$
where $\varepsilon_{n-1} = (-1)^{n-1+ |\mu_{1,\, i_{1}}^{c}| + \cdots + |\mu_{n-1,\, i_{n-1}}^{c}|}$. We compute\\
$(dh + hd)(s^{-1}\mu_{1,\, i_{1}}^{c} \otimes \cdots \otimes s^{-1}\mu_{n,\, i_{n}}^{c} \otimes \nu^{c} \otimes \sigma) =$
$$\begin{array}{ll}
= & 0 + h\left(\sum \cdots \otimes \nu'^{c} \otimes \sigma + \varepsilon_{n}\sum_{i = 1}^{m} s^{-1}\mu_{1,\, i_{1}}^{c} \otimes \cdots \otimes s^{-1}\mu_{n,\, i_{n}}^{c} \otimes s^{-1}\nu^{c}_{i} \otimes 1^{c} \otimes \sigma\right) \\
= & s^{-1}\mu_{1,\, i_{1}}^{c} \otimes \cdots \otimes s^{-1}\mu_{n,\, i_{n}}^{c} \otimes \nu^{c} \otimes \sigma,
\end{array}$$
where $\nu'^{c} \in \overline{\C}$ and $i$ is the emphasized entry of $\nu_{i}^{c}$, and\\
$(dh + hd)\left( \sum_{i_{n} = 1}^{m_{n}} s^{-1}\mu_{1,\, i_{1}}^{c} \otimes \cdots \otimes s^{-1}\mu_{n,\, i_{n}}^{c} \otimes 1^{c} \otimes \sigma \right) =$
$$\begin{array}{ll}
= & d(\varepsilon_{n-1} s^{-1}\mu_{1,\, i_{1}}^{c} \otimes \cdots \otimes s^{-1}\mu_{n-1,\, i_{n-1}}^{c} \otimes \mu_{n}^{c} \otimes \sigma)\\
+ & h\left(\sum_{i_{n} = 1}^{m_{n}} \sum_{j=1}^{n} \sum \varepsilon_{j-1}(-1)^{|{\mu'}^{c}_{j,\, i'_{j}}|} s^{-1}\mu_{1,\, i_{1}}^{c} \otimes \cdots \otimes s^{-1}{\mu'}^{c}_{j,\, i'_{j}} \otimes s^{-1}{\mu''}^{c}_{j,\, i''_{j}} \otimes \cdots \otimes s^{-1}\mu_{n,\, i_{n}}^{c} \otimes 1^{c} \otimes \sigma \right)\\
+ & h\left(\sum_{i_{n} = 1}^{m_{n}} \sum_{j=1}^{n} \varepsilon_{j-1} s^{-1}\mu_{1,\, i_{1}}^{c} \otimes \cdots \otimes (-s^{-1}d_{\C}(\mu_{j,\, i_{j}})) \otimes \cdots \otimes s^{-1}\mu_{n,\, i_{n}}^{c} \otimes 1^{c} \otimes \sigma\right)\\
= & \varepsilon_{n-1}\sum_{j=1}^{n-1} \varepsilon_{j-1}(-1)^{|{\mu'}^{c}_{j,\, i'_{j}}|} s^{-1}\mu_{1,\, i_{1}}^{c} \otimes \cdots \otimes s^{-1}{\mu'}^{c}_{j,\, i'_{j}} \otimes s^{-1}{\mu''}^{c}_{j,\, i''_{j}} \otimes \cdots \otimes s^{-1}\mu_{n-1,\, i_{n-1}}^{c} \otimes \mu_{n}^{c} \otimes \sigma\\
+ & \varepsilon_{n-1}\sum \varepsilon_{n-1} s^{-1}\mu_{1,\, i_{1}}^{c} \otimes \cdots \otimes s^{-1}\mu_{n-1,\, i_{n-1}}^{c} \otimes s^{-1}{\mu'}^{c}_{n,\, i'_{n}} \otimes {\mu''}^{c}_{n} \otimes \sigma\\
+ & \varepsilon_{n-1}\varepsilon_{n-1} \sum_{i_{n} = 1}^{m_{n}} s^{-1}\mu_{1,\, i_{1}}^{c} \otimes \cdots \otimes s^{-1}\mu_{n,\, i_{n}}^{c} \otimes 1^{c} \otimes \sigma\\
- & \varepsilon_{n-1}\sum_{j=1}^{n-1} \varepsilon_{j-1} s^{-1}\mu_{1,\, i_{1}}^{c} \otimes \cdots \otimes s^{-1}d_{\C}(\mu_{j,\, i_{j}}) \otimes \cdots \otimes s^{-1}\mu_{n-1,\, i_{n-1}}^{c} \otimes \mu_{n}^{c} \otimes \sigma\\
+ & \varepsilon_{n-1} \varepsilon_{n-1} s^{-1}\mu_{1,\, i_{1}}^{c} \otimes \cdots \otimes s^{-1}\mu_{n-1,\, i_{n-1}}^{c} \otimes d_{\C}(\mu_{n}^{c}) \otimes \sigma\\
&\\
- & \varepsilon_{n-1}\sum_{j=1}^{n-1} \varepsilon_{j-1}(-1)^{|{\mu'}^{c}_{j,\, i'_{j}}|} s^{-1}\mu_{1,\, i_{1}}^{c} \otimes \cdots \otimes s^{-1}{\mu'}^{c}_{j,\, i'_{j}} \otimes s^{-1}{\mu''}^{c}_{j,\, i''_{j}} \otimes \cdots \otimes s^{-1}\mu_{n-1,\, i_{n-1}}^{c} \otimes \mu_{n}^{c} \otimes \sigma\\
- & \varepsilon_{n-1}(-1)^{|{\mu'}^{c}_{n,\, i'_{n}}|}\sum \varepsilon_{n-1}(-1)^{|{\mu'}^{c}_{n,\, i'_{n}}|} s^{-1}\mu_{1,\, i_{1}}^{c} \otimes \cdots \otimes s^{-1}\mu_{n-1,\, i_{n-1}}^{c} \otimes s^{-1}{\mu'}^{c}_{n,\, i'_{n}} \otimes {\mu''}^{c}_{n} \otimes \sigma\\
+ & \varepsilon_{n-1}\sum_{j=1}^{n-1} \varepsilon_{j-1} s^{-1}\mu_{1,\, i_{1}}^{c} \otimes \cdots \otimes s^{-1}d_{\C}(\mu_{j,\, i_{j}}) \otimes \cdots \otimes s^{-1}\mu_{n-1,\, i_{n-1}}^{c} \otimes \mu_{n}^{c} \otimes \sigma\\
- & \varepsilon_{n-1} \varepsilon_{n-1} s^{-1}\mu_{1,\, i_{1}}^{c} \otimes \cdots \otimes s^{-1}\mu_{n-1,\, i_{n-1}}^{c} \otimes d_{\C}(\mu_{n}^{c}) \otimes \sigma\\
&\\
= & \sum_{i_{n} = 1}^{m_{n}} s^{-1}\mu_{1,\, i_{1}}^{c} \otimes \cdots \otimes s^{-1}\mu_{n,\, i_{n}}^{c} \otimes 1^{c} \otimes \sigma.
\end{array}$$
Thus $dh + hd = id$ and $h$ is a homotopy. Finally, $O_{\Omega(\C)}$ is acyclic and Theorem \ref{obstruction} gives the theorem since any quasi-free operad $\Po = \Omega(\C)$ is a PBW-operad (because a free operad has no relations).
$\cqfd$
\end{pf}

We conjecture that this theorem is true for any cofibrant operad.

When $\Po$ is a Koszul operad, the previous theorem applied to $\C = \Poa$ shows that the André-Quillen cohomology of a homotopy algebra is always an Ext-functor over its enveloping algebra.

Let $A$ be a $\Po$-algebra. The algebra $A$ is a $\Po_{\infty}$-algebra since there is a map of operads $\Po_{\infty} = \Omega(\Poa) \rightarrow \Po$. Similarly, an $A$-module over the operad $\Po$ is also an $A$-module over the operad $\Po_{\infty}$. This leads to the following result.

\begin{prop}\label{homotopyounon}
Let $\Po$ be a Koszul operad and let $A$ be a $\Po$-algebra. The André-Quillen cohomology of the $\Po$-algebra $A$ is equal to the André-Quillen cohomology of the $\Po_{\infty}$-algebra $A$. That is,
$$\textrm{\emph{H}}^{\bullet}_{\Po}(A,\, M) = \textrm{\emph{H}}^{\bullet}_{\Po_{\infty}}(A,\, M), \textrm{ for any $A$-module $M$ over the operad $\Po$.}$$
\end{prop}

\begin{pf}
A resolution of $A$ as a $\Po$-algebra is given by $\Po \circ \Poa (A)$ and a resolution of $A$ as a $\Po_{\infty}$-algebra is given by $\Po_{\infty} \circ \Poa(A)$. Thus, by Theorem \ref{iso2}, we have
$$\textrm{Hom}_{\M_A^{\Po_{\infty}}}(A\otimes^{\Po_{\infty}} \Poa(A),\, M) = \textrm{Hom}_{g\mathcal{M}od_{\K}}(\Poa(A),\, M) = \textrm{Hom}_{\M_A^{\Po}}(A\otimes^{\Po} \Poa(A),\, M).$$
Moreover, the differential on Hom$_{g\mathcal{M}od_{\K}}(\Poa(A),\, M)$ is the same in both cases since the higher products $\Poa(k)\otimes_{\So_k} A^{\otimes k} \rightarrow A$ for $k\geq 3$ are $0$.
$\cqfd$
\end{pf}

We showed that, for commutative algebras and Perm algebras, the André-Quillen cohomology of a $\Po$-algebra cannot always be written as an Ext-functor over the enveloping algebra $A\otimes^{\Po}\K$. However, by the following theorem, it can always be written as an Ext-functor over the enveloping algebra $A\otimes^{\Po_{\infty}}\K$.

\begin{thm}
Let $\Po$ be a Koszul operad, let $A$ be a $\Po$-algebra and let $M$ be an $A$-module over the operad $\Po$. We have
$$\textrm{\emph{H}}^{\bullet}_{\Po}(A,\, M) \cong \emph{Ext}_{A\otimes^{\Po_{\infty}}\K}^{\bullet}(\Omega_{\Po_{\infty}}(A),\, M).$$
\end{thm}

\begin{pf}
We make use of Theorem \ref{exthomotopy} and Proposition \ref{homotopyounon}.
$\cqfd$
\end{pf}

\section*{Acknowledgments}

I would like to thank my advisor Bruno Vallette for his useful and constant help. I am very grateful to Benoit Fresse for ideas and useful and clarifying discussions and to Martin Markl for discussions about his cohomology theory for homotopy associative algebras and for many relevant remarks on the first version of the paper. I am grateful to Henrik Strohmayer for his careful reading of my paper and for all his corrections and to Vladimir Dotsenko for pointing a mistake in a previous version of Theorem \ref{acyclicity trivial}. I also wish to thank the referee for his numerous comments.

\bibliographystyle{alpha}
\bibliography{bib}
\end{document}